\documentclass[10pt,reqno]{amsart}
\usepackage{epigraph}
\usepackage{fancyhdr}
\usepackage{amssymb}
\usepackage{amsmath}
\usepackage{tikz-cd}
\usepackage{array}

\input {xy}

\xyoption{all}
\usepackage[all]{xypic}
\usepackage{latexsym}
\usepackage{setspace}
\usepackage{enumitem}
\usepackage{tabularx}
\usepackage{booktabs}
\usepackage[labelfont=bf,format=plain,justification=raggedright,singlelinecheck=false]{caption} 
\usepackage[pagewise]{lineno}
\usepackage[widespace]{fourier}
\usepackage{frcursive}
\usepackage{xcolor}
\usepackage{hyperref}
\hypersetup{colorlinks=true,linkcolor=blue, citecolor=blue}

\begin{document}

\input{amssym.def}

\newsymbol \circledarrowleft 1309

 %\if@twoside \oddsidemargin 6pt \evensidemargin 6pt \marginparwidth 90pt
 %\else \oddsidemargin 18pt \evensidemargin 18pt \marginparwidth 68pt 
 %\fi
 %\marginparsep 10pt \topmargin -30pt \headheight 12pt \headsep 25pt 
 %\footheight 12pt \footskip 30pt \textheight 680pt \textwidth 455pt 
 %\columnsep 10.5pt \columnseprule 0pt

\newcommand{\ctext}[1]{\makebox(0,0){#1}}
\setlength{\unitlength}{0.1mm}

\newcommand{\wt}{\widetilde}
\newcommand{\Lr}{\Longrightarrow}
\newcommand{\Aut}{\mbox{{\rm Aut}$\,$}}
\newcommand{\ul}{\underline}
\newcommand{\ol}{\overline}
\newcommand{\lr}{\longrightarrow}
\newcommand{\bc}{{\mathbb C}}
\newcommand{\bp}{{\mathbb P}}
\newcommand{\bz}{{\mathbb Z}}

\newcommand{\bwt}{\boldsymbol\Omega}
\newcommand{\bbf}{{\mathbb M}}
\newcommand{\be}{{\mathbb E}}
\newcommand{\ba}{{\mathbb A}}
\newcommand{\bh}{{\tiny\tc{H}}}
\newcommand{\ca}{{\mathcal A}}
\newcommand{\ch}{{\mathcal H}}
\newcommand{\cu}{{\mathcal U}}
\newcommand{\cf}{{\mathcal F}}
\newcommand{\ce}{{\mathcal E}}
\newcommand{\co}{{\mathcal O}}
\newcommand{\cg}{{\mathcal G}}
\newcommand{\cm}{{\mathcal M}}
\newcommand{\fm}{{\mathfrak M}}
\newcommand{\fw}{{\mathfrak W}}
\newcommand{\fb}{{\mathfrak B}}
\newcommand{\fg}{{\mathfrak G}}
\newcommand{\fa}{{\mathfrak a}}
\newcommand{\fA}{{\mathfrak A}}
\newcommand{\fd}{{\mathfrak d}}
\newcommand{\cF}{\mathcal{F}}
\newcommand{\cz}{\mathcal{Z}}
\newcommand{\mcp}{{\mathcal P}}
\newcommand{\mcw}{{\mathcal W}}
\newcommand{\cp}{{\sf P}}
\newcommand{\tq}{{\tt q}}
\newcommand{\ta}{{\tt A}}
\newcommand{\gG}{\mathfrak{G}}
\newcommand{\cT}{\mathcal{T}}
\newcommand{\eF}{{\mathscr F}}

\newcommand{\cq}{{\sf Q}}
\newcommand{\cs}{{\mathcal S}}
\newcommand{\cl}{{\mathcal L}}
\newcommand{\mgl}{{\mathfrak {GL}}}

\newcommand{\hra}{\hookrightarrow}
\renewcommand\qedsymbol{\tt{Q.E.D}}

\newtheorem{guess}{\sc Theorem}[section]
\newcommand{\bth}{\begin{guess}$\!\!\!${\bf}~~\it}
\newcommand{\eeth}{\end{guess}}

\newtheorem{propo}[guess]{\sc Proposition}%[section]
\newcommand{\bprop}{\begin{propo}$\!\!\!${\bf }~~\it}
\newcommand{\eprop}{\end{propo}}
\newtheorem{rema}[guess]{\it Remark}%[section]
\newcommand{\brem}{\begin{rema}$\!\!\!${\bf }~~\it}
\newcommand{\erem}{\end{rema}}
\newtheorem{coro}[guess]{\sc Corollary}%[section]

\newcommand{\bcor}{\begin{coro}$\!\!\!${\bf }~~\it}
\newcommand{\ecor}{\end{coro}}
\newtheorem{lema}[guess]{\sc Lemma}%[section]

\newcommand{\blem}{\begin{lema}$\!\!\!${\bf }~~\it}
\newcommand{\elem}{\end{lema}}

\newtheorem{exam}[guess]{\it Example}%[section]
\newcommand{\beg}{\begin{exam}$\!\!\!${\bf }~~\rm}
\newcommand{\eeg}{\end{exam}}
\newcommand{\er}{\hfill {\Large $\bullet$}\linebreak}

\newtheorem{defe}[guess]{\sc Definition}
\newcommand{\bdefe}{\begin{defe}$\!\!\!${\bf }~~\rm}
\newcommand{\edefe}{\end{defe}}
\newtheorem{nota}[guess]{\sc Notation}

\newcommand{\bpr}{\begin{proof}}
\newcommand{\epr}{\end{proof}}
\numberwithin{equation}{guess}

%\numberwithin{equation}{subsection} 
\newcommand{\beqa}{\begin{eqnarray}}
\newcommand{\eeqa}{\end{eqnarray}}
\newcommand{\Hilb}{\operatorname{Hilb}}
\newcommand{\Gal}{\on{Gal}}
\newcommand{\spec}{{\rm Spec}\,}
\newcommand{\hspin}{\operatorname{HSpin}}
\newcommand{\spin}{\operatorname{Spin}}
\newcommand{\sym}{\operatorname{Sym}}
\newcommand{\lie}{\operatorname{Lie}}
\newcommand{\Cl}{\operatorname{C\ell}}
\newcommand{\GL}{\operatorname{GL}}
\newcommand{\PGL}{\operatorname{PGL}}
\newcommand{\SL}{\operatorname{SL}}
\newcommand{\hol}{\operatorname{hol}}
\newcommand{\End}{\operatorname{End}}
\newcommand{\Hom}{\operatorname{Hom}}
\newcommand{\tc}[1]{\text{\cursive{#1}}}
\newcommand{\tsf}[1]{\textsf{#1}}
\newcommand{\msf}{M^s_{\tsf{f}}}
\newcommand{\fr}[1]{\mathfrak{#1}}
\newcommand{\on}[1]{\operatorname{#1}}
%\newcommand{\PGL}{\operatorname{PGL}}
%\newcommand{\mathscr}[1]{\mathscr{#1}}

%\pagestyle{myheadings}

%\notitlepage

\title{On the desingularisation of moduli of principal bundles}
\author{Vikraman Balaji} \address {Chennai Mathematical Institute
  SIPCOT IT Park, Siruseri, Chennai-603103 INDIA}
\email{balaji@cmi.ac.in, vikramanbalaji@gmail.com}

%\date{Preliminary version}

\begin{abstract} In \cite{nr} Narasimhan and Ramanan and in \cite{desing}, Seshadri constructed desingularisations of the moduli space $M^{ss}_{_{\text{SL}(2)}}$ of semistable $\SL(2)$-bundles on a smooth projective curve $C$ of genus $g \geq 3$. Seshadri's construction was even modular and canonical. In this paper, we construct a smooth modular compactification of the moduli of stable principal $H$-bundles when $H$ is a simply connected almost simple algebraic group of type ${\tt B}_{_\ell}, {\tt D}_{_\ell}, {\tt G}_{_2}, {\tt F}_{_4}~~or~~{\tt C}_{_3}$. These spaces give canonical desingularisations of the moduli space $M^{ss}_{_H}$ of semistable principal $H$-bundles and thereby, a comprehensive generalisation of \cite{desing}.
  
\end{abstract}
%\vspace{-2mm}

%{\dn mA\2 a\7{n}-mr \7{y}@y c},
\maketitle

\epigraph{{\em{M$\overline{a}$\d{m} anusmara, yuddhya \c{c}a}},\\{\tt {\footnotesize Hold Me in mind and battle on}}.\\ BG, Chapter VIII, verse 7}

\small
{\tt\tableofcontents}
\setcounter{tocdepth}{2} % Set depth to 1
\normalsize

%\vspace{2cm}
%\noindent

%\vspace{2mm}
\section{Introduction} Let $k$ be an algebraically closed field of characteristic $0$ or large prime characteristic $p$, with bounds on $p$ in terms of $H$ (see \S\ref{charp}). Let $C$ be a smooth projective curve over $k$ of genus $g \geq 2$. Let $H$ be an almost simple, simply connected affine algebraic group over $k$, and $M^{ss}_{_H}$ denote the moduli space of slope semistable principal $H$-bundles over $C$. By the work of A.Ramanathan \cite{r1}, the moduli space $M^{ss}_{_H}$ is an irreducible normal projective variety over $k$. It is {\it invariably} singular unless $H = \text{SL}(2)$ and $g = 2$; in this very special case, it is even isomorphic to the projective space $\bp^{^3}$. 

Around 1978, Narasimhan-Ramanan \cite{nr} and Seshadri \cite{desing}, independently constructed desingularisations of $M^{ss}_{_{\SL(2)}}$.   An added feature in Seshadri's space  was that it solved a natural moduli problem. Moreover, Seshadri constructed birational models for $M^{ss}_{_{\SL(n)}}$ for all $n \geq 2$. These new models however remain singular for $n \geq 3$, as was shown in Le Bruyn and Reichstein \cite{zino}. Subsequently, F.Kirwan \cite{kirwan} constructed {\em partial desingularisations} of $M^{ss}_{_{\SL(n)}}$ for all $n$. These have at most finite quotient singularities but do not parametrize solutions of a  natural moduli problem; however, these have proven to be extremely useful in understanding the topology of the base moduli spaces.  In the paper by Kiem and Li \cite{kiemli}, the two smooth models, \cite{nr} and \cite{desing}, for $M^{ss}_{_{\SL(2)}}$, are compared and related to the one in \cite{kirwan}. We note that, for the moduli spaces of principal bundles for a general structure group, there is no construction analogous to the one given by Kirwan for the linear case.

The question of constructing  a  canonical desingularisation of $M^{ss}_{_H}$, which is also {\em modular} offers serious challenges and has remained open. By {\em modular and canonical}, we mean that the objects of the new space parametrize solutions of natural moduli problem upto normalizations. One also wonders as to why only the  case of $\SL(2)$ could be handled by the approaches in \cite{nr} and \cite{desing}. The aim of this paper is to address these issues and attempt to give some satisfactory answers. 

Let $H$ be a group of type ${\tt B}_{_\ell}, \ell \geq 1, {\tt D}_{_\ell}, \ell \geq 2, {\tt G}_{_2}, {\tt F}_{_4}, {\tt C}_{_3}$. We construct   smooth compactifications of the moduli space of stable principal $H$-bundles which give canonical desingularisations of $M^{ss}_{_H}$.  This also resolves the question as to why only $\SL(2) and \GL(2)$ appear in \cite{nr} and \cite{desing}. Indeed, by the low dimensional isomorphism $\SL(2) \simeq  \spin(3)$, the group $\SL(2)$ is now the first among our list of groups. The new list  now raises questions as to why only these  groups. We try to give a reasoning in terms of some natural spinorial plethysms related to the exceptional groups in our list (see \S8). It will be interesting to have something more satisfactory.

%A weak reason is that apart from the classical Spin groups, the exceptional groups which occur in our list come from a pretty spinor "plethysm", which we realized later has been discovered by Dynkin . 

 The classical compactifications of the moduli of stable princiapl $H$-bundles ({\`a la} Mumford-Seshadri-Ramanthan), is via GIT. In this  approach the boundary of the stable bundles in the compactification continue to remain principal $H$-bundles, i.e., they stay within the category of principal $H$-bundles, but the coarse moduli construction attaches $S$-equivalence classes of semistable $H$-bundles so as to compactify the moduli space of stable bundles. This makes these spaces singular. More precisely, the points in the boundary of the smooth  locus are natural equivalence classes, nevertheless the representatives remain principal $H$-bundles.  As we noted above, by the classical approach, barring one exception (namely, $\SL(2)$-bundles on curves of $g = 2$), one invariably obtains singular moduli spaces. 

The  paradigm shift  in the present approach is  (1) to allow the structure group $H$, of the stable $H$-bundles, to degenerate in  certain controlled  "flat families", and (2) to then realise the boundary  as isomorphism classes of torsors for these degenerate group schemes. Precisely speaking, we construct {\em flat} group schemes $\bh$ on complete regular local rings with generic fibre $H$, with natural versal properties,  but which degenerate to {\em non-reductive} limits $H'$, much in the sense of  a higher dimensional analogue of Bruhat-Tits group schemes (see \cite{bp} and \cite{torsors}). The boundary points of the new moduli space are natural equivalence classes under some relation, but the representative points in the moduli space are now isomorphism classes of principal $H'$-bundles. {\em In this sense, the present paper is therefore an application of the theory of Bruhat-Tits group schemes and torsors on higher dimensional base schemes as developed in \cite{bp}}

The moduli problem then becomes a classification of isomorphism classes of Ramanathan-stable $H$-torsors which degenerate to $H'$-torsors on the curve $C$. The limiting torsors  are also additionally {\em stable} in a certain natural sense. One needs to construct versal deformation spaces for this problem which are shown to be smooth, and this is somewhat delicate. This results in the smoothness of the moduli of torsors for the new group scheme. Proving the properness of the resulting moduli space is almost immediate and is  reduced to the compactness of the moduli spaces  in \cite{desing}. 

In \cite[Corollary 3.4]{bh} it is shown that the smooth locus of $M^{ss}_{_H}$ is precisely the locus $M^{{rs}}_{_H}$ of {\em regularly stable} bundles. I prove the following (see \eqref{main}, \eqref{maineven}, and \eqref{mainexceptional}).

\bth\label{mainintro} {\sf {\em Let $H$ be of type ${\tt B}_{_\ell}, \ell \geq 1, {\tt D}_{_\ell}, \ell \geq 2, {\tt G}_{_2}, {\tt F}_{_4}, ~or~{\tt C}_{_3}$. There is a canonical smooth compactification ${\fm}_{_H}$ of the moduli space $M^{rs}_{_H}$ of regularly stable principal $H$-bundles. There is natural birational morphism $p:{\fm}_{_H} \to M^{ss}_{_H}$, which is an isomorphism precisely over $M^{rs}_{_H}$}}. Further, ${\fm}_{_H}$ can be realized as the normalization of a space which represents a natural moduli functor.
\eeth
As an application, we  compute low degree Betti numbers of the new compactification ${\fm}_{_H}$ for the case when $H =  {\tt G}_{_2}$. This illustrates how one could get to know the structure of the smooth compactification,  which are described in terms of torsors under certain Bruhat-Tits group schemes.  A fuller comprehension of this structure will lead to the  fuller computation of the rational cohomology of these spaces which will be for the future.

In the {\underline{Appendix}}, we address the problem of desingularising the ``local model" of the singularity of the moduli space $M^{ss}_{_H}$ at its "worst singularity" i.e. the trivial bundle. This is a normal singularity which, by an easy application of the Luna \'etale slice theorem, can locally be described as a GIT quotient of a $g$-fold product of the Lie algebra of $H$ for the diagonal action of the adjoint group $\bar{H}$. As in Nori's appendix (cf. \cite{nori}), one could ask whether the construction of the desingularisation ${\fm}_{_H}$ suggests an analogous construction of a natural desingularisation of the local models. In collaboration with Rajarshi Ghosh, we  prove that  this is indeed the case in all the group types which appear in the main theorem. The higher dimensional analogues of the Lie algebras of the Bruhat-Tits group schemes plays the key role.

%\footnotesize
{\it Acknowledgments}:  I thank Rajarshi Ghosh (a fourth year CMI student) for several serious discussions. His collaboration on the material in the Appendix helped understand the issues which cropped up in the main paper. I thank Rohith Verma for many helpful discussions. I also thank Jagadish Pine for his interest in the work.
%\pagebreak
%\tableofcontents
%\normalsize
%\input{contents.tex}
%\footnotesize
\subsubsection{Terminology}
Throughout this paper, unless otherwise stated, we have the
following notations and assumptions:
%\small
{\renewcommand{\labelenumi}{{\rm (\alph{enumi})}}
\begin{enumerate}

\item We work over an algebraically closed field $k$ of characteristic zero (or of prime characteristic under some conditions \S\ref{charp}.

\item $C$ will be an irreducible, smooth projective curve over $k$.

\item $H$, $\fg$, will be an {\it almost simple}, connected, simply connected algebraic groups.
\item $M_{_{H}}^{ss}$ ($M_{_{H}}^{s}$, $M^{rs}_{_H}$) always stand for the moduli space of semi(stable, {\em regularly stable}) principal $H$-bundles on $C$.

\item Let $\ch$ be a $T$-group scheme. By an $\ch$-bundle or $\ch$-torsor on $C \times T$ we mean a $\ch$-torsor on $C \times T$ for the group scheme $\ch \times_{_T} (C \times T)$.

\item  Let $E$ be a principal $\ch$-bundle on $C \times T$. If $x \in C$ is a closed point then we shall denote by
 $E_{x,T}$ the restriction of $E$ to the
sub-scheme  $x \times T$, which will then be an $\ch$-torsor on $T$.  Similarly, $t \in T$ will denote the closed point of $T$ and the
restriction of $E$ to $C \times t$ will be denoted by $E_t$.

%\item In the case of $G = GL(V)$, when we speak of a principal $G$-bundle we identify it often with the associated vector bundle (and can therefore talk of the degree of the principal $G$-bundle with respect to the choice of $\Theta$).

%\item If $H_A$ is an $A$-group scheme, then by $H_A(A)$ (resp.$H_K(K)$) we mean its $A$ (resp $K$)-valued points. When $H_A = H \times \spec A$, then we simply write $H(A)$ for its $A$-valued points. We denote the closed fibre of the group scheme by $H_k$.

\item Let $Y$ be any $G$-variety and let $E$ be a $G$-principal bundle. For example $Y$ could be a $G$-module. Then we denote by $E(Y)$ the associated bundle with fibre type $Y$ which is the following object: $E(Y) := (E \times Y)/G$ where the equivalence relation is by the twisted action of $G$ on $E \times Y$, given by $\small\text{\cursive g}.(e,y)~=~(e.\small\text{\cursive g},~~\small\text{\cursive g}^{-1}.y)$.

%\item If we have a group scheme $H_A$ (resp $H_K$) over$\spec A$ (resp $\spec K$) an $H_A$-module $Y_A$ and a principal $H_A$-bundle $E_A$. Then we shall denote the associated bundle with fibre type $Y_A$ by $E_A(Y_A)$.

\end{enumerate}
\part{\sc The Construction of the Desingularisation}
\section{Recollections from Seshadri's construction}\label{seshadri} 
\subsubsection{The basic definitions}\label{basic1}

A finite dimensional unital associative algebra $k$-algebra $A_{_0}$ of $\dim A_{_0} = n^{^2}$ is called a {\em specialisation (or degeneration) of  the matrix algebra $\cm_{_n}(k)$}, denoted $A_{_0} \sim \cm_{_n}$, if there exists a discrete valuation ring $R$ with residue field $k$ and fraction field $K$, and an algebra $A$ over $R$ which  is free over $R$ as an $R$-module, such that $A \otimes_{_k} \overline{K}$ is isomorphic to $\cm_{_n}(\overline{K})$ and $A \otimes_{_R} k \simeq A_{_0}$. 

Let $V$ be a vector space over $k$ such that $\dim(V) = n^2$. Let $\text{Alg}(V)$ be the linear subspace of $\Hom(V \otimes V,V)$ whose elements are associative algebra structures on $V$. Let $e_0$ be a non-zero element of $V$ and let $\ca_{_n}$ denote the sub-variety of $\text{Alg}(V)$ consisting of specialisations of $\cm_{_n}$ with the identity element being $e_0$. We have, over $\text{Alg}(V)$, an $\co$-algebra $\fA$ such that for each point $z \in \text{Alg}(V)$, the $k$-algebra $\fA_z \otimes k$ is isomorphic to $z$. Let $\fA$ continue to denote the restriction of $\fA$ to $\ca_{_n}$.

\subsubsection{The versal property}\label{versalproperty}The following simple fact is from (\cite[page 111]{ast}): Let $Y$ be a reduced $k$-scheme, and $B$ an $\co_{_Y}$-algebra which is locally free of rank $n^2$. Further, suppose that there exists an open dense subset $Y' \subset Y$  such that for every point $y \in Y'$, the $k$-algebra $B_{_Y}\otimes k$ is a specialisation of the matrix algebra $\cm_{_n}$. Then for every point $y \in Y$, there exists a neighbourhood $U \subset Y$ and a morphism $f:U \to \ca_{_n}$ such that the $\co_{_U}$-algebras $f^{*}(\fA)$ and $B|_{_U}$ are isomorphic.

\subsubsection{The birational model in \cite{desing}} Let $N_{_n}$ be the moduli space of isomorphism classes of semistable vector bundles $W$ of rank $n^{^2}$ with trivial determinant such that (1) the endomorphism algebra $\End(W)$ is a specialisation of $\cm_{_n}$, and (2) it further  for admits a suitable {\em parabolic structure} which renders it {\em parabolic stable}. In \cite{desing}, Seshadri shows that $N_{_n}$ is a reduced and irreducible projective scheme over $k$, and further, there exists a birational morphsim $\pi_{_n}:N_{_n} \to M_{_{\SL(n)}}$ such that over the open subset 
\beqa\label{stableinN}
N_{_n}^s := \{W \in N_{_n} : \End(W) \sim \cm_{_n} \}
\eeqa
and $\pi_{_n}:N_{_n}^s \to M_{_{\SL(n)}}^{^{s}}$ is an isomorphism. In fact, $W \in N_{_n}^s$ if and only if $W \simeq \oplus W_1$, i.e.,  is isomorphic to a direct sum of $n$-copies of a stable $W_1 \in M_{_{\SL(n)}}^{^{s}}$.

It is shown in \cite{desing} that $N_{_n}$ is formally smooth to the scheme of specialisations of $\cm_{_n}$; this scheme of specialisations is shown in \cite{desing} to be smooth when $n = 2$; it is known to be singular for all $n \geq 3$ \cite{zino}. Thus, $N_{_n}$ is singular for all $n \geq 3$.

\subsubsection{Structure group reduction}\label{strgrp} It is shown in (\cite[Prop 4]{desing}, \cite[Prop 10]{ast}), that if $W$ is a point in $N_{_n}$, then the structure group of the principal $\GL(n^2)$-bundle underlying $W$ can be reduced to the group of units $\End(W)^{\circ,\times}$ of the opposite of the endomorphism algebra $\End^{\circ}(W)$ (see \cite[Propositions 3, 4 and 10]{ast}). Since $W$ is a semistable vector bundle of degree $0$, by choosing a base point $x \in C$, we get an embedding $\End(W) \hra \End(W_x)$. The opposite algebra is then identified with the centralizer $C(\End(W))$ of $\End(W)$ in $\End(W_x)$. This way, we get an inclusion of the group of units  $\End(W)^{\circ,\times}$ in $\GL(n^2)$, where the linear group is identified with $\Aut(W_x)$. This discussion can be carried over to an arbitrary reduced base scheme $T$ instead of $\spec k$. 
\subsubsection{Existence of a universal family}\label{univfamilyinseshadri} By embedding the moduli space $N_{_n}$ in the moduli space of parabolic vector bundles (with small weights), it is shown in \cite{desing} that $N_{_n}$ supports a universal family $\mcw \to C \times N_{_n}$ of rank $n^2$. We will denote the underlying principal $\GL(n^2)$-bundles by $P_{_n}$ and view it as a principal $\Aut(W_x)$-bundle, where $W$ is a point in $N_{_n}$. Seshadri also deduces the projectivity of $N_{_n}$ using this embedding.

%\begin{nota}\label{basepoint} We fix a base point $x \in C$ once for all and identify the structure groups $\GL(n^2)$ of a bundle $W \in N_{_n}$ with $\text{Aut}(W_x)$. \end{nota}

\section{Some remarks on the Clifford Algebra}
Let $\bbf$ be a finite dimensional vector space over $k$ and let $\dim \bbf = 2{\ell} + 1 =: m$ which is assumed to be odd. Set $n := 2^{^\ell}$. 

Let ${\tq}:\bbf \to k$ be a quadratic form and let $b_{_\tq}(x,y)$ be the associated symmetric bilinear form. Let $\Cl_\tq$ denote the Clifford algebra associated to $\tq$. It is well-known that $\Cl_\tq$ is a $\bz/2$-graded unital associative algebra, which is isomorphic to the exterior algebra $\Lambda \bbf$ as a $k$-module and therefore has dimension $\dim(\Cl_\tq) = 2^m = 2^{2{\ell} +1}$. The positive graded part 
\[
\Cl^{+}_\tq \subset \Cl_\tq
\]
has dimension $n^2 = 2^{2{\ell}}$; when $\tq$ is non-degenerate, then $\Cl^{+}_\tq$ is isomorphic to  the matrix algebra $\cm_{_n}$.

Let $\sf Q := \sym^2(\bbf)$, denote the affine space of quadratic forms on $\bbf$ or equivalently, the affine space of symmetric $m \times m$-matrices. 
\bprop\label{forversality} The map $\theta:{\sf Q} \to \ca_{_n}, \tq \mapsto \Cl^{+}_\tq$ is a closed embedding. \eprop
\bpr The proof follows the ideas in \cite[Theorem 1, Remark 3, page 165-167]{desing}, more specifically, some comments due to S.Ramanan in \cite{desing}. Let $e_0 \in \Cl^{+}_\tq$ be the identity element and let us fix a basis $\langle e_1, \ldots, e_n \rangle$ for the vector space $\bbf$. Recall that the multiplication rule for the Clifford algebra  $\Cl_\tq$ is given as follows:
for $x \in \Cl_\tq$, $x^2 = \tq(x).e_0$; which gives $xy +yx = b_\tq(x,y).e_0$, $\forall x,y$. For the basis vectors this gives 
\beqa\label{laws}
e_i^2 = \tq(e_i).e_0\\
e_i.e_j + e_j.e_i = b_\tq(e_i,e_j).e_0
\eeqa
Let $\Cl_\tq^{^{\leq r}} \subset \Cl_\tq$ be the $k$-submodule generated by products $a_1 \cdots a_j$, with $a_{_1}, \ldots, a_{_j} \in \bbf$, where $0 \leq j \leq r$. In particular, we concentrate on the case $r = 2$, i.e., the submodule $\Cl_\tq^{^{\leq 2}}$. Let 
\beqa\label{helmstetter}
L_\tq := \Cl^{+}_\tq \cap \Cl_\tq^{^{\leq 2}}.
\eeqa
Then by \cite[5.5.3, page 263]{helm}, the natural Lie algebra structure on 
$\Cl^{+}_\tq$ arising from the associative algebra structure, restricts to give a Lie subalgebra structure on $L_\tq$. Let $[.]$ denote the induced Lie bracket.

Let $\alpha_{ij} := e_i.e_j$. Then, in terms of the basis vector it can be seen that 
\beqa\label{liebasis}
L_\tq = \langle e_0, \alpha_{12}, \ldots, \alpha_{1n}, \alpha_{23}, \ldots, \alpha_{n-1n}\rangle
\eeqa
as a $k$-vector space. Since $[e_0,e_0] = 0$, the quotient $L'_\tq := L_\tq/k.e_0$ gets a natural Lie algebra structure, where $[\overline{e}_i,\overline{e}_j] = [e_i,e_j]~mod~ke_0$. Let $\overline{\alpha}_{ij} := \overline{e}_i.\overline{e}_j$, which gives a basis for 
$L'_\tq$.

{\em Claim}: The map $\theta$ is injective {\em if $2$ is invertible}: in fact, we claim that the multiplication table for the Lie algebra structure on $L'_\tq$ recovers the bilinear form $b_{\tq}$ and hence the quadratic form $\tq$. Observe that for $i < j < k < l$, we have:
\beqa 
\overline{\alpha}_{ij}.\overline{\alpha}_{kl} = \overline{e}_i.\overline{e}_j.\overline{e}_k.\overline{e}_l 
\eeqa
which, by using \eqref{laws}, plus that fact that we are computing modulo $k.e_0$,  can be easily checked to give:
\beqa\label{m1}
\overline{\alpha}_{ij}.\overline{\alpha}_{kl} = b_\tq(i,l) \overline{\alpha}_{jk} - b_\tq(j,l)\overline{\alpha}_{ik} - b_\tq(i,k)\overline{\alpha}_{jl} + b_\tq(j,k) \overline{\alpha}_{il}.
\eeqa
Similarly, computing $\overline{\alpha}_{kl}.\overline{\alpha}_{ij} = \overline{e}_k.\overline{e}_l.\overline{e}_i.\overline{e}_j$, we see easily that 
\beqa\label{m2}
[\overline{\alpha}_{ij},\overline{\alpha}_{kl}] = 2.b_\tq(i,l).\overline{\alpha}_{il} -2.b_\tq(i,k)\overline{\alpha}_{jl} + 2.b_\tq(i,l)\overline{\alpha}_{jk} - 2.b_\tq(j,l)\overline{\alpha}_{ik}. 
\eeqa
Again, if we allow $j = k$, then a similar computation yields:
\beqa\label{m3}
[\overline{\alpha}_{ij},\overline{\alpha}_{jl}] = 2 \tq(e_j) \overline{\alpha}_{il} - b_\tq(l,j)\overline{\alpha}_{ij} -b_\tq(i,l) \overline{\alpha}_{jl}.
\eeqa
It is clear that {\em if $2$ is invertible}, then we can recover the matrix for the bilinear form $b_\tq$ from the matrix above. This shows injectivity. Working with valued points, it is immediate that its is an {\em immersion}

{\em The morphism $\theta$ is proper}: Let $R$ be a discrete valuation ring and $K = \text{Fract}(R)$, with residue field $k$. Let $A$ be an algebra structure on $\Lambda(\bbf) \otimes_{_k} R$, with identity element $e_0 \otimes 1$. Suppose that the algebra structure $A \otimes_{_R} K$ on $\Lambda(\bbf) \otimes_{_k} K$ lies in the image of $\theta(K)$. This implies that the values $\{ b_\tq(i,j)\}_{_{i,j}}$ all lie in $K$.  Since $A$ is given to be an $R$-algebra structure, the coefficients in the multiplication relations \eqref{m1}, \eqref{m2} and \eqref{m3} with respect to the basis \eqref{liebasis} all lie in $R$, i.e.,  the entries $\{ b_\tq(i,j)\}_{_{i,j}}$ all lie in $R$. This implies that the bilinear form $b_\tq$ and hence the quadratic form $\tq$ take its values in $R$. Denote it by $\tq_{_R}$. It is immediate that the algebra structure $A$ on $\Lambda(\bbf) \otimes_{_k} R$ is isomorphic to the Clifford algebra $\Cl^{+}_{{\tq_{_R}}}$.   This shows  that $\theta$ is proper by the valuative criterion  (see \cite[page 166]{desing}). Thus, $\theta$ is a closed embedding.

\epr
\subsubsection{Some remarks on Clifford bundles}
Let us fix some notations which we will use throughout. Let 
\beqa\label{opposite}
\ta_\tq^{\circ} := \Cl^{+,\circ}_\tq
\eeqa be the opposite of the positive graded part of the Clifford algebra and let 
\beqa\label{oppositeunits}
\cu_\tq \subset \ta_\tq^{\circ}
\eeqa
denote the "group of units", namely, the invertible elements in the opposite associative algebra $\Cl^{+, \circ}_\tq$. This process can be globalised. We work over the affine space $\sf Q \simeq S^{^2}(\bbf^*)\simeq \ba^{^{{m(m-1)}\over {2}}}$, the space of all quadratic forms on $\bbf$. Let $\Lambda(\bbf)_{\sf Q}$ be the trivial vector bundle over $\sf Q$ based on the $k$-vector space $\Lambda(\bbf)$. Then we can endow this vector bundle with an algebra structure  and obtain the Clifford algebra-bundle $\Cl_{\tq_{_{\cq}}}$ over $\cq$  associated to a quadratic form ${\tq_{_{\sf Q}}}:{\bbf} \otimes_k \cq \to k_{_{\sf Q}}$. This algebra bundle has a $\bz/2$-grading and we have the positive graded algebra bundle $\Cl^{+}_{\tq_{_{\sf Q}}} \to \sf Q$. The whole construction can be carried out with the opposite algebra bundle $\ta_{\tq_{_{\sf Q}}}^{\circ}$; this  is  in fact isomorphic to $\Cl^{+}_{\tq_{_{{\sf Q}_{_o}}}}$ over ${\sf Q}_{_o}$,  since we have a non-degenerate form on ${\bbf} \otimes_k {\sf Q}_{_o}$. The resulting algebra bundle is an Azumaya algebra over $ {\sf Q}_{_o}$.

It is known that the group scheme of units of this unital associative algebra bundle $\ta_{\tq_{_{\sf Q}}}^{\circ}$ is a flat, in fact smooth, affine group scheme which we denote by $\cu_{_\cq} \to \cq$ (see for example \cite[Proposition 2.4.2.1]{cf}). In general, \'etale locally, the quadratic form ${\tq_{_{{\sf Q}_{_o}}}}$ is the standard split hyperbolic quadratic form $\tq^{\text{h}}$, and so the group scheme $\cu_{_{\cq_{_o}}}$ is a {\em form} of the split group scheme associated to $\tq^{\text{h}}$.

We observe that if the quadratic form $\tq$ over $k$ is non-degenerate we have the isomorphism $\ta_\tq^{\circ} \simeq \cm_n$, where $n = 2^{^\ell}$. Hence, the group of units $\cu_\tq$ is isomorphic to the linear group $\GL(2^{^{\ell}})$.  We have a closed embedding 
\beqa\label{unitsinthebiggroup}
j_{_{\cq}}: \cu_{_\cq} \hra \GL(n^2)_\cq
\eeqa (see Remark \ref{strgrp}). %We note that the group scheme $\cu_{_\cq}$ over degenerate quadratic forms $\tq$ is {\em non-reductive} in general. 

Let $\tq$ be a non-degenerate quadratic form on $\bbf$ over $k$, and let $\spin(\tq)$ be the standard (odd) Spin group. This simply connected group can be realised as a subgroup of the group of units of the positive graded $\Cl^{+}_\tq$; this {\em faithful irreducible representation} is just the  standard {\em spinor} representation. 

More precisely, we have a decomposition $\bbf = N \oplus P \oplus U$, where $N,P$ are isotropic subspaces of dimension $\text{dim}(N) = \text{dim}(P) = \ell$, while $U$ is a $1$-dimensional subspace orthogonal to $N$ and $P$.  This decomposition identifies the Clifford algebra $\Cl_\tq$ with the direct sum $\text{End}(\Lambda N) \oplus \text{End}(\Lambda P)$ and the positive graded part $\Cl^{+}_\tq$ gets identified with $\text{End}(\Lambda N)$. Thus, we can identify $\GL(S)$ with the units of $\Cl^{+}_\tq$.

Let $S := \Lambda N$, so that $\text{dim}(S) = 2^{^\ell}$. Then the spinor representation is denoted by:
\beqa\label{thespinorrep}
{\Large\text{\cursive s}}:\spin(\tq) \hra \GL(S). 
\eeqa
This is a {\em faithful and irreducible representation} if we work with the case when $m = 2\ell+1$ is odd.

By applying the canonical involution on the Clifford algebra, we can identify the Spin group as a subgroup \eqref{oppositeunits}
\beqa\label{spininoppunits}
{\Large\text{\cursive s}}:\spin(\tq) \hra  \cu_\tq.
\eeqa
This inclusion is again {\em irreducible} and we have made an abuse of notation to denote it by ${\Large\text{\cursive s}}$ and call this the {\em spinor representation}.

More generally, let $K$ be an extension of $k$ and let $\tq_{_K}:\bbf_{_K} \to K$ be a non-degenerate quadratic form. Then, the spinor representation is given by ${\Large\text{\cursive s}}:\spin(\bbf_{_K}, \tq_{_K})) \hra \cu_{_{\tq_{_K}}}$, where the $K$-group scheme $\cu_{_{\tq_{_K}}}$ is a twisted form of $\cu_{_{\tq^{^{\tt h}}_{_K}}}$, which as an algebraic group is isomorphic to $\text{SL}(n)$. 

The image of the spinor representation is given as follows: for every $K$-algebra $R$
\beqa
{\Large\text{\cursive s}}_{_R}\big(\spin(\bbf_{_K}, \tq_{_K})(R)\big) = \{g_{_R} \in \Cl^{+,\times}_\tq(R) \mid g_{_R}.\bbf_{_R}.g_{_R}^{^-1} \subset \bbf_{_R}, g_{_R}.\tau(g_{_R}) = 1\}
\eeqa
where $\tau$ is the standard involution on the Clifford algebra.

\subsubsection{On the universal family}\label{onunivfam} We return to the setting of \eqref{univfamilyinseshadri}. In the situation of the moduli space $N_{_n}$, with $n = 2^{^\ell}$, we observe that we have a group scheme $\cu_\tq$ on $N_{_n}$ where $\tq$ is a quadratic form on $\bbf$. This quadratic form is such that on the open subset $N_{_n}^s$ of stable bundles, is non-degenerate but {\em non-split}, i.e., it is not the standard hyperbolic form. Let us denote the standard hyperbolic quadratic form by the $\tq^{{\tt h}}$. Moreover, it is shown in \cite{desing} that the $\text{GL}(n^{^2})$-family on $C \times N_{_n}$ in fact gets a reduction to the {\em twisted} subgroup scheme $\cu_\tq$. When restricted to the stable locus we therefore get a {\em universal family} $E_{_{\tq}}$ on $C \times N_{_n}^{^{s}}$ which is a torsor for the group scheme $\cu_\tq$ restricted to the stable locus. 

By an abuse of notation, we will continue to denote this {\em non-split} non-degenerate quadratic form $\tq$ over $N_{_n}^{^{s}}$ by $\tq$ itself, and we will reserve the notation $\tq^{{\tt h}}$ for the standard split hyperbolic form. 

We make an important observation here on the family $E_{_{\tq}}$ on $C \times N_{_n}^{^{s}}$. Firstly, it is well-known that there is no universal family of rank $n$ on the moduli space of stable bundles. Hence, this group scheme $\cu_\tq$ over $N_{_n}^s$ cannot be the one coming from the split form  $\tq^{{\tt h}}$, for that would result in a universal torsor for the group scheme $\cu_{\tq^{{\tt h}}}$ which is in fact isomorphic to $\text{SL}(n)$. Thus, the group scheme $\cu_\tq$ is a {\em form} of $\cu_{\tq^{{\tt h}}}$. Moreover, if we use the classical Quot scheme construction of the moduli space $N_{_n}^s$, now viewed as the moduli space of stable vector bundles of rank $n$, we get a principal $\text{PGL}$-bundle 
\beqa
\psi:\text{Quot}_{_n}^s \to N_{_n}^s.
\eeqa
This quotient map has two key properties: 
\begin{itemize}
\item The group scheme $\cu_{_{\tq,N_{_n}^s}}$ pulls back under $\psi$ to the split group scheme $\cu_{_{\tq^{{\tt h}},\text{Quot}_{_n}^s}}$
\item The $\cu_{_{\tq,N_{_n}^s}}$-torsor  $E_{_{\tq}}$ on $C \times N_{_n}^s$ pulls back to $\mathcal E_{_{\tq}}$ on $C \times \text{Quot}_{_n}^s$, the universal family on the Quot scheme.
\end{itemize}
More precisely, the universal bundle $\mathcal E_{_{\tq}}$ does not descend as an $\text{SL}(n)$-torsor, but it descends as a torsor under the twisted group scheme $\cu_{_{\tq,N_{_n}^s}}$.
 
Let $\cq_{\tt o} \subset \cq$ as before be the {\em open sub-scheme} of non-degenerate quadratic forms. Then it is easily seen that the inclusion \eqref{spininoppunits} can be extended to a closed embedding 
\beqa\label{spin2}
H_{_{\cq_{\tt o}}} \subset \cu_{_{\cq_{\tt o}}}
\eeqa
of affine {\em reductive} group schemes.  

%Since $H_{_{\cq_{\tt o}}}$ is a reductive group scheme over $\cq_{\tt o}$, it is \'etale locally split \cite[XXII, Cor 2.3]{sga3} over $\cq_{\tt o}$; whence, if $t \in \cq_{\tt o}$, the group scheme is {\em split} over the completion $\hat{\co}_{_{\cq,t}}$. In particular, if $L$ is the fraction field of $\hat{\co}_{_{\cq,t}}$, the group scheme $H_{_L} = \spin(\tq_{_{\cq}})$ has a {\em split} maximal torus $T$.  We will assume we are working with the quadratic form over $L$ with maximal Witt index. 

\section{Towards the main Theorem \ref{main}}\label{towardsmain} We work with the notations of \S\ref{seshadri}. In this section we will denote by the letter $H$ the odd spin group $\spin (m)$ and we let $m = 2\ell +1$. 

Let  $E$ be a stable $H$-bundle. Let $\rho:\pi_{_1}(C,x) \to K_{_H}$ be the unitary representation which gives rise to $E$ under the Narasimhan-Seshadri correspondence, \`a la Ramanathan. Let $\text{Hol}(E) := \overline{\text{Im}(\rho)}$ be the Zariski closure of the image of $\rho$ in $H$. We say a stable bundle $E$ has {\em full holonomy} if we have the identification $\text{Hol}(E) = H$. This set is non-empty (for the proof, see \cite[Page 361]{bbn}; this fact holds for {\em arbitrary semisimple groups}, and hence in the proof of \eqref{mainexceptional}  for $\mathfrak G$ as well. An essential remark is due to Nori, which says that in a compact connected real semi-simple Lie group, there exists a dense subgroup generated by two general elements (for a proof see \cite[Lemma 3.1]{su}).

It is a fact that all points $E$ which correspond to stable bundles with full holonomy have only the center of $H$ as automorphisms and are hence smooth points in the moduli of semi-stable principal $H$-bundles.  {\em Let $M_{_{\tt f}}(H) \subset M_{_{H}}^s \subset M^{ss}_{_H}$ be the open smooth locus of stable bundles with full holonomy}.

%This induces an embedding $\text{Lie}({\Large\text{\cursive s}}):\text{Lie}(H) \hra \text{Lie}(\cu)$.

Let ${\Large\text{\cursive s}}:H \hra  \cu$ be the {\em spinor} representation as in \eqref{spininoppunits}; recall that $\text{dim}(S) = n = 2^{^\ell}$. Since ${\Large\text{\cursive s}}$ is irreducible, by an application of Schur's lemma and the Narasimhan-Seshadri theorem, we conclude that {\em if $E$ is stable with full holonomy, then the $\cu$-bundle ${\Large\text{\cursive s}}_{_*}(E) := E \times ^{^H} \cu$ obtained by extension of structure group by the homomorphism ${\Large\text{\cursive s}}$, continues to remain a stable $\cu$-bundle, i.e., a stable $\SL(n)$-bundle}. Thus, ${\Large\text{\cursive s}}_{_*}$ induces a morphisms 
\beqa
{\Large\text{\cursive s}}_{_*}: M_{_{\tt f}}(H) \to M_{_{\SL(n)}}^s\\
{\Large\text{\cursive s}}_{_*}: M^{ss}_{_H} \to M_{_{\SL(n)}}^{ss}.
\eeqa 
By \S2 (see \cite{desing}), we have a {\em proper and birational} morphism $N_{_n} \to M_{_{\SL(n)}}^{ss}$ and an identification $M_{_{\SL(n)}}^s \simeq N^s_{_n}$. We also have the following commutative diagram:
\beqa\label{seshadri1}
\begin{tikzcd}
	{N_{_n}} \\
	{M_{_{\SL(n)}}^{ss}} & {M_{_{\GL(n^2)}}^{ss}}
	\arrow[from=1-1, to=2-1]
	\arrow[from=1-1, to=2-2]
	\arrow[hook, from=2-1, to=2-2]
\end{tikzcd}
\eeqa

We also have the open immersion ${\Large\text{\cursive v}}:M_{_{\tt f}}(H) \hra M^{ss}_{_H}$. Therefore, we get a morphism $\phi = ({\Large\text{\cursive s}}_{_*}, {\Large\text{\cursive v}}): M_{_{\tt f}}(H) \to M^{ss}_{_H} \times N_{_n}$. Let $\Gamma_{_\phi}$ denote the base-change (with its reduced induced structure), and we have the diagram:
\beqa\label{gammaphi}
\begin{tikzcd}
	{\Gamma_{_\phi}} & {N_{_n}} \\
	{M^{ss}_{_H}} & {M^{ss}_{_{SL(n)}}}
	\arrow[from=1-1, to=2-1]
	\arrow["{\Large\text{\cursive s}}_{_*}", from=2-1, to=2-2]
	\arrow[from=1-1, to=1-2]
	\arrow[from=1-2, to=2-2]
\end{tikzcd}\eeqa
Thus, $\phi$ induces canonically an open embedding $\phi:M_{_{\tt f}}(H) \to \Gamma_{_\phi}$. Let $\fm_{_H}' := \text{normalisation~of}~\Gamma_{_\phi}$. Note that this is a disjoint union of irreducible components. Let
\beqa\label{thecandidate}
\fm_{_H}:= \text{component~containing}~M_{_{\tt f}}(H).
\eeqa
We get the following diagram, with projections $p$ and $q$. 
\beqa\label{thecorresp}
\begin{tikzcd}
	& \fm_{_H} \\
	M^{ss}_{_H} && N_{_n}
	\arrow["{^p}", from=1-2, to=2-1]
	\arrow["{^q}"', from=1-2, to=2-3]
\end{tikzcd}
\eeqa
Note that by construction, both projections are {\em projective morphisms}.

\subsubsection{On the image of $M_{_{\tt f}}(H)$ in $N_{_n}^s$}\label{imageofmfh} Let us take a closer look at the image, especially the representations which are involved in taking the associated bundles. Recall that a point in $N_{_n}^s$ is a vector bundle $W$ of rank $n^2$ which isomorphic to $V \otimes \co_{_C}^n$, where $V$ is a stable vector bundle of rank $n$. By \eqref{strgrp}, the structure group of $W$ can be reduced to the group of units $\cu$ of the endomorphism algebra of $W$. Hence $\text{Lie}(\cu) = {\mathcal L}(\End(W))$ \cite[2.8.0.41]{cf}, where ${\mathcal L}(\End(W))$ is the  natural Lie algebra structure from the associative algebra structure on $\End(W)$. 

%Let $F$ be the principal $\cu$-bundle thus obtained. We thus have the equality $H^{^0}(C, F(\text{Lie}(\cu))) = \End(W)$.

Suppose now that we have an $E \in M_{_{\tt f}}(H)$, which is a stable $\spin(\tq)$-bundle with full holonomy. Let us be more precise and we let ${\Large\text{\cursive s}}:\spin(\tq) \hra  \cu_{_\tq}$ be the {\em spinor} representation as in \eqref{spininoppunits}. Let $W_{_E}$ denote the image of $E$ in $N_{_n}^s$. Then, as we have observed above, we have a canonical identification of ${\mathcal L}(\End(W_{_E}))$  with $\text{Lie}(\cu_{_\tq}) =  {\mathcal L}(\Cl^{+,\circ}_\tq)$. From the proof of \eqref{forversality}, it is clear that as a point of $\ca_{_n}$, $\End(W_{_E})$ lies in the image $\theta(\cq_{_o})$. All this can also be done in families. 

\subsubsection{More on structure group reduction} Let $\xi \in \fm_{_H}$ be a point. Let $R_{_\xi} := \hat\co_{_{{\fm},\xi}}$, which is a  normal local domain, and let $Z = \spec R_{_\xi}$;  let $K_{_\xi} = \text{Fract}(R_{_\xi})$ be the field of fractions which corresponds to the {\em non-empty} open sub-scheme in $Z$ comprising of stable $H$-bundles with full holonomy. The second projection $1 \times q$ gives a $\GL(n^2)$-bundle $(1 \times q)^*(P_{_n})$ on $C \times Z$. This in turn gives us a family of $\co_{_Z}$-algebras $\fb_{_n}$, which is a  specialisation of the matrix algebra $\cm_{_n}$. We have the following proposition on reduction of structure groups over the generic point.

\bprop\label{stablereduction1} Let $K = K_{_\xi}$, and let $P_{_K}$ be the universal principal $\GL(n^2)_{_{K}}$-torsor restricted to $C_{_K} := C \times \spec K$. We assume that $P_{_K}$ is in fact a $\spec K$-valued point of the moduli space  $N_{_n}^s(K)$. Let $H_{_K}$ be a connected semisimple $K$-group of type $\tt B_{_n}$ or $\tt D_{_n}$. Then, we have a reduction of structure group of $P_{_K}$ to $H_{_K} \subset \GL(n^2)_{_K}$ through the chain of inclusions  $$H_{_K} \stackrel{{\Large\text{\cursive s}}_{_K}}  \hra \cu_{_{K}} \stackrel{j_{_K}}\hra \GL(n^2)_{_K}$$ \eqref{spininoppunits}.  \eprop 
\bpr Let $\spec L \to N_{_n}^s$ denote the generic point, where $n = 2^{^\ell}$. We recall that the group scheme of units, $\cu_{_\tq}$ is defined on $L$, where $(\bbf_{_L},\tq)$ is a quadratic space.  We will denote it by $\cu_{_L}:= \cu_{_{\tq,L}}$. We also have a chain of group scheme inclusions:
$$H_{_L} \stackrel{{\Large\text{\cursive s}}_{_L}}  \hra \cu_{_{L}} \stackrel{j_{_L}}\hra \GL(n^2)_{_L}$$
where $H_{_L}$ is the Spin group $\spin(\bbf_{_L},\tq_{_L})$.

We prove this in the case when the group scheme $H_{_L}$ is of type $\tt B_{_n}$, the other case being similar.  Recall \eqref{onunivfam}. Let $E = E_{_\tq}$ be the $\cu_{_{\tq,L}}$-torsors on $C \times \spec L$; recall that it is realized as a reduction of structure group of the universal family $P_{_n}$ to $C \times \spec L$. Let
\beqa
F:\text{Sch}\big/L \to \text{Sets}
\eeqa
be defined as follows: let $T \to \spec L$ be in $\text{Sch}\big/L$. Define 
$$F(T) := \left \{ \begin{array}{l} \mbox{isomorphism classes of
pairs $(E_{_T},s)$, $E_{_T} = \{ E_t\}_{t \in T}$}\\ \mbox{a family of
$\cu_{_{\tq,T}}$-bundles in \mbox{$N_{_n}^s$} and} \\ \mbox{$s = \{ s_t \}_{t \in T}$ a section of $E(\cu/H)$ on $C \times T$.}\\  
\end{array} \right \}$$
Since $H_{_L} \subset \cu_{_L}$ is an inclusion of connected reductive groups we conclude, using \cite[Lemma 4.8.1]{r1}, that $F$ is representable by a  $L$-scheme $\Gamma_{_L}$, which is of finite type and which gives the structure morphism $\epsilon:\Gamma_{_L} \to \spec L$; whence we get the following diagram of principal $\text{PGL}$-bundles \eqref{onunivfam}:
\beqa\label{nannairukku}
\begin{tikzcd}
	{\text{Quot}_{_{H_{_L}}}} & {\text{Quot}^s_{_L}} \\
	{\Gamma_{_L}} & {\spec L}
	\arrow[from=1-1, to=2-1]
	\arrow["\epsilon", from=2-1, to=2-2]
	\arrow[from=1-1, to=1-2]
	\arrow["\psi",from=1-2, to=2-2]
\end{tikzcd}\eeqa
It follows easily from the discussion in \eqref{onunivfam} that ${\text{Quot}_{_{H_{_L}}}}$ is in fact the scheme which represents the functor of reductions of structure group of the $\cu_{_{\tq^{^{\text{h}}}}}$-torsor $\ce := \psi^{^*}(E)$ on $C \times \text{Quot}_{_L}^s$ for the inclusion $\spin(n) \hra \cu_{_{\tq^{^{\text{h}}}}}$. In other words, the $\text{PGL}$-quotient of ${\text{Quot}_{_{H_{_L}}}}$ is precisely the inverse image of $\spec L$ in the moduli space $M_{_{\tt f}}(H)$ of $\spin(n)$-bundles with full holonomy. In other words we have an iinclusion $\Gamma_{_L} \hra M_{_{\tt f}}(H)$. Whence, we also have a universal reduction of structure group $E_{_H}$ of $E$ to the subgroup scheme $H_{_L} \hra \cu_{_{\tq,L}}$ ocer $C \times \Gamma_{_L}$. 

\epr

\subsubsection{The versality of the group scheme $\bh_{_\cq}$}\label{versality}  Let $\xi \in \fm_{_H}$ be a point. Let $R_{_\xi} := \hat\co_{_{{\fm},\xi}}$, which is a  normal local domain, and let $Z = \spec R_{_\xi}$;  let $K_{_\xi} = \text{Fract}(R_{_\xi})$ be the field of fractions which corresponds to the {\em non-empty} open sub-scheme in $Z$ comprising of stable $H$-bundles with full holonomy. The second projection $1 \times q$ gives a $\GL(n^2)$-bundle $(1 \times q)^*(P_{_n})$ on $C \times Z$. This in turn gives us a family of $\co_{_Z}$-algebras $\fb_{_n}$, which is a  specialisation of the matrix algebra $\cm_{_n}$. By \eqref{versalproperty}, we get a morphism $\vartheta:Z \to \ca_{_n}$. Moreover, since $(1 \times q)^*(P_{_n})$ has a reduction of structure group to $H_{_{K_{_\xi}}}$ on $C_{_{K_{_\xi}}}$ by Proposition \eqref{stablereduction1}, the morphism $\vartheta$ is such that the generic point $\spec K_{_\xi}$ maps to $\theta(\cq_{_o})$, under the embedding $\theta:\cq \hra  \ca_{_n}$ (see \eqref{forversality}). 

By \eqref{forversality} again, the rational map $\vartheta: \spec K_{_\xi} \to \cq$ extends at all height $1$ primes in $R$. Since $Z$ is a normal domain, by \cite[Lemma 2, page 109]{blr} the morphism extends to the whole of $Z$ to give the basic "versal morphism" 
\beqa\label{versalmorph}
\vartheta:Z \to \cq.
\eeqa

\subsubsection{The structure group reduction}\label{step1}
{\em  We will now use the results proven in Part II and the notation in \eqref{overq}}. Let $B = \hat{\co}_{_{\cq,t}}$ and $L:= \text{Fract}(B)$. Let $\vartheta:(Z,\xi) \to (\spec B,t)$ be the morphism \eqref{versalmorph}. Note that $\vartheta$ maps $\spec K_{_\xi}$ to $\spec L$. Put $R = R_{_\xi}$ and $\bh_{_{R}} := \vartheta^*(\bh_{_{B}})$ and let $H_{_{K}} := \vartheta^*(H_{_{L}})$. Let 
\beqa\label{hbt2}
j_{_R} \circ {\tilde{\Large\text{\cursive s}}}_{_R}:\bh_{_{R}} \hra \cu_{_{R}} \hra \GL(n^2)_{_{R}}
\eeqa 
be the pull-back of the chain of inclusions \eqref{hbt1}.

That we do indeed get such a chain of inclusions after pull-back by $\vartheta$ needs to be verified. We argue as follows. The only thing which needs checking is the inclusion $\bh_{_{R}} \hra \GL(n^2)_{_{R}}$, the second one coming from the group scheme of units of the algebra bundle on $Z = \spec R$. That we do get a morphism of flat group schemes on $Z$ is clear. Also, over the generic point $\spec K$ of $Z$, we clearly get an inclusion $\bh_{_K} \hra  \GL(n^2)_{_{K}}$. Let $\mathfrak K_{_R} := \text{Ker}(\bh_{_{R}} \to \GL(n^2)_{_{R}})$. We see that $\mathfrak K_{_K}$ is the identity group scheme. Note that the Lie algebra of $\bh_{_K}$ extends to the whole of $Z$ simply as $\vartheta^{^*}(\lie(\bh_{_B}))$. This follows what we did over $\spec B$  and the Lie algebra extension from \eqref{lieextension}. 

Once the Lie algebra bundle extends, the proof of Theorem \ref{hbt} goes through over $Z$ as well since the only assumption which was need for this extension was the noetherianness of $Z$. In other words, we get a smooth connected quasi-affine subgroup scheme $\bh'_{_R}$ extending $\bh_{_K}$. By the uniqueness of such group schemes \eqref{hbt}, it follows that we have an isomorphism between $\bh_{_R}$ and $\bh'_{_R}$.

%We can conclude the same by noting that have the same Lie algebras which gives an identification of a neighbourhood of the identity sections and then we simply appeal to \cite[1.2.14, page 20]{bt2}.
The aim now is to prove the following central result.

\bth\label{stablereduction} Let $Z = \spec R$ and let $P_{_R}$ be a principal $\GL(n^2)_{_{R}}$-torsor on $C_{_R} := C \times Z$. Under the assumption that $P_{_R}$ is a $\spec R$-valued point of the moduli space $N_{_n}$ and is such that the restriction $P_{_K}$ lies in the stable locus $N_{_n}^s(K)$. 
Suppose that we have a reduction of structure group of $P_{_K}$ to $H_{_K}$ through the chain of inclusions  $ H_{_K} \stackrel{{\Large\text{\cursive s}}_{_K}}  \hra \cu_{_{K}} \stackrel{j_{_K}}\hra \GL(n^2)_{_K}$ \eqref{spininoppunits}. Then,  the structure group of $P_{_R}$ can be reduced to $\bh_{_R}$.  \eeth 

\bcor\label{stablereductionBD} In the case when $H$ is of type $B_{_n}$ or $D_{_n}$, the structure group of the  $\GL(n^2)_{_{R}}$-torsor $P_{_R}$ reduces to the group scheme $\ch_{_R}$. \ecor
%We begin by making some key observations. Let us denote by $G$ the linear group $\GL(n^2)$. 

%remarks on twisting the group scheme $H_{_L}$ in terms of the torsor. 

%\subsubsection{When the base is a smooth curve} {\em We work with the fixed base point $x \in C$ as in \eqref{basepoint}}. We have the following: 
We prove the following proposition before proving the Theorem \ref{stablereduction} and derive Corollary \ref{stablereductionBD}.
\bprop\label{overht1} Let $\mathfrak p$ be a prime ideal of height $1$ in $Z$. Let $(R, \mathfrak m)$ be a discrete valuation ring with $R = \co_{_{\mathfrak p}}$, so that $\text{Fract}(R) = K = K_{_\xi}$. Let $R/\mathfrak m = \kappa$ and let $c:\spec R \to Z$ be the induced morphism such that $c(K)$ lies in $U \subset Z$ and $c(\kappa) = \mathfrak p$.  Put $P_{_R} := c^*(P_{_Z})$. Then, we get a reduction of structure group of $P_{_R}$ to the schematic closure of  $H_{_K}$ in $\GL(n^2)_{_R}$ extending the given reduction over $C_{_K}$. \eprop
\bpr By assumption, since $c(K)$ lies in $U$,  we have a reduction of structure group of $P_{_R}$ to $H_{_K}$ on $C_{_K} = C \times \spec K$.

%As a first step we get a reduction of structure group  of the restriction $P_{_{x,R}}$ on $x \times \spec R$ to the schematic closure of $H_{_K}$ in $\GL(n^2)_{_R}$, which extends the given reduction over $x \times \spec K$. 

Let  $E_{_K} \subset P_{_K}$ be this given $H_{_K}$-reduction. The extension of $P_{_K}$ to $P_{_R}$ gives a $K$-automorphism:
\beqa
\phi_{_K}:P_{_K} \to P_{_K}
\eeqa
such that after modifying by $\phi_{_K}$, the transition functions of $P_{_K}$ extend to transitions functions of $P_{_R}$ over $C_{_R}$. 

Let $\theta_{_{ij,K}}:U_{_{ij,K}} \to H_{_K}$ be transition functions of $E_{_K}$. The automorphism $\phi_{_K}$ is given in the open cover $U_{_{i,K}}$ be $\{\phi_{_i}\}$ which act on the transition functions $\theta_{_{ij}}$, viewed as transition functions of $P_{_K}$. Then we have new transition functions $\delta_{_{ij,K}}$ of $P_{_K}$ which satisfy the following relations:
\beqa\label{trfns}
\delta_{_{ij,K}} = \phi_{_i}.\theta_{_{ij,K}}.\phi_{_j}^{^{-1}}.
\eeqa
Let us keep the notation $j_{_K} \circ {\Large\text{\cursive s}}_{_K}:H_{_K} \hra \cu_{_K} \hra \GL(n^2)_{_K}$ as in \eqref{hbt2}. 

By assumption, the $\GL(n^2)_{_R}$-bundle $P_{_R}$ is the frame bundle of a point $W_{_R} \in N_{_n}(R)$. Here the vector bundle $W \in N_{_n}$ is as in \S2. Thus, we can view $\phi_{_K}$ as an automorphism of $W_{_K}$ on $C_{_K}$ and the transition functions $\delta_{_{ij,K}}$ and $\theta_{_{ij,K}}$ as transition functions of $W_{_K}$. The added property of $\delta_{_{ij,K}}$ is that, there exist transition functions $\delta_{_{ij,R}}$ over $C_{_R}$, which extend $\delta_{_{ij,K}}$. 

Let $A = A(W)$ be the endomorphism algebra $\text{End}(W)$ of $W \in N_{_n}$ and let $A_{_r}^{^{\vee}}$ be the dual of the {\em right regular representation} of the algebra $A$. Then one sees easily (see \cite[Proposition 3]{desing}) that the opposite algebra $A^{^{\circ}}$ can be identified with $\text{End}_{_{A-mod}}(A_{_r}^{^{\vee}})$. Equivalently, if we take the inclusion $A \hra \text{End}_{_{k-mod}}(A_{_r}^{^{\vee}})$, we may view $A^{^{\circ}}$ as the centralizer $\tiny{\tc{C}}(A)$ of $A$ in the algebra $\text{End}_{_{k-mod}}(A_{_r}^{^{\vee}})$. 

Let $x \in C$ be any point on the curve. Recall that, since $W$ is a semistable vector bundle of $\text{deg}(W) = 0$, we have the inclusion $\text{End}(W) \hra \text{End}(W_{_x})$ given by the evaluation at $x$. Thus, we get an inclusion $A \hra \text{End}(W_{_x})$. As is observed in \cite{desing}, we can identify the $A$-module structure on $W_{_x}$ with $A_{_r}^{^{\vee}}$, and this is independent of the point $x \in C$.  All these identifications can be done over a parameter and therefore for an $x \in C$, we have the inclusion $\text{Aut}(W_{_K}) \hra \text{Aut}(W_{_{x,K}})$. As in \S2, we are able to make the following identifications: 
\beqa
\GL(n^2)_{_K} = \text{Aut}(W_{_{x,K}})\\
\cu_{_K} = \tiny{\tc{C}}(A_{_K})^{\times}.
\eeqa
Rewriting $j_{_K} \circ {\Large\text{\cursive s}}_{_K}:H_{_K} \hra \cu_{_K} \hra \GL(n^2)_{_K}$, we get the chain of inclusions:
\beqa\label{hbt3.1}
j_{_K} \circ {\Large\text{\cursive s}}_{_K}:H_{_K} \hra  {\tiny{\tc{C}}}(A_{_K})^{\times} \hra \text{Aut}(W_{_{x,K}}).
\eeqa
Evaluating the relation \eqref{trfns} at the point $x \in C$, we get the relation:
\beqa\label{trfns1}
\delta_{_{{ij},K}}(x) = \phi_{_i}(x).\theta_{_{ij,K}}(x).\phi_{_j}^{^{-1}}(x).
\eeqa
By the discussion above, we see that the automorphism $\phi_{_K} \in \text{Aut}(W_{_K})$ is determined by its evaluation at $x \in C$ and also $\phi_{_i}(x) = \phi_{_j}(x) \in \text{Aut}(W_{_{x,K}})$. Further, by \eqref{hbt3.1}, the subgroup $H_{_K}$ centralizes all automorphisms of $W_{_K}$ inside $\text{Aut}(W_{_{x,K}})$. In other words, the relation \eqref{trfns1} gives 
\beqa\label{thiseq}
\delta_{_{ij,K}}(x) = \theta_{_{ij,K}}(x).
\eeqa 
The identification of the $A$-modules $W_{_x}$ with $A_{_r}^{^{\vee}}$ shows in fact, that the equality \eqref{thiseq} {\em  holds for all $x \in C$}.

Since $\delta_{_{ij},K}$ are extendable,  we conclude that by  restriction to $x \times \spec K$, the transitions functions $\theta_{_{ij,K}}(x)$ are extendable to $\spec R$ and this holds {\em for all $x \in C$}. {\em A priori}, these take their values in $\GL(n^2)_{_R}$, while generically taking it in $H_{_K}$. But since $R$ is a discrete valuation ring, we conclude that the transition functions $\theta_{_{ij,K}}$ take their values in the schematic closure of $H_{_K}$ in $\GL(n^2)_{_R}$, i.e., in $\bh_{_R}$. In other words, we conclude that {\em for each $x \in C$}, the principal $\GL(n^2)_{_R}$-bundle $P_{_{x,R}}$  on $ x \times \spec R$ has a reduction of structure group to the schematic closure $\bh_{_R}$ of $H_{_K}$ in $\GL(n^2)_{_R}$. 

To complete the proof of \eqref{overht1}, i.e. globalize these structure group reductions,  let us recall a couple of facts, with notations as in Proposition \ref{overht1}:

(a) (\cite[Lemma 2.10, page 328]{basa}): Let $\mcw_{_R}$ be a family of semistable vector bundles of degree $0$ on $C_{_R}$ and let $s_{_K}:C_{_K} \to \mcw_{_K}$ be a section with the property that at any one point $x \in C$, the evaluation of the section $s_{_K}(x)$ extends to a section of $\mcw_{_{x,R}}$. Then the section $s_{_K}$ extends to a section $s_{_R}$ of  $\mcw_{_R}$ on $C_{_R}$. 

(b) (\cite[Lemma 9, page 18]{bapa} ("Chevalley embedding")): Consider the embedding $\bh_{_R} \subset \GL(n^2)_{_R}$. Then, there exists a finite dimensional $\GL(n^2)_{_R}$-module $\cm_{_R}$, such that the homogeneous space $\GL(n^2)_{_R} \big/ \bh_{_R}$ has an $\GL(n^2)_{_R}$-immersion in $\cm_{_R}$.

 Let $\eta:\GL(n^2)_{_R} \to \GL(\cm_{_R})$ denote the representation   which gives $\cm_{_R}$ the $\GL(n^2)_{_R}$-module structure in (b).  We now consider the associated vector bundle $\mcw_{_R} = P_{_R}(\cm_{_R})$ on $C_{_R}$. Since the principal bundle $P_{_R}$ lies in $N_{_n}$, this vector bundle gives a family of semistable vector bundles of degree $0$ on $C_{_R}$ (this follows immediately in characteristic zero, but for positive characteristics, we would require bounds on the charactertistic in terms of the Dynkin height of $\eta$, see \S\ref{charp}, \cite{bdp} and \cite{bapa}).  
 
 We are given a section $s_{_K}:C_{_K} \to P_{_K}\big(\GL(n^2)_{_K}\big/\bh_{_K}\big) \subset P_{_K}(\cm_{_K})$, given by the reduction of structure group to $\bh_{_K}$. By what has already been shown, for each $x \in C$, this section $s_{_K}(x)$ on $x \times \spec K$ extends  to  $x \times \spec R$ as a section $\spec R \to P_{_{x,R}}\big(\GL(n^2)_{_R}\big/\bh_{_R}\big) \subset P_{_{x,R}}(\cm_{_R})$.  By (a), we see that $s_{_K}$ therefore extends to a global section $s_{_R}:C_{_R} \to P_{_R}(\cm_{_R})$, while the values of this section lie in  $P_{_{x,R}}\big(\GL(n^2)_{_R}\big/\bh_{_R}\big)$ for each $x \in C$. This implies that $s_{_R}$ in fact gives a section $s_{_R}:C_{_R} \to P_{_R}\big(\GL(n^2)_{_R}\big/\bh_{_R}\big)$. This gives the required reduction of structure group of $P_{_R}$ to $\bh_{_R}$ on $C_{_R}$, completing the proof of Proposition \ref{overht1}.

\epr

\bpr (Completion of the proof of Theorem \ref{stablereduction} and Corollary \ref{stablereductionBD}). The principal $\GL(n^2)_{_{R}}$-torsor $P_{_R}$ on $C_{_R} := C \times Z$ has a reduction of structure group to $H_{_K}$ on $C_{_K}$. Using the language of transition function in Proposition \ref{overht1}, we see that the transition functions $\theta_{_{ij,K}}$ take their values in $H_{_K}$ and by Proposition \ref{overht1}, these functions extend to all height $1$ primes to take their values in $\bh_{_R}$. 

The scheme $Z$ is normal, the projection $C_{_R} \to Z$ is a smooth morphism, and the group scheme $\bh_{_R}$ is smooth by Theorem \ref{hbt}. Hence by Weil's theorem \cite[Theorem 1, page 109]{blr}, it follows that the transition functions in fact take their value in $\bh_{_R}$ over the whole of $C_{_R}$. In other words, we have a global reduction of structure group of $P_{_R}$ to $\bh_{_R}$.

Corollary \ref{stablereductionBD} follows immediately from Proposition \ref{stablereduction1} and Theorem \ref{stablereduction}.
\epr

Let $\ce_{_R}$ be the $\bh_{_R}$-torsor obtained by \eqref{stablereduction}.
%\blem\label{keyprop} Consider the canonical projection $p:Z \to M^{ss}_{_H}$. Let $z \in Z$ be a point and $\ce_{_z}$ be the $\bh_{_z}$-torsor on $C \times \spec k(z)$. Then there is a natural homomorphism $\mathfrak r: \bh_{_z} \to H$ over $\spec k(z)$ which factors as $\bh_{_z} \to L_{_z} \hra H$, where $L_{_z}$ is the Levi quotient of $\bh_{_z}$; furthermore, the associated $H$-bundle $\mathfrak r_{_*}(\ce_{_z})$ on $C \times \spec k(z)$ obtained by from  $\ce_{_z}$ by extension of structure via $\mathfrak r$ is semistable. In fact, this gives an explicit choice in the equivalence class determined by the image $p(z)$.

%\elem
\subsubsection{Some remarks on the fibres of the group scheme $\bh_{_B}$}
Following the notation in \eqref{overq}, let $B = \hat{\co}_{_{\cq,t}}$ and $L:= \text{Fract}(B)$ and let $\bh_{_B}$ be the higher Bruhat-Tits group scheme constructed on $\spec B$. 

Let $t$ be a point in $\spec B$. We now take a complete discrete valuation ring ${R} :=\mathbb k\llbracket s \rrbracket$ with residue field $\mathbb k$ and fraction field $K:= \mathbb k(\!(s)\!)$. Suppose we have a map $\spec R \to (\spec B, t) \subset \tsf{Q}$ such that the closed point maps to $t$, i.e., $k(t) = \mathbb k$ and the generic point $\spec K$ maps to the open locus $\tsf{Q}_{_o}$. Consider the group scheme $\bh_{_R}$, the pullback of $\bh$ from $\spec B$. Note that the base change $\bh_{_R}$ is such that over $\spec K$, we have an isomorphism $\bh_{_K} \simeq H_{_K}$.
\blem\label{oldwine} Under these assumptions, there is a ramified cover $R'/R$ obtained by taking $d^{^{th}}$-root of the uniformizer $s$ for some $d$, and a natural morphism
\beqa\label{faprime1}
{\mathfrak r}_{_{R'}}: \bh_{_{R'}} \to H_{_{R'}}
\eeqa
which gives a morphism 
\beqa\label{faprime1.5}
{\mathfrak r}_{_{\mathbb k}}:\bh_{_{\mathbb k}} \to  H_{_{\mathbb k}}
\eeqa
over the residue field $\mathbb k$.
\elem
\begin{proof} By \cite[Theorem 2.3.1]{base} (see also \cite[Proposition 8]{bapa}), there is a totally ramified extension $K'/K$, obtained by taking roots of the uniformizer $s$, such that 
if $R'$ the integral closure of $R$ in $K'$ then we have an inclusion
\beqa\label{faprime0} 
\bh_{_{R'}}(R') \subset M^{^{hyp}}.
\eeqa
where $M^{^{hyp}}$ is a hyperspecial parahoric arising as the stabilizer of a  hyperspecial vertex of the apartment determined by the maximal torus $T \subset H$. Thus we get a {\em semisimple} group scheme $H_{_{R'}}$ over $\spec R'$ such that $M^{^{hyp}} = H_{_{R'}}(R')$.

 Let $R^{^{sh}}$ be the {\em strict henselization} of $R$. Then by \cite[Proposition 10, page 50]{blr}, we see that $R^{^{sh}}$ is a strictly henselian discrete valuation ring with uniformiser $s$. The arguments above and \eqref{faprime0} then give a totally ramified cover $R'{^{^{sh}}}$ and an inclusion $\bh_{_{R^{^{sh}}}}(R'{^{^{sh}}}) \subset  H_{_{R^{^{sh}}}}(R'{^{^{sh}}})$. Since we have an $K$-isomorphism, $\bh_{_K} \simeq H_{_K}$, by (see \cite[1.7.6, page 39]{bt2}), this map in fact extends to a morphism \eqref{faprime1} and \eqref{faprime1.5}. \end{proof}

%When the residue field $\mathbb k$ of $R$ is {\em algebraically closed} or even  {\em separably closed}, the argument in {\it loc cit} shows that we have a morphism $\bh_{_{R'}}\to H_{_{R'}}$ extending the map $\bh_{_K} \to H_{_K}$. 

By an abuse of notation, let $\tc{s}$ also denote the composite  
\begin{equation}
\tc{s}:H \hookrightarrow \SL(n) \hra \GL(n^2)
\end{equation}
 By Lemma \eqref{oldwine} we have the morphism \eqref{faprime1}. On the other hand, if $L := \text{Frac}(R')$, then we have morphisms $\tc{s}_{_L}: H_{_L} \to \GL(n^2)_{_L}$ and $\tc{s}_{_{\mathbb k}}: H_{_{\mathbb k}} \to \GL(n^2)_{_{\mathbb k}}$. To globalize these two morphisms over $R'$, we may have to go to a \'etale extension $R''$ of $R'$, for by going to an \'etale extension $R''$ of $R'$, we see that $H_{_{R''}} = H \times R''$. We will therefore assume that we have a {\em spinor} inclusion $\tc{s}_{_{R'}}:H_{_{R'}}\hra \GL(n^2)_{_{R'}}$. 
 
By composing  \eqref{faprime1}  and the inclusion $\tc{s}$ we get a morphism of group schemes 
\beqa\label{newinclusionviamaxpar}
\psi_{_{R'}}:\bh_{_{R'}} \stackrel{{\mathfrak r}_{_{R'}}}\to H_{_{R'}} \stackrel{\tc{s}_{_{R'}}}\hra \GL(n^2)_{_{R'}}.
\eeqa
which is such that, if $L := \text{Frac}(R')$, then over $L$, the inclusion $\psi_{_L}:\bh_{_L} \simeq H_{_L} \hra \GL(n^2)_{_L}$ is the base change $\tc{s}_{_L}$ of $\tc{s}$ to $L$.
On the other hand, since $\bh_{_R}$ is the pull-back of $\bh_{_B}$, it comes with a natural immersion \eqref{hbt1}
\begin{equation}\label{degspinor0}
\tc{w}_{_R} := j_{_R} \circ \tilde{\tc{s}}_{_{R}}:\bh_{_R}  \hookrightarrow \GL(n^2)_{_R}.
\end{equation}
Over the closed fibre, this gives the {\em degenerate} spinor representation 
\begin{equation}\label{degspinor}
\tc{w}_{_\mathbb k}:\bh_{_\mathbb k}  \rightarrow \GL(n^2)_{_\mathbb k}.
\end{equation}
The morphisms $\tc{w}_{_{{R'}}}, \psi_{_{R'}}:\bh_{_{R'}}  \hookrightarrow \GL(n^2)_{_{R'}}$ coincide over the fraction field $L$ with $\tc{s}_{_L}$, and hence by \cite[1.2.4]{bt2}, they coincide everywhere. In particular, the if $\psi_{_\mathbb k}$ is 
\beqa\label{secondincl}
    \psi_{_\mathbb k}: \bh_{_\mathbb k}  \stackrel{{\mathfrak r}}\longrightarrow  H_{_\mathbb k} \stackrel{\tc{s}}\hookrightarrow \GL(n^2)_{_\mathbb k}
\eeqa
then, we get the identification $\psi_{_\mathbb k} =\tc{w}_{_\mathbb k}$.

\subsubsection{A key property of the $\bh_{_Z}$-torsor} 
As a corollary, we have the following bundle-theoretic statement. 
\bcor\label{verykeyprop}
    Suppose $\mathbb k$ is a field and we have a morphism $\spec \kappa \to \spec B$. Let $\bh_{_{\mathbb k}}$ and $\cu_{_{\mathbb k}}$ denote the pull backs  along with their associated morphism $\mathfrak r:\bh_{_{\mathbb k}} \to H_{_{\mathbb k}}$ \eqref{faprime1.5}. Suppose on $C_{_{\mathbb k}}:=C\times \spec \mathbb k$, we have an $\bh_{_{\mathbb k}}$-torsor $E$ whose extension of structure group to $\GL(n^2)$ via $ \tilde{{\tc{s}}}_{_\mathbb k}$ is semistable. Then
    \begin{itemize}
        \item The $H_{_{\mathbb k}}$-bundle obtained by extending structure group via $\mathfrak r$ is itself semistable.  
        \item $\mathfrak r_*(E)$ further extended to $\GL(n^2)$ and $\tilde{\tc{s}}_*(E)$ gives the same point in the moduli space $M^{ss}_{_{\GL(n^2)}}$ and hence in $M^{ss}_{_{\SL(n)}}$ (see \eqref{seshadri1} and \eqref{gammaphi}).
    \end{itemize}
    \label{cor}
\ecor
\begin{proof} Recall that the degenerate spinor map  $\tc{w}_{_\mathbb k}:\bh_{_\mathbb k} \to \GL(n^2)$ \eqref{degspinor} has the property that, starting with $E$, one knows that $\tc{w}_*(E)$ is a semistable $\GL(n^2)$-bundle of degree $0$. 

The identification $\psi_{_\mathbb k} = \tc{w}_{_\mathbb k}$ shows that $\psi_{_{\mathbb k,*}}(E)$ is semistable of degree $0$. Whence, by \eqref{secondincl}, $(\tc{s} \circ \mathfrak r)_*(E)$ is semistable of degree $0$. It follows that ${\mathfrak r}_{_{*}}(E)$ gives a semistable $H$-bundle. This gives a point in the moduli space $M^{ss}_{_H}$. It is now obvious that it gives the same point in $M^{ss}_{_{\GL(n^2)}}$. The embedding in the bottom arrow in \eqref{seshadri1} completes the proof of \eqref{cor}. \end{proof}

\section{The smooth compactification when $H$ is of type ${\tt B}_{_\ell}$} The notations in this section are as in \S\ref{towardsmain}.
The following is the first among the main results of this article and its proof  will occupy this entire section.
\bth\label{main} Let $H$ be an almost simple, connected, simply connected group of type ${\tt B}$. 
\begin{enumerate}
\item The scheme $\fm_{_H}$ is  smooth, irreducible, and projective.
\item The morphism $p:\fm_{_H} \to M^{ss}_{_H}$  is an {\em isomorphism} over the open locus $M_{_H}^{{rs}} \subset M_{_{H}}^s$ of smooth points of $M^{ss}_{_H}$, namely the regularly stable bundles.\end{enumerate}  \eeth
\subsubsection{The proof of the first part of Theorem \ref{main} (\text{a})} 
\bpr We need to verify the smoothness.  Let $D_{_1}$ be an Artin local ring and let $D_{_o}$ be a quotient of $D_{_1}$ with residue field $k$. Let $i_{_o}: \spec D_{_o} \to \spec D_{_1}$ be the induced closed immersion. Suppose that we are given a morphism $\phi_{_o}:\spec D_{_o} \to (Z, \xi)$, with the closed point mapping to $\xi$. Here $\xi \in {\fm}_{_H}$ and $Z$ as in \eqref{versality}. To prove that the scheme ${\fm}_{_H}$ is smooth at $\xi$, we need to find a lift $\phi_{_1}$ such that the following diagram commutes:
\[\begin{tikzcd}
	{\spec D_{_1}} \\
	{\spec D_{_o}} & {(Z,\xi)}
	\arrow["{_{\phi_{_o}}}", from=2-1, to=2-2]
	\arrow["{^{i_{_o}}}"', from=2-1, to=1-1]
	\arrow["{^{\phi_{_1}}}"', from=1-1, to=2-2]
\end{tikzcd}\]
By the Cohen structure theorem, we can find a complete regular local ring $D$, such that $D_{_1}$ is a quotient of $D$. Thus we have closed immersions $i:\spec D_{_o} \stackrel{i_{_o}} \to \spec D_{_1} \to \spec D$. It suffices therefore to lift $\phi_{_o}$ to a morphism $\phi:\spec D \to (Z, \xi)$. In other words, it is enough if  we get a $\phi$ such that the following  diagram commutes:
\[\begin{tikzcd}
	{\spec D} \\
	{\spec D_{_1}} & {^{\phi}} \\
	{\spec D_{o}} & {^{\phi_{_o}}} & {(Z, \xi)}
	\arrow[from=2-1, to=1-1]
	\arrow[from=3-1, to=2-1]
	\arrow[from=3-1, to=3-3]
	\arrow["{^{\phi_{_1}}}"', from=2-1, to=3-3]
	\arrow[from=1-1, to=3-3]
\end{tikzcd}\]

\bprop\label{case1} Let $\phi_{_o}:\spec D_{_o} \to (Z, \xi)$ be as above. Then, there exists $\phi:\spec D \to Z$ lifting $\phi_{_o}$. \eprop

\begin{proof}: Let $\vartheta:(Z,\xi)\to (\spec B, t)$ be the morphism obtained by ``versality", where $Z = \spec R$. Upon composing with $\phi_{o}$, one has the following commutative diagram.

\[\begin{tikzcd}
	& {(Z,\xi)} \\
	{\spec D_{o}} & {(\spec B,t)}
	\arrow["{\phi_{_o}}", from=2-1, to=1-2]
	\arrow["\vartheta", from=1-2, to=2-2]
	\arrow["{\alpha_0}"', from=2-1, to=2-2]
\end{tikzcd}\]
Set $\alpha_{_o}:=\vartheta\circ\phi_{_o}:\spec D_{_o}\to (\spec B,t)$. The closed point of $\spec D_{o}$ is mapped to $t\in \spec B$.

 By expressing $D$ as an inverse limit of Artinian quotients, and using smoothness of $\spec B$, the morphism $\alpha_{_o}$ lifts to $\alpha:\spec D\to (\spec B,t)$ as follows.
% https://q.uiver.app/?q=WzAsMyxbMCwxLCJcXHNwZWMgRF8wIl0sWzAsMCwiXFxzcGVjIEQiXSxbMSwxLCIoXFxzcGVjIEIsdCkiXSxbMCwyLCJcXGFscGhhXzAiLDJdLFsxLDIsIlxcYWxwaGEiXSxbMCwxLCJpIl1d
\[\begin{tikzcd}
	{\spec D} \\
	{\spec D_{_o}} & {(\spec B,t)}
	\arrow["{\alpha_{_o}}"', from=2-1, to=2-2]
	\arrow["\alpha", from=1-1, to=2-2]
	\arrow["i", from=2-1, to=1-1]
\end{tikzcd}\]
The immersion $j_{_B}\circ \tilde{s}_{_B}:\bh_{_B}\hookrightarrow \GL(n^2)_{_B}$ of group schemes on $\spec B$ is pulled back by $\alpha$, whereby one obtains the immersion 
\beqa
    j_{_D}\circ \tilde{s}_{_D}:\bh_{_D}\hookrightarrow \GL(n^2)_{_D}
    \label{6.2.2}
\eeqa
of smooth group schemes on $\spec D$.

By Cor \ref{stablereductionBD}, the $\GL(n^2)_{_R}$-torsor $P_{_R}$, together with its reduction of structure group to $\bh_{_R}$, is pulled back to $C_{_{D_{_o}}} = C \times \spec D_{_o}$ via via $(1\times \phi_{_o})$. Let $E_{_R}$ denote the $\bh_{_R}$-torsor on $C \times Z$ and let the afoerementioned pullback by $(1\times \phi_{_o})$ be denoted by by $E_{_{D_{_o}}}$. By \eqref{6.2.2}, one has the group scheme $\bh_{_D}$ on $\spec D$. Recall the theorem of Grothendieck \cite[Lemma 1, page 175]{desing}. On curves, the obstruction, which lies in   in $H^2(C, \operatorname{Ad} \bh_{_{D_{_o}}})$, vanishes, and the $\bh_{_{D_{_o}}}$-torsor $E_{_{D_{_o}}}$ lifts to an $\bh_{_D}$-torsor $E_{_D}$ on $\spec D$. 

To summarize, we have on $C_{_D}$ the $\GL(n^2)_{_D}$-torsor $P_{_D}$ together with a reduction of structure group $E_{_D}$ to the subgroup scheme $\bh_{_D}$. 

Now, by construction, the associated vector bundle $P_{_D}(n^2)$ gives a $D$-valued point of $N_{_n}$, say $q':\spec D\to N_{_n}$ such that the closed point maps to the image of $\xi$ under the projection $\fm_{_H}\to N_{_n}$. Our aim now is obtain a morphism $\spec D\to M_{_H}^{ss}$ such that the composites $\spec D\to M_{_H}^{ss}\to M_{_{\operatorname{SL}_{_n}}}^{ss}$ and $\spec D\to N_{_n}\to M_{_{\operatorname{SL}_{_n}}}^{ss}$ coincide. This entails a morphism into $\Gamma_{_\phi}$ \eqref{gammaphi}, and the normality of component $\fm_{_H}$ will give the desired lift to $(Z, \xi)$.

Let $K_{_D}$ denote the fraction field of $D$ and let $E_{K_{_D}}$ denote the restriction of the torsor to $C\times \spec K_{_D}$. Corollary \ref{cor} with $\mathbb k$ being $K_{_D}$, provides a semistable $H$-bundle $E_{K_{_D}}(H)$, which in fact comes as a reduction of structure group of a  semistable $\GL(n^2)$-bundle (and also from a semistable $\SL(n)$-bundle \eqref{cor}), and we thus have a morphism 
\beqa
\spec K_{_D}\to M_{_H}^{ss}\to M_{_{\SL(n)}}^{ss}.
\eeqa
The second statement of \eqref{cor} says that this map agrees with the composite $\spec K_{_D}\to N_n\to M_{_{\SL(n)}}^{ss}$. In other words, in the diagram \eqref{gammaphi} the map $\spec K_{_D}\to M_{_H}^{ss}$ lifts to a morphism $\spec K_{_D}\to \Gamma_{_\phi}$. The goal is to extend this map to all of $\spec D$. 

Since $M_{_H}^{ss}$ and hence $\Gamma_{_\phi}$ is proper, by the valuative criterion, the morphism $\spec K_{_D} \to\Gamma_{_\phi}$ extends to all prime ideals in $D$ of height $1$. Let $U \subset \spec D$ be the open subset containing all height $1$ primes. Thus, we have a diagram:
\beqa\label{affinestuff0}
\begin{tikzcd}
	& {\Gamma_{_\phi}} & {N_{_n}} \\
	U & {M_{_H}^{ss}} & {M_{_{\SL(n)}}^{ss}}
	\arrow[from=1-2, to=2-2]
	\arrow[from=2-2, to=2-3]
	\arrow[from=1-2, to=1-3]
	\arrow[from=1-3, to=2-3]
	\arrow[from=2-1, to=2-2]
	\arrow[from=2-1, to=1-2]
\end{tikzcd}
\eeqa
extending the morphism from $\spec K_{_D}$. We also have a base-change diagram:
\beqa\label{affinestuff}
\begin{tikzcd}
	{M_{_{H,D}}^{ss}} & {M_{_H}^{ss}} \\
	{\spec~D} & {M_{_{\SL(n)}}^{ss}}
	\arrow[from=1-1, to=2-1]
	\arrow[from=1-1, to=1-2]
	\arrow[from=1-2, to=2-2]
	\arrow[from=2-1, to=2-2]
\end{tikzcd}
\eeqa
By \cite{bapa}, the morphism $M_{_H}^{ss} \to M_{_{\SL(n)}}^{ss}$ induced by the spinor representation is an {\em affine} morphism. Hence, the vertical morphism $M_{_{H,D}}^{ss} \to \spec D$ is affine. By \eqref{affinestuff0}, the inclusion $U \hra \spec D$ lifts to $M_{_{H,D}}^{ss} $ and hence, by Hartogs' lemma, the morphism $U \to M_{_{H,D}}^{ss}$ extends to a morphism $\spec D \to M_{_{H,D}}^{ss}$. In other words, the arrow $U \to \Gamma_{_\phi}$ in \eqref{affinestuff0} extends to an arrow $\spec D \to \Gamma_{_\phi}$. This completes the proof.

\end{proof}

%{\em Completion of the proof of the first part of Theorem \ref{main}}:   

\subsubsection{Proof of the second part of Theorem \ref{main}}\label{secondpart} We need to complete the part which claims that we have a precise desingularisation, i.e., we have not blown-up any smooth point in $M^{ss}_{_H}$. It is known \cite[Corollary 3.4]{bh} that the locus of regularly stable bundles is precisely the smooth locus of $M^{ss}_{_H}$. 

Let $F_{_o}$ be a stable $H$-bundle such that $\Aut(F_{_o}) \subset Z(H)$, i.e. that $F_{_o}$ is {\em regularly stable}. Then we wish to show that this point in fact lies in $\fm_{_H}$, i.e., the smooth locus is not getting blown-up. In fact, we will prove that if $s \in \fm_{_H}$ is a closed point such that $p(s) = F_{_o}$, then the  the $H$-torsor $F_{_o}$ itself represents the point $s$ in the moduli space $\fm_{_H}$.

As in the discussion in \eqref{cor}, let $(R,\mathfrak m)$ be any complete dvr at the closed point $s \in \spec R$ with $K = \text{Fract}(R)$, and let $\bh_{_R}$ be the smooth group scheme and $E_{_R}$ a $\bh_{_R}$-torsor such that the point $s$ corresponds to the $\bh_{_o}$-torsor $E_{_o}$. 

%The group scheme $\bh_{_R}$ gives rise to a bounded subgroup $\bh_{_R}(R) \subset H_{_K}(K)$. 

%and hence, there is a maximal bounded subgroup $M$ sandwiched as $\bh_{_R}(R) \subset M \subset H_{_K}(K)$. Since $R$ is a complete discrete valuation ring and since $H$ is almost simple and simply connected, every maximal bounded subgroup $M \subset H_{_K}(K)$ is {\em schematic} and in fact, is a {\em maximal parahoric} in $H_{_K}(K)$. That is, there is a maximal parahoric group scheme $\bh'_{_R}$ with generic fibre $H_{_K}$ such that $\bh'_{_R}(R) = M$.

%Thus, we can  work for the moment with a maximal parahoric group scheme $\bh'_{_R}$ and view $E_{_o}$ as a $\bh'_{_o}$-torsor. 

As we have observed in \eqref{oldwine}, we get a morphism:
\beqa\label{morphatclosedfibre}
\gamma:\bh_{_o} \to H
\eeqa
and we have the identification $E_{_o}(H) \simeq F_{_o}$. Observe that if the closed fibre $\bh_{_o}$ is reductive, then by the rigidity of reductive groups, it will coincide with $H$ and this is reflected by the fact that the morphism \eqref{morphatclosedfibre} is an isomorphism. In other words, $\bh_{_R}$ would become {\em hyperspecial}.

{\em We claim that this is indeed the case}. Let us examine the morphism $\gamma$. Let $q:\bh_{_o} \to L_{_o}$ be the canonical Levi quotient map. Observe that the rank of the reductive quotient $L_{_o}$ is the same as the rank of $H$.  The image $\gamma (R_{_u}(\bh_{_o}))$ of the unipotent radical of $\bh_{_o}$ is either trivial, in which case, the map $\gamma$ factors as $\gamma \circ q:\bh_{_o} \to L_{_o} \hra H$. Else, $\gamma (R_{_u}(\bh_{_o}))$ is non-trivial and is in fact the unipotent radical of the image $\gamma(\bh_{_o}) \subset H$, which is therefore {\em non-reductive}. 

{\em We will deal with the second case first}. By the theorem of Borel-Tits \cite[Cor 3.9]{bort}, there is a canonically defined parabolic $P$ such that $\gamma(\bh_{_o}) \subset P \subset H$ and also $\gamma (R_{_u}(\bh_{_o})) \subset R_{_u}(P)$. Consider the extension of structure group of the $\bh_{_o}$-bundle $E_{_o}$ to $H$. Then the above discussion shows that this comes with a reduction of structure group $E_{_{o}}(P)$ to $P$. Furthermore, if $L = L_{_P}$ is the Levi quotient of $P$, then the associated $L$-bundle $E_{_{o}}(P)(L)$ has degree zero since it comes from $E_{_o}$. Since $P$ is a parabolic subgroup, it has a Levi splitting and so we get an inclusion $L \hra H$ and a morphism:
\beqa\label{vialevi}
\gamma(\bh_{_o}) \subset P \to L \hra H.
\eeqa
We now observe that the bundle $E_{_o}$ can be extended to an $H$-bundle via this new morphism. Let us call this bundle $E_{_o}(H)'$. We wish to compare the two $H$-bundles $E_{_o}(H)$ and  $E_{_o}(H)'$. Since $k$ is perfect,  there is a $1$-PS $\lambda: \mathbb G_{_m} \to H$ such that $P = P(\lambda) \subset H$. Furthermore, the Levi subgroup $L_{_P}$ can be realized as the centralizer of this $1$-PS. All in all, we can use this to conclude that the two principal $H$-bundles  $E_{_o}(H)$ and  $E_{_o}(H)'$ are $S$-equivalent in the moduli space of $H$-bundles, in the sense of Seshadri.   Since $F_{_o} = E_{_o}(H)$ is assumed to be stable, it follows that the new bundle $E_{_o}(H)'$ obtained by extending structure group via \eqref{vialevi} is isomorphic to $F_{_o}$. 

{\em The upshot is that in this case,  we deduce that the stable bundle $F_{_o}$ we started with gets a degree $0$ reduction of structure group to a Levi subgroup of a parabolic subgroup, i.e. to a {\em maximal rank subgroup of $H$}}.

In the first case, since the map $\gamma$ factors via $L_{_o}$, it again follows that we have a degree $0$ reduction of structure group to a {\em maximal rank subgroup} of $H$. Reformulating this discussion, we see that in either case, we have a reductive subgroup $L \subset H$ of maximal rank such that $F_{_o}$ has a reduction of structure group to $L$. Let us denote this reduction by  $E_{_L}$. Then  we have an isomorphism:
\beqa\label{vialevi}
F_{_o} \simeq E_{_L}(H).
\eeqa
%Since $\bh'_{_o}$ is the closed fibre of a maximal parahoric group scheme, the Levi quotient of the closed fibre {\em is not just reductive but is even semi-simple}. Furthermore, its root system is obtained by the Borel-de Siebenthal algorithm, namely by the omission of a simple root from the extended Dynkin diagram. Thus, 

Since $L$ is a connected reductive subgroup of $H$ of maximal rank,  by \cite{bdes} we have the identification: 
\beqa\label{locnontr}
L = Z_{_H}(Z(L))^{^\circ}.
\eeqa
Whence, if $\beta \in Z(L)$ is a non-trivial element (which exists by \eqref{locnontr}), then as an immediate consequence of \eqref{vialevi}, $\beta$ gives a non-central automorphism of $F_{_o}$. This contradiction completes the proof of the claim. It shows that the group scheme $\bh_{_R}$ was hyperspecial in the first place and $E_{_o} \simeq F_{_o}$, i.e. $F_{_o} \in {\fm}_{_H}$. This therefore completes the proof of Theorem \ref{main}.
\epr
%Since $\bh'_{_R}$  is a maximal parahoric group scheme, it is shown in \cite[page 662-663]{strongly} that at the closed fibre, the canonical Levi quotient map $q:\bh'_{_o} \to L'_{_o}$ has a splitting; $L'_{_o}$ is in fact realized as a suitable centralizer in the closed fibre. Thus, we get the composite $L'_{_o} \hra \bh'_{_o} \stackrel{\gamma}\to H$. 

%In other words, along with the morphism \eqref{morphatclosedfibre}, we get another morphism factors as $\gamma \circ q:\bh'_{_o} \to L'_{_o} \hra H$. 

\brem\label{comparewithses} It is checked easily enough that  when $H = \spin(3) \simeq \SL(2)$, we obtain an identification $\fm_{_H} = N_{_2}$. \erem
 
\section{The smooth compactification when $H$ is of type $\tt{D}_{_\ell}$} In this section we assume that $m = 2\ell$, so that $\text{dim}(\bbf) = 2\ell$. The Clifford algebra $\Cl_\tq$ associated to $\tq$ therefore has dimension $\dim(\Cl_\tq) = 2^{2\ell}$. The positive graded part 
\[
\Cl^{+}_\tq \subset \Cl_\tq
\]
therefore has dimension $2^{2\ell-1}$.  We shall also assume throughout that the quadratic form $\tq$ is of maximal Witt index, i.e., in the even case, of index $\ell$. 

In the even case, we have a decomposition $\bbf = N \oplus P$, where $N$ and $P$ are totally isotropic subspaces of dimension $\text{dim}(N) = \text{dim}(P) = \ell$.  This decomposition identifies the Clifford algebra $\Cl_\tq$ with $\text{End}(\Lambda N)$. We have a decomposition:
\beqa
\Lambda N = \Lambda ^{^{even}} N \oplus \Lambda ^{^{odd}} N
\eeqa
and the positive graded part $\Cl^{+}_\tq$ is compatible with this splitting. Thus, we get an identification:
\beqa
\Cl^{+}_\tq \simeq \text{End}(\Lambda ^{^{even}} N) \oplus \text{End}(\Lambda ^{^{odd}} N),
\eeqa
i.e., as a direct sum of two copies of the matrix algebra $\cm_{_{2^{^{\ell-1}}}}(k)$.

Let $S^{^{+}} := \Lambda^{^{even}} N$ and $S^{^{-}} := \Lambda^{^{odd}} N$, and let $S := S^{^{+}} \oplus S^{^{-}}$, where $\text{dim}(S) = 2^{^\ell}$.
Then we have the {\em half-spin representations}  denoted by:
\beqa
{\Large\text{\cursive s}}_{_{\pm}}:\spin(\tq) \to \SL(S^{^{\pm}})
\eeqa
which are {\em irreducible representations}; if $m = 2\ell$ and $\ell$ is odd  then ${\Large\text{\cursive s}}_{_{\pm}}$ {\em are both faithful and irreducible}, else ${\Large\text{\cursive s}}_{_{\pm}}$ {\em are not faithful}.  The representation:
\beqa
{\Large\text{\cursive s}}:\spin(\tq) \to \GL(S)
\eeqa
given by ${\Large\text{\cursive s}} := {\Large\text{\cursive s}}_{_+} \oplus {\Large\text{\cursive s}}_{_-}$ is however {\em faithful}. 
Let $n = \text{dim}(S) =  2^{^\ell}$ so that ${n \over 2} = \text{dim}(S^{^{\pm}}) = 2^{^{\ell-1}}$.

\brem In the case when $m = 2\ell$, and $\ell$ is {\em odd}, we could work with either of the half-spin representations and carry out the procedure that has been executed for the case when $m$ itself was odd. The construction of the smooth compactification for $M_{_H}^s$ and hence a desingularisation of $M_{_H}^{ss}$, in this case now follows. \erem

For this reason, unless otherwise mentioned, we will henceforth assume that {\em $s$ is even}. As is shown in \cite[\S3.6, III.6.1, page 165]{chevalley}, the kernels of the half-spin representations ${\Large\text{\cursive s}}_{_{\pm}}$ are given by the subgroups $\{1, w\}$ and $\{1,-w\}$, where $w$ is the involution lying in the center of $\Cl_\tq$ which anti-commutes with every element of $\bbf \subset \Cl_\tq$. Nevertheless, the images
$\spin(m)/\text{ker}({\Large\text{\cursive s}}_{_{\pm}})$ are isomorphic and the resulting group is called the {\em half-spin group}, and denoted by $\hspin(m)$. Moreover, the representations ${\Large\text{\cursive s}}_{_{\pm}}$ become equivalent on $\hspin(m)$ and we denote this {\em faithful, irreducible} representation by:
\beqa
{\Large\text{\cursive s}}_{_h}:\hspin(m) \to \SL(S^{^{+}}).
\eeqa
The spin and half-spin representations have some nice behaviour with respect to restrictions. We faithfully follow \cite{chevalley} for this.

%From $\spin(2s)$ to $\spin(2s +1)$: Consider the embedding of $k^{^{2s}} = N \oplus P$ in $k^{^{2s +1}} = N \oplus P \oplus U$, where $U$ is a line orthogonal to $W$ and $W'$. This induces a corresponding embedding of the Clifford algebras 

From \underline{$\spin(2\ell)$ to $\spin(2\ell -1)$}: let us denote by $\overline{\bbf}$ a vector space of dimension $2\ell -1$ over $k$ and let $\overline{\tq}$ be a quadratic form on $\overline{\bbf}$ which is non-degenerate and of maximal Witt index $\ell-1$. We denote the full Clifford algebra $\Cl_{\overline{\tq}}$ of $\overline{\tq}$ by $\overline{\Cl}$,  by $\overline{\Cl}^{+}$ the even graded part of  $\overline{\Cl}$, by $\overline{{\Large\text{\cursive s}}}:\spin(2\ell -1) \to \SL(\overline{S})$, the spin representation. Let $N',P'$ be two maximal totally singular subspaces of $\overline{\bbf}$ such that the sum $N' + P'$ is a direct sum and a non-isotropic subspace of $\overline{\bbf}$. Choose a vector $v \in \overline{\bbf}$ orthogonal to $N' \oplus P'$, and set 
\beqa
a := \overline{\tq}(v).
\eeqa
We embed $\overline{\bbf}$ in $\bbf$ (which is of dimension $2\ell$) and which comes with the decomposition $\overline{\bbf} \oplus \langle v' \rangle$; we extend the quadratic form $\overline{\tq}$ to a quadratic form $\tq$ on $\bbf$ by setting
\beqa\label{tqfrombar}
\tq(\overline{x} + c.v') := \overline{\tq}(\overline{x}) - a.c^2,~~~ c \in k
\eeqa
if $\overline{x} \in \overline{\bbf}$. Then $\tq$ is of rank $2\ell$ and of index $\ell$, for $\tq(v + v') = 0$ and $v + v'$ is orthogonal to every element of $N'$ (relative to $\tq$). This shows that $N = N' + \langle v + v' \rangle $ is totally singular for $\tq$. Also, $P = P' + \langle v - v' \rangle$ is totally singular for $\tq$ and $\bbf = N \oplus P$.

Either one of the half-spin representations ${\Large\text{\cursive s}}_{_{\pm}}:\spin(2\ell) \to \SL(S^{^{\pm}})$ induces the  spin representation $\overline{{\Large\text{\cursive s}}}:\spin(2\ell + 1) \to \SL(\overline{S})$.

By applying the canonical involution on the Clifford algebra  we can identify the Spin group as a subgroup of the group of units of the opposite algebra \eqref{opposite}, and as in   \eqref{oppositeunits}, we get
\beqa\label{spin1}
{\Large\text{\cursive s}}:\spin(\tq) \hra \cu_\tq.
\eeqa
This inclusion gives rise to analogues ${\Large\text{\cursive s}}_{_{\pm}}$, and by an abuse of notation, we will henceforth call these the half-spin representations. 

\bcor\label{overqeven} 
Let $\overline{\cq} := \cq(\overline{\bbf})$, the affine space of quadratic forms on $\overline{\bbf}$. Let $t \in \overline{\cq}$, $B = \hat{\co}_{_{\overline{\cq},t}}$, and $L = \text{Fract}(B)$.  
\begin{enumerate}
\item Then we have a {\em smooth quasi-affine} group scheme ${\bh_{_{B}}^{^{\spin}}}$ on $\spec B$ which {\em extends} $\hspin_{_L}$ on $\spec L$. Further, we have an inclusion of smooth group schemes:
\beqa\label{hbt1.1}
{\bh_{_{B}}^{^{\spin}}} \stackrel{{\Large\text{\cursive s}}_{_{h,B}}}\hra \cu_{_{B}} \stackrel{j_{_B}}\hra \text{GL}(S^{^{+}})_{_{B}}.
\eeqa
where $j_{_B}$ is a closed embedding while ${{\Large\text{\cursive s}}_{_{h,B}}}$ is not.
\item Let $\tq_{_L}$ be the split form as in \eqref{tqfrombar}. Then the faithful representation $\spin(\tq_{_L}) \hra \GL(S^{^{+}} \oplus S^{^{-}})_{_L}$ extends to give a quasi-affine subgroup scheme $\bh_{_B} \hra  \text{GL}(S)_{_{B}}$ together with inclusions:
\beqa\label{hbt3}
\bh_{_B} \hra {\bh_{_{B}}^{^{\spin}}}\times {\bh_{_{B}}^{^{\spin}}} \hra  \text{GL}(S)_{_{B}}
\eeqa
\end{enumerate}
\ecor
\bpr The existence of the group schemes and their inclusions in $\text{GL}(S)_{_{B}}$ follow verbatim from Theorem \ref{hbt}. The only conclusion which requires a bit more work is the last statement, namely that we have a canonical morphisms ${\Large\text{\cursive s}}_{_{\pm, B}}:\bh_{_B} \to {\bh_{_{B}}^{^{\spin}}}$, extending the half-spin maps ${\Large\text{\cursive s}}_{_{\pm, L}}: \spin(\tq_{_L}) \to \hspin(\tq_{_L})$. This map is well-defined upto height $1$ primes since it comes by taking schematic closures of $\bh_{_L} \hra  \text{GL}(S)_{_{L}}$ which factor through the images of $\hspin(\tq_{_L}) \times \hspin(\tq_{_L})$. Once, height $1$ primes are handled, then we appeal to the extension theorem due to Raynaud \cite[IX, Corollary 1.4]{raynaud} to get the morphism $\bh_{_B} \hra {\bh_{_{B}}^{^{\spin}}} \times {\bh_{_{B}}^{^{\spin}}}$. 

\epr

\subsubsection{The construction of the space $\fm_{_H}$}

As in \S6, we again take a stable $H$-bundle $E$ with full holonomy and by extending structure group using the irreducible homomorphisms ${\Large\text{\cursive s}}_{_{\pm}}$, we obtain two stable $\hspin(m)$-bundles  and corresponding stable $\SL(S^{^{\pm}})$-bundles $E(S^{^{\pm}})$. 

Note that the groups $\SL(S^{^{\pm}})$ are the special linear groups contained in the group of units of the matrix algebras $\cm^{^{\pm}} \simeq \cm_{_{{n \over 2}}}(k)$, with ${n \over 2} = 2^{^{\ell-1}}$. By the discussion above, we can identify $\cm^{^{\pm}}$ with the even graded Clifford algebra $\overline{\Cl}^{+}$ associated to the quadratic form $\overline{\tq}$ on $\overline{\bbf}$. The Clifford algebra $\overline{\Cl}^{+}$ is isomorphic to  $\cm_{_{{n \over 2}}}(k)$ and therefore, we view  $\GL({n \over 2}))$ as the group of units of the opposite algebra $\overline{\Cl}^{+\circ}$.

Thus, ${\Large\text{\cursive s}}_{_{\pm}}$ induces morphisms ${\Large\text{\cursive s}}_{_{\pm,*}}: M_{_{\tt f}}(H) \to M^s(\SL({n \over 2}))$. As in the odd case, we get a pair of points $(W_{_+}, W_{_-}) \in N_{_{n \over 2}}^s \times N_{_{n \over 2}}^s$. We can allow independent specializations of the matrix algebra $\cm^{^{\pm}}$ and follow the procedure in \S6.

 We also have the open immersion ${\Large\text{\cursive v}}:M_{_{\tt f}}(H) \hra M^{ss}_{_H}$. Thus, we get a morphism $\phi = ( {\Large\text{\cursive v}},~{\Large\text{\cursive s}}_{_{+,*}}, {\Large\text{\cursive s}}_{_{-,*}}): M_{_{\tt f}}(H) \to M^{ss}_{_H} \times N_{_{n \over 2}} \times N_{_{n \over 2}}$. 
 
 Let $\Gamma_{_\phi}$ denote the base-change (with its reduced induced structure) given by:
\beqa\label{gammaphitypeD}
\begin{tikzcd}
	{\Gamma_{_\phi}} & {N_{_{n \over 2}} \times N_{_{n \over 2}}} \\
	{M^{ss}_{_H}} & {M^{ss}_{_{\SL({n \over 2})}} \times M^{ss}_{_{\SL({n \over 2})}}}
	\arrow[from=1-1, to=2-1]
	\arrow["{\Large\text{\cursive s}}_{_*}", from=2-1, to=2-2]
	\arrow[from=1-1, to=1-2]
	\arrow[from=1-2, to=2-2]
\end{tikzcd}\eeqa
Let
\beqa\label{thecandidateD}
\fm_{_H} := component~of~the~normalization~of~\Gamma_{_\phi}~containing~M_{_{\tt f}}(H).
\eeqa
As in \eqref{thecandidate}, we get the following diagram, with projections $p$ and $q$. 
\beqa\label{thecorrespeven}
\begin{tikzcd}
	& \fm_{_H} \\
	M^{ss}_{_H} && N_{_{n \over 2}} \times N_{_{n \over 2}}
	\arrow["{^p}", from=1-2, to=2-1]
	\arrow["{^q}"', from=1-2, to=2-3]
\end{tikzcd}
\eeqa
The specializations of the individual matrix algebras are controlled by the two independent degenerations of the non-degenerate quadratic form $\overline{\tq}$ on $\overline{\bbf}$. 

Let $\xi \in \fm_{_H}$ and let $(Z,\xi)$ be a formal neighbourhood of $\xi$. As in \eqref{step1}, let $\vartheta:(Z,\xi) \to (\spec B,t)$ be the map induced by the "versality", with $B$ as in \eqref{overqeven}. By pulling back, using $\vartheta$, the diagram \eqref{hbt3}, we get:
\beqa\label{hbt4}
\bh_{_Z} \hra {\bh_{_{Z}}^{^{\spin}}} \times {\bh_{_{Z}}^{^{\spin}}} \hra  \text{GL}(S)_{_{Z}}
\eeqa

To complete the proof of the theorem analogous to Theorem \ref{main}, we need to prove \eqref{stablereduction}. This involves two steps. Let us label the pull-backs under $(1 \times {\Large\text{\cursive s}}_{_{\pm,*}})$ of the universal familiy $P_{_{n \over 2}}$ on $C \times N_{_{n \over 2}}$ by $P_{_{n \over 2}}^{^{\pm}}$. These are rank $({n \over 2})^{^2}$-bundles whose transition functions lie in the group of units of the corresponding degenerate Clifford algebra. 

Following the proof of \eqref{stablereduction}, we derive the conclusion that the structure groups of $P_{_{n \over 2}}^{^{\pm}}$ can be reduced firstly to the half-spin group scheme  ${\bh_{_{Z}}^{^{\spin}}}$. Then, by the same set of arguments, it can be further reduced to the subgroup scheme  $\bh_{_Z}$ since the generic bundles have a reduction of structure group to $\spin(\tq_{_L})$. All else goes as before and we summarise this in the following:
\bth\label{maineven} The conclusion of Theorem \ref{main} holds when $H$ is of type ${\tt D}$.\eeth

\section{The exceptional groups ${\tt G}_{_2}$ and ${\tt F}_{_4}$ and E.Dynkin's "plethysm"} 
The following important results of Dynkin \cite[Theorem 6.2, page 308, and Table 16, page 372]{dynkin} helps us bring out a connection between apparently unrelated phenomena! 

\begin{enumerate}\label{dynkin}
\item (\footnote{We found (\eqref{dynkin} (a)) independently, and verified it using the tables for branching rules. This operation is a plethysm of sorts. Searching for similar features for other exceptional groups, we discovered Dynkin's general theorem!})
 Let $j:{\tt G}_{_2} \hra \spin(\mathfrak g_{_2})$ be the adjoint representation. The adjoint representation lifts to the Spin group by simple connectedness of the group $H$ and the Killing form on $\mathfrak g_{_2}$ is the non-degenerate quadratic form. Then, the half-spin representations ${\Large\text{\cursive s}}_{_{\pm}}:\spin(\mathfrak g_{_2}) \hra \SL(S^{^{\pm}})$ restricted to ${\tt G}_{_2}$ are equivalent to the standard irreducible representation $\pi_{_\rho}:{\tt G}_{_2}\to \SL(V_{_\rho})$, where $\rho$ is the half-sum of positive roots. Note that $\text{dim}(V_{_\rho}) = 2^{^6}$ 
\item Let $j:{\tt F}_{_4} \hra \spin(26)$ {\em be the irreducible representation fundamental of dimension $26$}. Then, the half-spin representations ${\Large\text{\cursive s}}_{_{\pm}}:\spin(26) \hra \SL(S^{^{\pm}})$ restricted to ${\tt F}_{_4}$ are equivalent to the standard irreducible representation $\pi_{_\rho}:{\tt F}_{_4} \to \SL(V_{_\rho})$, where $\rho$ is the half-sum of positive roots. Note that $\text{dim}(V_{_\rho}) = 2^{^{6}}$.
\item Let $j:{\tt C}_{_3} \hra \spin(14)$ {\em be the irreducible representation fundamental of dimension $14$}. Then, the half-spin representations ${\Large\text{\cursive s}}_{_{\pm}}:\spin(14) \hra \SL(S^{^{\pm}})$ restricted to ${\tt C}_{_3} $ are equivalent to the standard irreducible representation $\pi_{_\rho}:{\tt C}_{_3}  \to \SL(V_{_\rho})$, where $\rho$ is the half-sum of positive roots. Note that $\text{dim}(V_{_\rho}) = 2^{^{12}}$.
\end{enumerate}
In this section let $H$ stand for any of these three types of groups listed above. Let us fix the representations $H \stackrel{\phi}\hra \spin(2\ell) \stackrel{{\Large\text{\cursive s}}} \hra \SL(S)$ as in \eqref{dynkin}, i.e., $\ell = 7~or~13$. The composite ${\Large\text{\cursive s}} \circ \phi$ is identified with the fundamental representation $\pi_{_\rho}$ of $H$. 

Let $E$ be a stable $H$ bundle on the curve $C$ with full holonomy, i.e., the group $H$ coincides with the Zariski closure of the representation of the fundamental group of $C$ which gives $E$. Since the genus $g \geq 2$, such bundles always  exists and give a smooth open locus of the moduli space $M_{_{H}}^{s}$ of stable $H$-bundles. By the arrows above, we deduce the fact that the associated bundle $\pi_{_{\rho,*}}(E)$ is a {\em stable} $\SL(S)$-bundle. Furthermore, this stable bundle is isomorphic to ${\Large\text{\cursive s}}_{_*} (\phi_{_*}(E))$. Since $\phi_{_*}(E)$ is a stable $\spin(26)$-bundle and  $26 = 2 \times 13$, and $13$ being odd, the half-spin representations are both faithful and irreducible. Hence, by the earlier construction in \S6, the bundle ${\Large\text{\cursive s}}_{_*} (\phi_{_*}(E))$ gives a point in the moduli space $N_{_n}^s$, for $n = 2^{^{\ell-1}}$. Therefore, we obtain a morphism:
\beqa
\pi_{_*}:M_{_{\tt f}}(H) \to N_{_n}.
\eeqa
As before \eqref{thecandidate}, \eqref{thecandidateD}, we have the open immersion ${\Large\text{\cursive v}}:M_{_{\tt f}}(H) \hra M_{_{\fg}}^{ss}$ and  the scheme $\fm_{_H}$ together with the two projections:
\beqa\label{thecorrespeexcep}
\begin{tikzcd}
	& {\fm_{_{H}}} \\
	{M_{_{H}}^{ss}} && {N_{_n}}
	\arrow["{^p}", from=1-2, to=2-1]
	\arrow["{^q}"', from=1-2, to=2-3]
\end{tikzcd}
\eeqa
The "versal" family is again the smooth quasi-affine group scheme ${\bh}_{_{\spec B}}$ on $\spec B \subset \cq$, with generic fibre $H_{_L}$ (obtained from Theorem \ref{hbt}). This now plays the role of the group scheme $\bh_{_B}$ in \eqref{step1}.  
\bprop\label{stablereduction2} With the hypothesis as in \eqref{stablereduction1}, let further $H_{_K}$ be split of type ${\tt G}_{_2}$ , ${\tt F}_{_4}$ or ${\tt C}_{_3}$. Then the conclusion of \eqref{stablereduction1} holds good.\eprop
\bpr With the notations as in \eqref{stablereduction1}, in the cases in discussion in this proposition, it is known \cite{biswas-holla}, that the Brauer group $\text{Br}(Z_{_L}) = 0$. It follows that $\epsilon^{^*}(E)$ in fact is a principal $\cu_{_{\tq^{^{\tt h}}}}$-torsor, i.e. for the split group scheme. Further, this gets a reduction of structure group to $H_{_L}$ in all these cases since the universal torsor already exists on $C \times Z_{_L}$ in these case. This proves the proposition. \epr

Modulo the previous proposition, the rest of the proof goes through all the steps of Theorem \ref{main} and we get the following:
\bth\label{mainexceptional} Let $H$ be of type ${\tt G}_{_2}$, ${\tt F}_{_4}$ or ${\tt C}_{_3}$. The scheme $\fm_{_H}$ is irreducible, smooth and projective, and the morphism $p:\fm_{_H} \to M_{_{H}}^{ss}$  is an {\em isomorphism} over the open smooth locus $M_{_{H}}^{rs} \subset M_{_{H}}^{s}$ consisting of regularly stable principal $H$-bundles (i.e., with trivial automorphisms).  \eeth
%\brem The Table 16 of \cite{dynkin} also gives us one more interesting case, namely $\fg = {\tt C}_{_3}$ for which we obtain a theorem analogous to Theorem \ref{mainexceptional}. However, the third centerless exceptional group ${\tt E}_{_8}$ is unfortunately not in the list! \erem
\section{On the characteristic of the base field $k$}\label{charp}
In this section $\text{char}(k) = p > 0$. I will briefly indicate the bound on $p$ needed for the main theorems to hold, including the ones in the Appendix for the local models. There are two places where the bounds will be needed. Firstly in \S4, where such bounds will be needed for exponentiation of unipotent subgroups. This follows by \cite{bdp} and the Coxeter number of $H$ is good enough as a lower bound for $p$.

The second and more delicate place is prime bounds imposed by semistability of associated bundles. Let $\rho:H \hra \GL(V)$ be the following list of representations: ${\tt B}_{_\ell}, {\tt D}_{_\ell}, {\tt G}_{_2}, {\tt F}_{_4}~~or~~{\tt C}_{_3}$
\begin{itemize}
\item When $H$ is of type ${\tt B}_{_\ell}$, $\rho$ is simply the {\em spinor} representation $\tc{s}$ \eqref{thespinorrep}.
\item When $H$ is of type ${\tt D}_{_\ell}$, $\rho$ will be the representations $\tc{s}_{_\pm}$ (defined after \eqref{spin1}).
\item When $H$ is of type ${\tt G}_{_2}, {\tt F}_{_4}~~or~~{\tt C}_{_3}$, then the representation $\rho$ is simply $\tc{s}_{_\pm} \circ j$ with notations as in  \eqref{dynkin}.

\end{itemize}
With this notation, let $p > \text{ht}_{_H}(\rho)$, where  $\text{ht}_{_H}(\rho)$ is the Dynkin height of the representation $\rho$ (\cite[Definition 4.4]{bdp}. We note that this bound on $p$ suffices for the existence of a desingularisation of $M_{_{H}}^{ss}$.
\section{The moduli functor} Recall the functor $F: \text{Sch/k} \to \text{Sets}$, where $F(T)$ is the set of isomorphism classes of {\em rigidified families} of vector bundle $(V,\Delta)$ of rank $n^2$ on $C \times T$, such that the $T$-algebra $p_{_*}({\mathcal E}nd(V))$ is a {\em specialization} of the matrix algebra (\cite[Definition 2, page. ]{desing}). It is shown in \cite{desing} that the moduli space $N_{_n}$ represents $F$. It is known that $N_{_n}$ is a normal projective variety. The construction in Part II extends without much difficulty to show that the existence of group scheme $\bh_{_N} \subset \text{GL}(n^2)_{_N}$ on $N_{_n}$ with the defining properties and inclusions as in \eqref{}. The inclusion is not a {\em closed embedding} and hence in general the quotient $\text{GL}(n^2)_{_N}/\bh_{_N}$ need not be representable.  

For a point $(V,\Delta)$ in $F(T)$, let $E_{_V}$ denote the underlying frame bundle with structure group $\text{GL}(n^2)_{_T}$.

 \beqa
F_{_H}:\text{Sch/k} \to \text{Sets}
\eeqa
be defined as $F_{_H}(T) := \text{Isom classes of pairs $(E_{_V}, \sigma)$}$, where $\sigma$ is a reduction of structure group of $E_{_V}$ to $\bh_{_T} \subset \text{GL}(n^2)_{_T}$. Notice that, by the representability of $F$ by the scheme $N_{_n}$, we have a structure morphism $T \to N_{_n}$ and the group scheme $\bh_{_T}$ is the pull-back. The following theorem is proven along the same lines as in \cite{desing}.
\bth The scheme $\fm_{_H}$ represents $F_{_H}$. \eeth

\part{\sc Flat degenerations of semi-simple group schemes}
\section{The construction of the versal group schemes}
 This section is of a general nature, so we have the following notations and conventions. Let $Y = \spec(B)$, where $B$ is a regular local domain with residue field $k$ (which is algebraically closed) and with fraction field $L$. Let $\mathcal A_{_B}$ be a unital associative algebra over $B$ which is locally free as a $B$-module of dimension $n^2$. Further, let $\mathcal A_{_L}$ be an Azumaya algebra. Suppose further that we have an inclusion of algebras
 \beqa
 \mathcal A_{_B} \hra \text{End}(\mathcal A_{_B})
 \eeqa
 where $\mathcal A_{_B}$ is viewed as a projective $B$-module.  This inclusion induces an inclusion of the group scheme of units:
 \beqa
 j_{_B}: \mathcal U_{_B} \hra \text{Aut}(\mathcal A_{_B}).
 \eeqa
 
 Let $H = H_{_L}$ be a  connected, semisimple algebraic group over $L$ such that:
 \beqa
 H_{_L} \hra \mathcal U_{_L} \hra \text{GL}(n^2)_{_L}
 \eeqa 
is a closed embedding over $L$.  Let $T$ be a {\em maximal $L$-split torus} of $H_{_L}$. Let $\bf \Phi$ be the system of roots of $H_{_L}$ relative to $T$. For $a \in \bf \Phi$, let $U_{_a} \subset H$ be the root groups (which are normalised by $T$ and such that the weights of $T$ in $U_{_a}$  belong to $a$, and let $\bf \Phi^+$ (resp. $\bf \Phi^-$) be the subset of positive (resp. negative) roots for the choice of an order on the character group $X^{*}(T)$. 
 
The question we wish to address in this part of the paper is a very general one. {\it ``Does there exist a smooth quasi-affine group scheme $\mathcal H_{_B}$ over $\spec B$ together with a diagram of inclusions:
\beqa
H_{_B} \hra \mathcal U_{_B} \stackrel{j_{_B}}\hra \text{Aut}(\mathcal A_{_B}).
 \eeqa
 such that its fibre over $\spec L$ is isomorphic to $H_{_L}$? Can one describe the fibres of this extended group scheme over the closed points of $\spec B$? How unique or canonical is this extension?} We note that the inclusion $j_{_B}$ is a closed embedding but the first inclusion which extends the closed embedding over $L$ is not necessarily a closed embedding and hence gives, in general, only a quasi-affine group scheme $H_{_B}$.
 
We answer this question in the cases when the Dynkin type of $H$ is $\tt B_{_n}, \tt D_{_n}, \tt G_{_2}, \tt F_{_4}, \tt C_{_3}$. This list is central to the many body of the paper. We observe that for Dynkin type $\tt A_{_n}$, {\em we have assumed} that such an extension already exists over $B$, by virtue of the assumption that the algebra $\mathcal A_{_B}$ is an Azumaya algebra over $L$. In the situations where this will get applied, such an assumption will be shown to hold.

This answer has general interest in the question of Bruhat-Tits group schemes over higher dimensional regular rings, especially since in the case of $\tt B_{_n}, \tt D_{_n}$, we encounter {\em non-split} group $H_{_L}$ and this therefore goes out of the considerations of Bruhat-Tits theory even over discrete valuation rings. 
 
\subsection{When $H$ is split} Let $H_{_L}$ be {\em split}, in which case,  $T$ is assumed to be a split maximal $L$-torus and hence the root groups $U_{_a}$ are all $1$-dimensional.  The product $\prod U_{_a} \times T$ in $H$ is an open immersion and its image is the so-called {\em big cell} which we denote by $\fb$. We have the following lemma:
\blem\label{lieextension} When $H_{_L}$ is split,  the Lie algebra $\lie(H_{_L})$ of $H_{_L}$ extends as a Lie algebra bundle $\mathfrak R$ to the whole of $\spec B$.\elem
\bpr Consider the Lie algebra $\lie(H_{_L})$ of $H_{_L}$. The tangent space of the big cell $\fb_{_L} \subset H_{_L}$ gives the root space decomposition of the Lie algebra $\lie(H_{_L})$ as a direct sum of the Lie algebras of the root groups $U_{_a}$  and maximal torus $T$. The extension ${\fb}_{_{Y'}}$ of the big cell to the open subset $Y'$ therefore gives a decomposition
\beqa\label{bigcell 2}
\text{Lie}(H_{_{Y'}}) = \text{Lie}(\mathfrak U_{_{Y'}}^{^+}) \oplus \text{Lie}(\mathfrak T_{_{Y'}}) \oplus \text{Lie}(\mathfrak U_{_{Y'}}^{^-}) .
\eeqa
of $\lie(H_{_{Y'}})$  {\em as a direct sum of line bundles}. By taking double duals $\lie(H_{_{Y'}})^{{\ast\ast}}$, we get a reflexive sheaf on $\spec B$. {\em Since we are over a regular ring}, reflexive sheaves of rank $1$ are locally free (by \cite[Proposition 1.9]{hartshorne}). The bundle  $\lie(H_{_{Y'}})$ is a direct sum of line bundles, therefore, $\lie(H_{_{Y'}})^{{**}}$ is a locally free $B$-module. Because $\text{codim}(Y \setminus Y') \geq 2$, we can, by an application of Hartogs lemma, conclude that the Lie algebra structure on $\lie(H_{_{Y'}})$ also extends to  $\lie(H_{_{Y'}})^{{**}}$. We will simply denote this Lie algebra bundle over $\spec B$ by $\mathfrak R$. \epr

\subsubsection{The Group scheme extension} We now make assumptions on $H_{_L}$ in accordance with the previous discussion. As noted above, for the present, we have restricted ourselves to the {\em split case}. 

Let $W$ be a finite dimensional $k$-vector space and let $H_{_L} \subset \GL(W \otimes L)$ be a closed embedding of  $H_{_L}$ and let $\fw := W \otimes B$.  In the case we are in, we can take $W_{_B}$ to be the projective $B$-module $\mathcal A_{_B}$.  

Then the schematic closures $\mathfrak U_{_a}$ of the root groups $U_{_a}$ and the schematic closure $\mathfrak T$ of the maximal torus $T$ in $\text{Aut}(\fw) :=\mathfrak{GL}({\fw})$, are $B$-schemes which in general are not {\em flat} since $\text{dim}(B) \geq 1$. Let $\mathfrak P$ be the set of primes of height $1$ in $B$ which defines the complement of $\spec L$ in $\spec B$. By taking the schematic closure of $H_{_L}$ in $\mathfrak {GL}(W \otimes B_{_{\mathfrak p}})$ at each of the $\mathfrak p \in \mathfrak P$, we get a closed subgroup scheme  $H_{_{Y'}} \subset \mathfrak {GL}({\fw})_{_{Y'}}$, where $Y'$ is an open subset of $\spec B$ which contains points of height $1$. This group scheme is therefore {\em canonical}, determined completely by the $Y'$-lattice $\fw_{_{Y'}} \subset W_{_L}$.

Note that the schematic closures at the height $1$ primes of both, the root groups $U_{_a}$ and the maximal torus $T$, are {\em flat} and hence we get {\em smooth} subgroup schemes $\mathfrak U_{_{\bf \Phi^-,Y'}}$, $\mathfrak U_{_{\bf \Phi^+,Y'}}$ and $\mathfrak T_{_{Y'}}$ of $H_{_{Y'}}$ (we remind ourselves that the characteristic of the base field $k$ is either zero or large, which ensures smoothness). 

Firstly, the group schemes  $\mathfrak U_{_{\bf \Phi^-,Y'}}$, $\mathfrak U_{_{\bf \Phi^+,Y'}}$ are flat {\em unipotent} closed sub-group schemes of $H_{_{Y'}}$. Secondly, since $T$ is a split torus, by \cite[Example 5, p.291]{blr}, it follows that the extension $\mathfrak T_{_{Y'}}$ of $T$ is a multiplicative subgroup scheme of $H_{_{Y'}}$ on  the whole of $Y'$. Let  
\beqa\label{bigcell1}
{\fb}_{_{Y'}} := \mathfrak U_{_{\bf \Phi^-,Y'}} \times \mathfrak T_{_{Y'}} \times \mathfrak U_{_{\bf \Phi^+,Y'}}
\eeqa
denote the so-called {\em big cell} of the group scheme $H_{_{Y'}}$. Note that this gives a neighbourhood of the identity of the group scheme $H_{_{Y'}}$. We denote by $\mathfrak U^{^{\pm}}$ the group schemes $\mathfrak U_{_{\bf \Phi^{\pm}}}$. 
 
We now make a few important observations. Firstly, since $\lie(H_{_L})$ is a Lie sub-algebra of $\lie(\GL(W \otimes L))$, the roots of $H_{_L}$ with respect to $T_{_L}$ are the roots of $\GL(W \otimes L)$ with respect to $T_{_L}$ and hence $\prod_{_{a \in \bf \Phi}} U_{_a}$ is a closed sub-scheme of the corresponding product of root groups in $\GL(W \otimes L)$. We have the following basic observation.

Let us denote by $\mathfrak Y$ the schematic closure of $H_{_L}$ in $\mathfrak{GL}({\fw})$. We note that $\mathfrak Y$ need not be a group scheme nor does it need to be flat in general. Let  $\mathfrak U_{_a}$ for $a \in \bf \Phi^{\pm}$ denote the schematic closure of the root groups $U_{_a}$ and $\mathfrak T$ the schematic closure of $T_{_L}$, namely the fixed  maximal split $L$-torus. {\em Note that the root groups are necessarily $1$-dimensional since we are working with the root system coming from  the split maximal torus $T_{_L}$}.
\bprop\label{one} When the characteristic of $k$ is zero (or large prime characteristics), the $B$-schemes $\mathfrak U_{_a}$ for $a \in \bf \Phi^{\pm}$, and $\mathfrak T$ are flat (and hence smooth) subgroup schemes of $\mathfrak{GL}({\fw})$. Hence, the big cell ${\fb}_{_{Y'}} \subset H_{_{Y'}}$ extends as an open subscheme  $\fb$ of the schematic closure $\mathfrak Y$. \eprop

\bpr We work on each factor in the big cell \eqref{bigcell1}.  By our assumption, $T_{_L}$ is a maximal split $L$-torus and canonically extends to $Y'$. Since the complement of $Y'$ in $Y$ has codimension   $\geq 2$,  by \cite[Proposition  IX, 2.4]{raynaud}, the group scheme inclusion $\mathfrak T_{_{Y'}} \subset H_{_{Y'}}  \subset \mathfrak{GL}({\fw})_{_{Y'}}$ extends canonically to a smooth group scheme $\mathfrak T'$ together with an closed embedding $\mathfrak T' \hra \mathfrak{GL}({\fw})$, i.e., as a sub-torus of the group scheme $\mathfrak{GL}({\fw})$  on the whole of $Y$.   Since $\mathfrak T'$ is {\em flat and hence torsion-free}, we deduce that $\mathfrak T'$ is indeed the schematic closure $\mathfrak T$ of  $\mathfrak T_{_{Y'}}$ in $\mathfrak Y$ and hence in $\mathfrak{GL}({\fw})$ (see \cite[1.2.6]{bt2}).  In particular $\mathfrak T$ is {\em flat}.

Consider the decomposition \eqref{bigcell 2}.  As we have seen, the Lie algebra $\text{Lie}(H_{_L})$ extends as a locally free Lie algebra bundle, $\mathfrak R$  to the whole of $\spec B$, we may identify $\mathfrak R$  as the reflexive closure of $\text{Lie}(H_{_Y'})$. 

%In particular, we see that $\text{Lie}(\mathfrak T)$ is the reflexive closure of $\text{Lie}(T_{_L})$ in $\mathfrak R$. 

Observe that $\text{Lie}(\mathfrak U_{_{Y'}}^{^\pm})$ is a direct summand of $\text{Lie}(H_{_{Y'}})$. {\em Let ${\mathfrak{L}^{^{\pm}}}$ denote the double duals of $\text{Lie}(\mathfrak U_{_{Y'}}^{^\pm})$}.
So ${\mathfrak{L}^{^{\pm}}}$ get canonical Lie sub-algebra structures of $\mathfrak R$  by a Hartogs lemma for Lie sub-bundles of $\mathfrak R$ . It is easily checked that the extensions ${\mathfrak{L}^{^{\pm}}}$ are nilpotent Lie algebra bundles on $\spec B$.

Recall that for a coherent $S$-module $\mathcal F$, we recall \cite[Exp I, 4.6, 4.6.3, 4.6.5]{sga3}, the functor 
\beqa\label{thewconstruct}
{\tt W}:S' \mapsto \Gamma(S', \cF \otimes_{_{\co_{_S}}} \co_{_{S'}}).
\eeqa 
When $\cF$ is locally free, this  is represented by a group scheme ${\tt W}(\cF)$ defined by the affine scheme $\text{Sym}(\cF^{^\vee})$. Since ${\mathfrak{L}^{^{\pm}}}$ are  nilpotent Lie algebra bundles, we can endow 
 the schemes ${\tt W}({\mathfrak{L}^{^{\pm}}})$ over $Y$ with the group law given by Baker-Campbell-Hausdorff formula  which makes ${\tt W}({\mathfrak{L}^{^{\pm}}})$ {\em smooth unipotent group schemes over $Y$ with connected fibres}. 
  
The hypothesis of  characteristic $0$ (or large prime characteristics, see \cite{bdp}) shows that, over $Y'$, we have a canonical identification  $\mathfrak U_{_{Y'}}^{^{\pm}} \simeq {\tt W}(\text{Lie}(\mathfrak U_{_{Y'}}^{^{\pm}}))$ of affine unipotent group schemes. In other words, the unipotent group scheme ${\tt W}({\mathfrak{L}^{^{\pm}}})$ extends the $Y'$-unipotent group scheme $\mathfrak U_{_{Y'}}^{^{\pm}}$.

By \cite[IX, Corollary 1.4]{raynaud}, the $Y'$-group scheme inclusions $\mathfrak U_{_{Y'}}^{^{\pm}} \subset \mathfrak{GL}({\fw})_{_{Y'}}$  extend to $Y$-group scheme inclusions ${\tt W}({\mathfrak{L}^{^{\pm}}}) \hookrightarrow \mathfrak{GL}({\fm})$ on $Y = \spec B$.

Set $$\mathfrak U_{_{Y}}^{^{\pm}} := {\tt W}({\mathfrak{L}^{^{\pm}}})$$
which are {\em smooth unipotent group schemes} over $Y$. This is clearly the product of the schematic closure $\mathfrak U_{_a}$, of the root groups $U_{_a}$, implying {\em a fortiori}, that $\mathfrak U_{_a}$ are {\em flat}.  This proves the main claims of the proposition. 

Note that as Lie-algebra bundles over $Y$, we have the identification $$\text{Lie}(\mathfrak U_{_{Y}}^{^{\pm}} ) = {\mathfrak{L}(\mathfrak U_{_{Y'}}^{^{\pm}})}.$$

The open immersion $j_{_{Y'}}: \mathfrak U_{_{Y'}}^{^+} \times \mathfrak T_{_{Y'}} \times \mathfrak U_{_{Y'}}^{^-}   \hookrightarrow H_{_{Y'}}$ clearly extends as an open immersion $j_{_{Y}}: \mathfrak U_{_{Y}}^{^+} \times \mathfrak T \times \mathfrak U_{_{Y}}^{^-}   \hookrightarrow \mathfrak Y$ over $Y$. This gives the desired open sub-scheme   $\fb \subset \mathfrak Y$.

\epr
Let $$\mathfrak C := \text{Im} (j_{_Y})$$ be the open subscheme of the schematic closure $\mathfrak Y$ of $H_{_L}$ in $\mathfrak{GL}({\fw})$. Then, by \cite[Theorem 2.2.3 (ii)]{bt2}, we derive the important conclusion that the schematic closure  $\mathfrak T\mathfrak U_{_{Y}}^{^+}$ of $T. U_{_{\bf\Phi^{^+}}}$  and $\mathfrak U_{_{Y}}^{^-} \mathfrak T$ of $U_{_{\bf\Phi^{^-}}}. T$  in $\mathfrak{GL}({\fm})$, are flat, closed sub-group schemes of  $\mathfrak{GL}({\fw})$. Since we are in either characteristic zero or large prime characteristics, by Cartier's theorem, one knows that {\em flat group schemes are smooth}. We are now precisely in the setting of \cite[Proposition 2.2.10]{bt2} (see also \cite[Section 6]{bp}) and we get the important consequence that the image  of $\mathfrak C \times \mathfrak C$ in $\mathfrak{GL}({\fw})$ under the product map of $\mathfrak{GL}({\fw})$ is open and dense in $\mathfrak Y$ and is a subgroup scheme of $\mathfrak{GL}({\fw})$, with generic fibre $H_{_L}$, which is flat (and hence smooth). 
The extension of the group scheme is unique because the group scheme is uniquely determined by $H_{_{Y'}}$ by \cite[IX,1.5]{raynaud}.
We summarise this as the following: 

\bth\label{hbt} Under the hypothesis of \eqref{one}, the group scheme $H_{_L}$ extends uniquely as a {\em smooth quasi-affine} subgroup scheme $\bh_{_B} \subset \mathfrak{GL}({\fw})$.
\eeth
We record the basic corollary to this theorem.
\bcor\label{overq} Let $t \in \cq$ and let $B := \hat{\co}_{_{\cq,t}}$ and $L = \text{Fract}(B)$.  Then we have a {\em smooth quasi-affine} group scheme $\bh_{_{B}}$ on $\spec B$ which {\em extends} $H_{_L}$ on $\spec L$. Further, we have an inclusion of smooth group schemes:
\beqa\label{hbt1}
\bh_{_{B}} \stackrel{{\tilde{\Large\text{\cursive s}}}_{_B}}\hra \cu_{_{B}} \stackrel{j_{_B}}\hra \text{GL}(n^2)_{_{B}}.
\eeqa
where $j_{_B}$ (from \eqref{unitsinthebiggroup}) is a closed embedding while ${\tilde{\Large\text{\cursive s}}_{_B}}$ is not.
\ecor
\brem The basic idea of proof of this key theorem has many similarities with the ideas developed in \cite[Section 6]{bp} which are more general. The idea that a desingularisation is at all possible came from the earlier proof. The proof is however shorter here and in turn relies on the basic structure of results in \cite[Section 2.2]{bt2}, in particular on \cite[Proposition 2.2.10]{bt2}.   \erem

\subsection{When $H$ is non-split,  of type $\tt B_{_\ell}$ or $\tt D_{_\ell}$} 

Let $(\bbf,\tq)$ be a quadratic module over $B$, where $\bbf$ is assumed to be a free $B$-module of rank $n = 2\ell, or 2 \ell+1$, such that over the fraction field $L$, the quadratic form $\tq_{_L}$ is non-degenerate and {\em non-hyperbolic}, i.e., we are in the setting of the {\em non-split} case. The result and its proof go through a somewhat different analysis in this situation. 

Consider the positive graded Clifford algebra $\Cl^{+}_\tq(\bbf)$ which is a unital graded associative algebra over $B$. We make no assumptions on the characteristic of the base field except that $2$ is invertible in $k$. With notations as above, in this section we prove the following:
\bth\label{hbt,non-split} The group scheme $H_{_L}$ extends uniquely as a {\em smooth, connected, quasi-affine} subgroup scheme $\bh_{_B}$ such that
\beqa\label{hbt1.5}
\bh_{_{B}} \stackrel{{\tilde{\Large\text{\cursive s}}}_{_B}}\hra \cu_{_{B}} \stackrel{j_{_B}}\hra \text{GL}(n^2)_{_{B}},
\eeqa
where where $j_{_B}$ (from \eqref{unitsinthebiggroup}) is a closed embedding while ${\tilde{\Large\text{\cursive s}}_{_B}}$ is not. The closed fibres of $\bh_{_B}$ can be described as natural subgroups of the generalized Clifford-Lipschitz group associated to the quadratic module over the closed point. The Lie algebra of these groups can be identified with the submodule $\Cl_\tq^{^{+\leq 2}}/B$ \eqref{helmstetter}.
\eeth
The proof is a consequence of the work of Jacques Helmstetter \cite{helm0}, \cite{helm}. For this reason, we will recall facts from these papers which when put together gives the one of the key ingredients for proof of the theorem. An essential point to be kept in mind in this subsection is that, the group scheme $H_{_L}$ is not only non-split but is not {\em quasi-split} as well. This makes these constructions somewhat orthogonal to the ones considered in \cite{bt2}, where quasi-split-ness is an essential assumption.

\subsubsection{On the Clifford-Lipschitz monoid} 
\begin{itemize}

\item Let $\tau$ denote the principal anti-automorphism of the full Clifford algebra $\Cl_\tq$. Recall that $\tau(a) = a$ if $a \in \bbf$ and $\tau(xy) = \tau(y)\tau(x)$ if $x,y \in \Cl_\tq$.
\item We have a $B$-module isomorphism $\Lambda \bbf \simeq \Cl_\tq$ of the exterior algebra with the full Clifford algebra.
\item If $\tq_1$ and $\tq_2$ are two quadratic forms on $\bbf_{_1}$ and $\bbf_{_2}$, then the full Clifford algebra $\Cl(\tq_1 \oplus \tq_2)$ is canonically isomorphic to the {\em graded tensor product} $\Cl(\tq_1) \hat{\otimes} \Cl(\tq_2)$, where the multiplication goes as follows:
\beqa
(x_1 \otimes x_2) . (y_1 \otimes y_2) = \sigma(x_2,y_1).(x_1y_1 \otimes x_2y_2)
\eeqa
where for homogeneous elements $x,y \in \Cl_\tq$, the parity function is given by $\sigma(x,y) = 1$ when $x$ or $y$ is even and $\sigma(x,y) = -1$ if $x$ and $y$ are odd.
\item The module $\Cl(\tq_1) \hat{\otimes} \Cl(\tq_2)$ is equipped with a multiplication which makes it an algebra isomorphic to $\Cl(\tq_1 \oplus -\tq_2)$; the graded tensor product of $\Cl_\tq$ and the opposite graded $\Cl_\tq^\circ$,  which gives the following multiplication:
\beqa
(x \otimes y).(x' \otimes y') =\sigma(y,y').xx' \otimes yy'.
\eeqa
\item In $\bbf \oplus \bbf$, there are two sub-modulies $\Delta$ and $\Delta'$ which are totally isotropic for $\tq \oplus -\tq$, namely:
\beqa
\Delta := \{(a,a) \mid a \in \bbf\} \\
\Delta' := \{(a, -a) \mid a \in \bbf \};
\eeqa
and $\bbf \oplus \bbf$ is the direct sum  $\Delta \oplus \Delta'$ which determines a grading of $\Cl(\tq \oplus -\tq)$.
\item This grading can be transported onto $\Cl_\tq \otimes \Cl_\tq^{^{op}}$, by virtue of the isomorphism which extends the map $(a,b) \mapsto a \otimes 1 + 1 \otimes b$ with $a,b \in \bbf$.
\item Recall that an element $x \in \Cl_\tq$ is called homogeneous if it is odd or even.
\item Since, $\bbf$ is a free $B$-module, as we noted above, we can identify $\Cl_\tq \simeq \Lambda \bbf$. We now work with $(\bbf \oplus \bbf, \tq \oplus -\tq)$ and the totally isotropic sub-modules $\Delta$ and $\Delta'$ for $\tq \oplus -\tq$. Define:
\beqa
\Cl^{^0}(\tq \oplus -\tq,\Delta,\Delta') := \sum _{_j} \Lambda^j(\Delta) . \Lambda ^j(\Delta')
\eeqa
in $\Lambda(\bbf \oplus \bbf)$.

\end{itemize}
\bdefe\label{clifflip} (\cite[Proposition (9), page 22]{helm0}, \cite[5.3.1, 5.3.10]{helm})  Recall that $B$ is a regular local domain.  The Clifford-Lipschitz monoid associated to $(\bbf, \tq)$ is defined as:
\small
\beqa
X(\bbf, \tq):=  \Big\{\mbox{$x \in \Cl_\tq \Big\vert \mbox {x homogeneous and } \big(x \otimes \tau(x)\big)$ lies in} \mbox{ $\text{Im}\big(\Cl^{^0}(\tq \oplus -\tq;\Delta,\Delta')\big)$}\Big\}.
\eeqa
\normalsize
More precisely, for a general base, the Lipschitz monoid is defined as {\em locally homogeneous} elements in $\Cl_\tq$ which satisfies the above conditions. 
Define the Clifford-Lipschitz group  associated to a quadratic module $(\bbf,\tq)$ as follows. For every $B$-algebra $R$,
\beqa\label{lipgroupfunctor}
\text{GLip}(\bbf,\tq) (R) = \Big\{x \in X(\bbf,\tq) \otimes R \big\vert x ~~\text{invertible in}~ \Cl_\tq \Big\}
\eeqa
\edefe
Its properties go as follows:
\begin{enumerate}
\item Both conditions which define $X(\bbf,\tq)$ are closed conditions in $\Cl_\tq$ and therefore give a Zariski closed subset of $\Cl_\tq$. In fact, identifying $\Cl_\tq^{^{op}}$ with $\Cl_\tq$ as modules, we can view the various projections onto the graded summands of $\Cl_\tq \otimes \Cl_\tq^{^{op}}$ as giving a family of quadratic forms on $\Cl_\tq$. The condition that a Lipschitzian element lies in the $0^{th}$-graded piece can be viewed as a vanishing condition on a family of quadratic forms, which defines the Lipschitz monoid.
\item Via the inclusion $X(\bbf, \tq) \subset \Cl_\tq$ we get the inclusion 
\beqa\label{keyinclu}
X(\bbf, \tq) \hra \Cl_\tq \hra \text{End}(\Cl_\tq)
\eeqa
of multiplicative monoids. Define
\beqa\label{posLipmonoid}
X^{+}(\bbf, \tq):= X(\bbf, \tq) \cap \Cl^{+}_\tq 
\eeqa

\item $X(\bbf,\tq)$ contains all elements $a \in \bbf$, and thus, all the elemsnt $\lambda + ab$, where $\lambda \in B$ and $a,b \in \bbf$.
\item If $x \in X(\bbf,\tq)$ then $\tau(x) \in X(\tq)$. 
\item If $x \in X(\bbf,\tq)$ is invertible in $\Cl_\tq$, then $z = x.\tau(x) = \tau(x)x \in B$ is also invertible. Moreover, the inverse $x{^{-1}} = z{^{-1}}.\tau(x)$
and hence the inverse $x{^{-1}}$ also lies in $X(\tq)$. 
\item Via the embedding \eqref{keyinclu}, we conclude that the subgroup $\text{GLip}(\bbf, \tq) \subset X(\bbf, \tq)$ of invertible elements in $X(\bbf, \tq)$ is a  
an open subset of $X(\bbf, \tq)$ in the Zariski topology, and hence an affine group scheme. The functor defined in \eqref{lipgroupfunctor} is simply the functor of points of the group scheme $\text{GLip}(\bbf, \tq)$.

Similarly, we have $\text{GLip}^{+}(\bbf, \tq):= \text{GLip}(\bbf, \tq) \cap X^{+}(\bbf, \tq)$, the subgroup of invertible elements in $X^{+}(\bbf, \tq)$.
\item As seen above, the map $\text{GLip}^{+}(\bbf, \tq) \to B^*$ sending $x \mapsto x.\tau(x)$ is a group homomorphism. Define the subgroup:
\beqa\label{degspin}
{\mathfrak S}(\bbf,\tq) := \text{Ker}\big(\text{GLip}^{+}(\bbf, \tq) \to B^*\big)
\eeqa
\item Being a closed subgroup of $\text{GLip}^{+}(\bbf, \tq)$, ${\mathfrak S}(\bbf,\tq)$ is an affine group scheme over $B$.

\end{enumerate}

We consider the following submodule \eqref{helmstetter}:
\beqa\label{helmstetter1}
\Cl_\tq^{^{+\leq 2}} := \Cl^{+}_\tq \cap \Cl_\tq^{^{\leq 2}}.
\eeqa
Then by \cite[5.5.3, page 263]{helm}, the submodule $\Cl_\tq^{^{+\leq 2}}$ gets a natural Lie algebra structure over $B$ obtained from the associative algebra structure on $\Cl^{+}_\tq$. By taking $k[\epsilon]$-points of $\text{GLip}(\bbf, \tq)$, they can be realized precisely as the Lie algebra of {\em infinitesimal Lipschitzian} elements of the Clifford algebra (see \cite[Section 5.5, page 263]{helm}).

It is seen easily enough that $\Cl_\tq^{^{+\leq 2}}$ is a direct summand of $\Cl^{+}_\tq$, whence it is also projective as an $B$-module. Furthermore, the base ring $B$ is itself the zeroth term in the filtration and is therefore a summand of $\Cl_\tq^{^{+\leq 2}}$. Also, $B$ being commutative, it is closed under the Lie bracket induced from $\Cl^{+}_\tq$. Whence, the quotient $${\tt{spin}_{_B}}:= \Cl_\tq^{^{+\leq 2}} \text{mod} ~B$$ gets a natural Lie algebra bundle structure over $\spec B$.

By the discussion in \eqref{helmstetter}, it can be seen that over $L$, we can identify $$\text{Lie}(\spin(\bbf_{_L}, \tq_{_L}) \simeq {\tt{spin}_{_L}}.$$ Thus, we have shown  the following:
\blem\label{lieextension2} When $H_{_L}$ is non-split of type $B$ or $D$, the Lie algebra $\lie(H_{_L})$ of $H_{_L}$ extends as a Lie algebra bundle $\tt{spin}_{_B}$, to the whole of $\spec B$.\elem

Let $\mathfrak S^{\circ}_{_B} \subset {\mathfrak S}(\bbf,\tq)$ denote the connected component of ${\mathfrak S}(\bbf,\tq)$. As we have noted, over $\spec L$, we have the identification $\mathfrak S^{\circ}_{_L} \simeq \spin(\bbf_{_L}, \tq_{_L})$.

By \cite[Theorem 5.5.6, page 265]{helm}, we deduce the Lie algebra identification:
\beqa\label{lieidenti}
\text{Lie}(\mathfrak S^{\circ}) = {\tt{spin}_{_B}}.
\eeqa

\bth\label{thenonsplitbt} Let $\text{char}(k) = 0$. The connected component ${\mathfrak S}(\bbf, \tq)^\circ$ of ${\mathfrak S}(\bbf, \tq)$ \eqref{degspin} is a smooth group scheme which extends the group scheme $H_{_L}$ over $\spec B$. \eeth
\bpr 

By the very construction of the group scheme ${\mathfrak S}(\bbf, \tq)_{_B}$, it is a {\em closed subgroup scheme} of $\cu_{_{B}} \stackrel{j_{_B}}\hra \text{GL}(n^2)_{_{B}}$ and hence of finite type over $\spec B$. For each $s \in \spec B$, the fibre ${\mathfrak S}(\bbf, \tq)_{_s}$ is a reduced and hence smooth algebraic group over the residue field $k(s)$ (we are over characteristic $0$). Further, by the Lie algebra computations \eqref{lieidenti}, it follows that $\text{dim}({\mathfrak S}(\bbf, \tq)_{_s})$ is constant in $s$ and is in fact $\frac{n(n-1)}{2} -1$. 

Since the base $B$ is reduced, \cite[Theorem 3.10, and Corollary 4.4, Expose VI(B)]{sga3} immediately applies, and it follows that the connected component of the identity ${\mathfrak S}(\bbf, \tq)_{_{B}}^\circ$ is {\em representable and smooth} over $\spec B$. 

%can be identified with the {\em schematic closure} of $H_{_L}$ in $\text{GL}(n^2)_{_{B}}$. Let $L'/L$ be a finite extension which splits $H_{_L}$, and let $B'$ be the normalization of $B$ in $L'$. Hence the morphism $\spec B' \to \spec B$ is finite integral. By base change to $B'$, we get a closed subgroup scheme
%\beqa
%{\mathfrak S}(\bbf, \tq)_{_{B'}} \hra \text{GL}(n^2)_{_{B'}}.
%\eeqa
%Since ${\mathfrak S}(\bbf, \tq)_{_{L'}} \simeq H_{_{L}}$, we may identify the group scheme ${\mathfrak S}(\bbf, \tq)_{_{B'}}$ with the schematic closure of $H_{_{L'}}$ in $\text{GL}(n^2)_{_{B'}}$. Thus we are in the earlier setting of the split case, except that $B'$ is only a semi-local ring. But this is easily handled. In other words, we get a smooth, connected {\em open subgroup scheme} $\bh_{_{B'}} \subset {\mathfrak S}(\bbf, \tq)_{_{B'}}$ which extends $H_{_{L'}}$ to the whole of $\spec B'$. By \cite[Lemma 3.10.1, Expos\'e VI-B]{sga3}, we get the identification:
%\beqa
%\bh_{_{B'}}  \simeq {\mathfrak S}(\bbf, \tq)_{_{B'}}^\circ.
%\eeqa
%Whence, we see that the base change ${\mathfrak S}(\bbf, \tq)_{_{B'}}^\circ$ of ${\mathfrak S}(\bbf, \tq)_{_{B}}^\circ$ is smooth and hence {\em flat}. By \cite{ferrand}, we deduce that the group scheme ${\mathfrak S}(\bbf, \tq)_{_{B}}^\circ$ is itself {\em flat}. Since we have large prime characteristics, it follows that ${\mathfrak S}(\bbf, \tq)_{_{B}}^\circ$ is smooth and connected. This completes the proof of the theorem.

\epr

%\appendix{}
\part{\sc A topological application}
\section{Low degree Betti numbers of $\fm_{_{{\tt G}_{_2}}}$}
In this section, we will briefly touch on the structure of the new moduli spaces $\fm_{_\fg}$. This will enable us to compute some small Betti numbers. {\em We will illustrate this in the case when $\fg$ is the exceptional group ${\tt G}_{_2}$}. 

Recall that by the Borel-de Siebenthal algorithm \cite{bdes}, in the case of $\fg = {\tt G}_{_2}$, we have two maximal rank subgroups which have been seen in \eqref{locnontr}. In fact, in this case the two maximal rank subgroups are $\SL(3)$ and $\text{SO}(4) = \SL(2) \times \SL(2)/\text{center}$. In  \cite{grp}, we find  morphisms 
\beqa
i:M_{_{\SL(3)}} \to M_{_{\fg}}\\ 
j:M_{_{\SL(2) \times \SL(2)}} \to M_{_{\text{SO}(4)}} \to M_{_{\fg}}
\eeqa
The moduli space $M_{_{\text{SO}(4)}}$ has two connected components which we denote by $M^{+} \cup M^{-}$. The complement of the smooth locus (which is the locus of regularly stable bundles), is the union of two connected components, namely $i(M_{_{\SL(3)}}) \cup j(M^{+})$, and $j(M^{-})$, where, by our convention, the trivial ${\tt G}_{_2}$-torsor lies in $j(M^{+})$. 

Our limited purpose is to determine small Betti numbers and so standard codimension computations allow us to assert that we can work with an open subset $M^{^o} \subset M_{_{\fg}}$ whose complement has codimension strictly bigger than $1$. This open subscheme consists of the regularly stable bundles $M_{_{\fg}}^{rs} \subset M_{_{\fg}}^{s}$,  and the complement $M^{^o} \setminus M_{_{\fg}}^{rs}$ has exactly  three connected components which we will denote by $U(\SL(3))$, $U^+$ and $U^-$.  

Since $M_{_{\fg}}^{rs}$  is known to be {\em unirational}, so is $\fm_{_\fg}$ and moreover since $\fm_{_\fg}$ is smooth and projective it is simply connected. Hence, it is immediate that $H^{^0}(\fm_{_\fg}, \mathbb Q) = \mathbb Z$, $H^{^1}(\fm_{_\fg}, \mathbb Q) = (0)$. 

Let us denote the three components of the exceptional divisor coming from the subsets $U(\SL(3))$, $U^+$ and $U^-$ by  $E(\SL(3))$, $E^+$ and $E^-$ respectively. We will briefly describe these components when we compute the third Betti number.

It is known (see \cite{grp} for other references) that the Picard group of $M_{_{\fg}}$, and  hence of the smooth locus $M_{_{\fg}}^{rs}$, is $\mathbb Z$. Let ${\tt Z}_{_{\fg}} := M_{_{\fg}}^{rs} \cup E$, where $E = E(\SL(3)) \bigsqcup E^+ \bigsqcup E^-$. Whence, it follows that $H^{^2}({\tt Z}_{_{\fg}}, \mathbb Z)  = H^{^2}(\fm_{_\fg}, \mathbb Z) = \mathbb Z^4$, where the extra three copies of $\mathbb Z$ are contributed by the three connected components of the exceptional divisor.

\subsubsection{The third Betti number} We will now determine the third Betti number $h^3(\fm_{_\fg})$ and will also indicate why  $H^{^3}(\fm_{_\fg}, \mathbb Z)$ is {\em torsion free}. 

To compute the third cohomology, we need to work a bit harder. Let $\mathcal M(\fg)$ denote the moduli stack of principal $\fg$-bundles on $C$. By \cite[Theorem II.6]{faltings}, it follows that the {\em regularly stable locus} $\mathcal M(\fg)^{rs}$ is an open substack of $\mathcal M(\fg)$, and the complement is of codimension strictly bigger than $1$. By \cite{bishol}, it follows that $\mathcal M(\fg)^{rs}$ is a gerbe over the moduli space $M_{_{\fg}}^{rs}$ banded by the center $Z(\fg)$. Since we deal with $\fg = {\tt G}_{_2}$, the center is trivial. It is now easy to deduce that in order to determine the third Betti number of $M_{_{\fg}}^{rs}$, we need to determine it for the moduli stack  $\mathcal M(\fg)$. This will be the first step to determine $h^3(\fm_{_\fg})$.

To do this, we naturally appeal to \cite{ab} and we follow its notations. We fix a holomorphic principal $\fg$-bundle $E$ on $C$ with fixed Chern classes. Let $\mathcal C$ be the space of holomorpihc structures  on $E$. This is a complex affine space modeled on the space $\Omega^{0,1}(E(\mathfrak g))$, where $E(\mathfrak g)$ is viewed as a  $C^{\infty}$-bundle. The concept of a Harder-Narasimhan type associated to a canonical parabolic reduction is well-defined and is indexed by the characters of the Levi subgroups of the parabolic subgroup. These are denoted by $\mu$ and the Shatz stratification gives a stratification $\mathcal C_{_\mu}$ of $\mathcal C$. The codimension of the strata are denoted by $d_{_\mu}$ and the stratum associated to $d_{_\mu} = 0$ is the {\em semi-stable} stratum. It is shown in \cite[Theorem 10.10]{ab} that this stratification is equivariantly perfect and as a result, we get the following relation:
\beqa\label{atiyahbott}
P_{_t}(\mathcal C) = \sum_{_{\mu}} t^{^{2.d_{_\mu}}} P_{_t}(\mathcal C_{_\mu})
\eeqa
The equivariant Poincar\'e polynomial gives the Poincar\'e polynomial of the moduli stack $\mathcal M(\fg)$. The codimension statement allows us to conclude that, to compute the third Betti number of $M_{_{\fg}}^{rs}$, we need simply compute the coefficient of $t^{^3}$ on the left side of \eqref{atiyahbott}. Since $\mathcal C$ is contractible, the computation of the equivariant Poincar\'e polynomial reduces to the computing the Poincar\'e polynomial of the stack of maps $\text{Maps}(C, B\fg)$. We now use the re-interpretation of the Atiyah-Bott formula in \cite[Theorem 1.5.2.3]{gail}. An easy computation gives us:
\bprop The third rational cohomology of the moduli space of regularly stable $\fg $-bundles with $\fg ={\tt G}_{_2}$ is $H^{^3}(M_{_{\fg}}^{rs}, \mathbb Q) \simeq \mathbb Q^{^{2g}}$, where $g = \text{genus}(C)$.  \eprop
What remains is the computation of the third Betti number of the smooth compactification $\fm_{_\fg}$. By codimension computations, we will work with the open subset ${\tt Z}_{_{\fg}} \subset \fm_{_\fg}$, where we delete the subsets of codimension bigger than $1$. 

\subsubsection{The exceptional divisor} As noted earlier, we have ${\tt Z}_{_{\fg}} := M_{_{\fg}}^{rs} \cup E$, where $E = E(\SL(3)) \bigsqcup E^+ \bigsqcup E^-$. The subscheme $E$ is an open subset of the full exceptional divisor and as we have observed, for the computation of low degree Betti number, we can afford to discard high codimension subsets.

To describe the $E(\SL(3))$, we need to consider the locus of regularly stable $\SL(3)$-bundles on $C$. By extension of structure group via $i:\SL(3) \hra {\tt G}_{_2}$, this gives stable ${\tt G}_{_2}$-torsors.  Let $\ce \in M_{_{\fg}}^{s}$ be a point which is stable as a  ${\tt G}_{_2}$-torsor but not regularly stable as a ${\tt G}_{_2}$-torsor. Following the argument of \eqref{secondpart}, the fibre over $\ce$ are torsors under the closed fibre of a maximal parahoric group scheme whose closed fibre is $\bh'_{_o}$. This non-reductive group is such that its Levi quotient is semisimple and in this case isomorphic to $\SL(3)$. It is known, (see \cite[2.21]{pr}) that the unipotent radical $U_{_o}$ of  $\bh'_{_o}$ has a filtration whose composition factors are {\em vector groups}. The fibre is now easily seen to be a projective space of the vector space $H^{^1}(U_{_o})$ whose dimension is $\text{dim}(U_{_o})(g-1)$ (\cite[Cor 8.1.5]{behrend}). Observe that all the bundles are of degree $0$. Hence, the open locus $E(\SL(3))$ is a projective bundle with fibre as above. 

Similar descriptions can be given for the other components $E^+$ and $E^-$.

We now use the Thom-Gysin sequence for integral cohomology and obtain the following exact sequence:
\beqa
\ldots \to H^{^1}(E, \mathbb Z) \to H^{^3}({\tt Z}_{_{\fg}}, \mathbb Z) \to H^{^3}(M_{_{\fg}}^{rs}, \mathbb Z) \to H^{^2}(E, \mathbb Z) \to H^{^4}({\tt Z}_{_{\fg}}, \mathbb Z) \ldots
\eeqa
where the maps $H^{^j}(E, \mathbb Z) \to H^{^{j+2}}({\tt Z}_{_{\fg}}, \mathbb Z)$ is the Gysin map induced by cupping with the Euler class of the submanifold $E \subset {\tt Z}_{_{\fg}}$. Computations following the ones in \cite{balajithesis}, shows that $H^{^1}(E, \mathbb Z) = 0$ and the map $ H^{^2}(E, \mathbb Z) \to H^{^4}({\tt Z}_{_{\fg}}, \mathbb Z)$ is {\em injective}. Whence we conclude that 
\beqa
H^{^3}({\tt Z}_{_{\fg}}, \mathbb Z) \to H^{^3}(M_{_{\fg}}^{rs}, \mathbb Z).
\eeqa
is an isomorphism. It is shown very generally in \cite{bishol}  that the cohomological Brauer group of $M_{_{\fg}}^{rs}$ for any simply connected $\fg$ is trivial. Hence we conclude by the computations as above that the cohomological Brauer group of ${\tt Z}_{_{\fg}}$ and hence of $\fm_{_\fg}$ is trivial. These are unirational, smooth proper schemes and hence the third cohomology is torsion free and we get $H^{^3}(\fm_{_\fg}, \mathbb Z) \simeq \mathbb Z^{^{2g}}$.

\part{\sc Appendix(with Rajarshi Kanta Ghosh)}
\section{Local Models following Nori}
\subsubsection{Recollection of Nori's construction from \cite{nori}}
Let $k$ be an algebraically closed field and $A = k\{X_1, \ldots, X_n\}$ be the non-commuting polynomial ring over $k$. In \cite{artin}, M.Artin   constructs a moduli for $A$-modules $M$ which are of dimension $n$ over $k$.

Once a basis is chosen for $M$, the action of $X_i$ defines an element of $\cm_{_n}(k)$, which we denote again by $X_i$. A change of basis for $M$ replaces $X_i$ by $h.X_i. h^{-1}$ where $h$ runs through all invertible matrices. We consider the action of $\PGL(n)(k)$ on the $g$-fold product of $\cm_{_n}(k)$ given by $h \cdot (X_1, \ldots, X_g)=(h.X_1.h^{-1}, \ldots , h.X_g. h^{-1})$ and let  
\beqa\label{artinquotient}
Z_{_n}:= \cm_{_n}(k)^g\parallelslant \PGL(n)
\eeqa 
denote the GIT quotient. If $\cm_{_n}^{g,s}$ denotes the sub-scheme of stable points then we have a principal $\PGL(n)$-bundle
\beqa\label{basicpglbundle} 
 \cm_{_n}^{g,s}\to Z_{_n}^{^s}
\eeqa 

 %Let $H_1, H_2, \ldots , H_l$ be arbitrary monomials in the $X_i$, and consider the expression $P=\mbox{\rm det }(t_0+t_1H_1+\cdots +t_lH_l)= \sum_I t^If_I$, where $t^I$ denotes monomials in $t_0, \ldots ,t_l$ and $f_I4$ are functions on $M(n)^g$, clearly invariant under the group action.

%Artin proves the following statements (modulo Mumford's conjecture).

The points of the GIT quotient $Z_{_n}$ can be described as $S$-equivalence classes, or equivalently as isomorphism of associated graded of Jordan-Holder filtrations. More precisely, for any $A$-module $M$, we let $M=M_0\supset M_1\supset \ldots \supset M_r=0$ be a Jordan-Holder filtration, and define the associated graded $\text{gr} (M):=\oplus_iM_i/M_{i+1}$. Let $X=(X_i)$ and $Y=(Y_i)$ be points of $\cm_{_n}(k)^g$ which represent modules $M$ and $M'$ respectively. Then the orbit closures of $X$ and $Y$ intersect if and only if gr $(M)$ and gr $(M')$ are isomorphic. The corresponding equivalence classes give the points of $Z_{_n}$.

%(2) Consider the scheme $M(n)^g$ with the action of $Gl(n)$ defined over $\bz$, and let $B_n$ be the ring of invariant functions with $B'_n$ the subring generated by the $f_I$ for all possible choices of the monomials $H_j$. Then $B_n$ and $B'_n$ are both finitely generated $\bz$-algebras and $B_n$ is integral over $B_n'$.

%(3) Spec $B_n=Z_n$ is the coarse moduli space for the corresponding functor.

In \cite{nori}, Nori constructs a natural desingularisation of $Z_{_2}$ based on the construction of Seshadri \cite{desing}; in fact, as in \cite{desing}, he constructs birational models for $Z_{_n}$ for all $n$ but these  are singular (for the very same reason), for all $n \geq 3$.   

Let $A = \bz\{X_1, X_2, \ldots , X_g\}$ be the non-commuting polynomial ring over the integers. For any commutative ring $R$, let $\ch\text{ilb}_n (R)$ denote the isomorphism classes of pairs $(M, m)$ where $m \in M$ and $M$ is a $R\otimes_\bz A$-module, locally free of rank $n$ as an $R$-module, such that $m$ generates $M$. Equivalently, $\ch\text{ilb}_n(R)$ is the set of left ideals $I$ in $R\otimes_\bz A$ such that $M= R\otimes_\bz A/I$ is locally free of rank $n$ as an $R$-module. In fact, there is a  quasi-projective scheme $\Hilb_n$ which represents the functor $\ch\text{ilb}_n$, i.e.,   $\text{Mor}(\spec R, \ch\text{ilb}_n) = \Hilb_n(R)$.

Let $W$ denote the vector space $k^{^n}$, so that we can identify $\cm_{_n}(k)$ with $\End(W)$. Let $w \in W$ be a fixed vector and let $Y_{_n} := \Hom(\langle w \rangle, W) \times \End(W)^g$. We have an action of $\Aut(W) = \GL(n)$ on $Y$ given by $h. (\phi, M_{_1}, \ldots, M_{_g}) := (h \circ \phi, h.M_{_1}.h^{-1}, \ldots,h.M_{_g}.h^{-1})$. Consider the subset of points $Y_{_n}^{^s} \subset Y_{_n}$, consisting of $(\phi, M_{_1}, \ldots, M_{_g})$ such that $k\{X_1, X_2, \ldots , X_g\}.\phi(w) = W$. Then the action of $\GL(n)$ on $Y_{_n}^{^s}$ is {\em free} and Nori shows that there is a principal $\GL(n)$-bundle
\beqa\label{Hn}
h_{_n}:Y_{_n}^{^s} \to \Hilb_{_n} (= Y_{_n}^{^s} \big/\GL(n))
\eeqa
and hence $\Hilb_{_n}$ is smooth and irreducible of dimension $n + (g-1)n^2$.

We now work with $\Hilb_{_{n^2}}$ and let $\Hilb^{{\boldsymbol\tau}}_{_{n^2}}$ be the closed sub-scheme of $\Hilb_{_{n^2}}$ that corresponds to {\em two-sided} ideals. Denote, by $P_{_{n^2}}$, the restriction of the principal $\GL(n^2)$-bundle $Y_{_{n^2}}^{^s}$ on $\Hilb_{_{n^2}}$ to $\Hilb^{{\boldsymbol\tau}}_{_{n^2}}$. Then we can identify, $\text{Mor}(\spec R, P_{_{n^2}})$ with the set of all two-sided ideals $I$ of $R \otimes A$ with a specified basis $e_{_1}, e_{_2}, \ldots, e_{_{n^2}}$ of $R\otimes A/I$.

%Let $Q_n$ be the scheme of all algebra structures on a free module with basis $e_1, e_2, \ldots, e_n$; to be precise, Mor (Spec $R, Q_n$)$=$isomorphism classes of paris $(B, e)$ where $B$ is an $R$-algebra and $e=(e_1, e_2, \ldots , e_n)$ is a (free) basis for $B$ as an $R$-module. Then $Q_{n^2}$ naturally inherits a sheaf of algebras $B$ and the points of $Q_{n^2}$ where $B$ is a matrix algebra is an open set $W_n$ of $Q_{n^2}$, (this is so because $W_n$ is precisely the set of points of $Q_{n^2}$ where the module homomorphism $B\otimes B^0\to \End B=B\otimes B^*$ has maximal rank).

%There is a natural map $g_n \colon P_n \to Q_n$,. Let $H_n=g_{n^2}^{-1}(W_n)$. Then $H_n$ is an open set of $P_{n^2}$ which is clearly $Gl(n^2)$-invariant, and so $H_n$ is the restriction of $P_{n^2}$ to an open subset $V_n$ of $\Hilb^{{\boldsymbol\tau}}_{_{n^2}}$. 

Define the {\em open} subset $V_{_n} \subset \Hilb^{{\boldsymbol\tau}}_{_{n^2}}$ as follows. A geometric point $\spec k\to V_{_n}$ corresponds to a two-sided ideal $I$ of $k\otimes A$ such that $k\otimes A/I=\cm_{_n}(k)$. Let $Z_{_n}^{^s}$ denotes the open subset of $Z_{_n}$ which corresponds to simple modules. Then one sees that we can make the identification $Z_{_n}^{^s} \simeq V_{_n}$. Let $\overline{V}_{_n}$ denote the closure of $V_{_n}$ in the closed subscheme $\Hilb^{{\boldsymbol\tau}}_{_{n^2}}$ of $\Hilb_{_{n^2}}$, with its reduced induced structure.

The projection $\pi_{_{n^2}}:Y_{_{n^2}} \to \cm_{_{n^2}}^g$ induces a morphism $\pi_{_{n^2}}:\Hilb_{_{n^2}} \to Z_{_{n^2}}$. We also have a closed immersion $\Delta:Z_{_n} \hra Z_{_{n^2}}$, given by $\text{gr}(V) \mapsto \text{gr}(V) \oplus \ldots \oplus \text{gr}(V)$. The restriction of the morphism $\pi_{_{n^2}}$ to $\overline{V}_{_n} \subset \Hilb_{_{n^2}}$  gives a projective morphism $\pi_{_{n^2}}:\overline{V}_{_n} \to Z_{_{n^2}}$ which factors as:
\beqa\label{noridesings}
\begin{tikzcd}
	{\overline{V}_{_n}} \\
	{Z_{_n}} & {Z_{_{n^2}}}
	\arrow[from=1-1, to=2-1]
	\arrow["{\Delta_{_n}}", hook, from=2-1, to=2-2]
	\arrow[from=1-1, to=2-2]
\end{tikzcd}
\eeqa 
In {\em loc cit},  Nori proves the following.
\bth\label{nori} The morphism  $\pi_{_{n^2}}: \overline{V}_{_n} \to Z_{_n}$ is {\em proper and birational}. Furthermore,  for $n = 2$, $\pi_{_4}:\overline{V}_{_2} \to Z_{_2}$ gives a {\em desingularisation}.  \eeth
\brem As we have remarked earlier, since $\ca_{_n}$ is invariably non-smooth for $n > 2$, the scheme $\overline{V}_{_n}$ is also invariably {\em non-smooth}.\erem
\brem In the proof of \eqref{nori}, it is shown that given a point $\xi' \in \overline{V}_{_n}$, there is an analytic  neighbourhood $U(\xi')$ of $\xi'$ and a smooth morphism $\vartheta:U(\xi') \to \ca_{_n}$. Let $U(\xi')_{_{\cq}}$ be defined as the base-change:
\beqa\label{localversalhere}
\begin{tikzcd}
	{U(\xi')_{_{\sf Q}}} & {U(\xi')} \\
	{\sf Q} & {\mathcal A_{_n}}
	\arrow["{^\vartheta}", from=1-2, to=2-2]
	\arrow[from=1-1, to=1-2]
	\arrow["{^\vartheta}"', from=1-1, to=2-1]
	\arrow["{^\theta}"', from=2-1, to=2-2]
\end{tikzcd}
\eeqa
where $\theta$ is as in  \eqref{forversality}.

\erem

\subsubsection{Preliminaries}
Let $\tc{s}: H \hra \GL_n$ denote the irreducible and faithful representation \eqref{thespinorrep}. This gives us a map $\lie (H)^g \parallelslant \overline{H}\to \cm_{_n}^g \parallelslant \PGL_n$, where $\overline{H} = H /\text{center}$, and the action of groups on their respective Lie algebras is via the adjoint action. Let $\lie (H)^{g,s}$ denote the   the open subset {\em where the action of  $\overline{H}$ on $\lie (H)^g$ is free}.

Let $Z_{_n} :=  \cm_{_n}^g \parallelslant \PGL(n)$ and $Z(H):= \lie (H)^g \parallelslant \overline{H}$. Then we have the canonical map:
\beqa\label{Hquottoartinquot}
\tc{s}_{_*}: Z(H) \to Z_{_n}.
\eeqa
% https:\parallelslantq.uiver.app/?q=WzAsMixbMCwwLCJaKEgpIl0sWzAsMSwiWihIKS8vXFxvdmVybGluZXtIfSJdLFswLDEsInAiXV0=

Observe that $\lie (H)^g$ represents the functor $F:\text{k-algebras} \to \on{Sets} $ given by 
\begin{align*}
    R&\mapsto \{\text{Lie algebras homomorphisms }\Lambda_g(R)\to \lie H\otimes R\},
\end{align*}
where $\Lambda_g(R)$ denotes the free Lie algebra in $g$ variables. Maps from $\Lambda_g(R)$ to $\lie H\otimes R$ can also be viewed as maps from $R\{X_1,\ldots, X_g\}$ to $\cm_{_n}\otimes R$, such that the images of $X_1,\ldots, X_g$ lie in the image of the spinor representation of the Lie algebra $\lie H$ in $\cm_{_n}$. We will usually adopt this perspective and denote the data of maps from $\Lambda_g(R)$ to $\lie H\otimes R$  as shown below\footnote{Here and elsewhere, the dotted arrows are meant to indicate where the images of the $X_j$'s land.}.
% https:\parallelslantq.uiver.app/?q=WzAsMyxbMCwwLCJSXFx7WF8xLFxcbGRvdHMsWF9nXFx9Il0sWzIsMCwiTV9uXFxvdGltZXMgQSJdLFsxLDEsIlxcbGllIEhcXG90aW1lcyBBIl0sWzAsMV0sWzAsMiwiIiwyLHsic3R5bGUiOnsiYm9keSI6eyJuYW1lIjoiZG90dGVkIn19fV0sWzIsMV1d
\beqa\label{pointsofzh}
\begin{tikzcd}
	{R\{X_1,\ldots,X_g\}} && {\cm_{_n}\otimes R} \\
	& {\lie H\otimes R}
	\arrow[from=1-1, to=1-3]
	\arrow[dotted, from=1-1, to=2-2]
	\arrow[from=2-2, to=1-3]
\end{tikzcd}
\eeqa
{\em The goal in this appendix is to desingularise the quotient scheme $Z(H)$.}

Observe that the action of $\overline{H}$ on $\lie (H)^{g,s}$ being free, we have a principal $\overline{H}$-bundle 
\beqa\label{honestprbundle}
{\lie (H)^{g,s}} \stackrel{\overline{H}}\to Z(H)^s
\eeqa
which in turn implies the smoothness of the subscheme  $Z(H)^s$ in $Z(H)$. 

% https:\parallelslantq.uiver.app/?q=WzAsMixbMCwwLCJaKEgpXnMiXSxbMCwxLCJaKEgpXnMvXFxvdmVybGluZXtIfSJdLFswLDEsIlxcb3ZlcmxpbmV7SH0iXV0=
%\beqa\label{honestprbundle}
%\begin{tikzcd}
	%{\lie (H)^{g,s}} \\
	%{Z(H)^s}
	%\arrow["{\overline{H}}", from=1-1, to=2-1]
%\end{tikzcd}
%\eeqa

Let us denote by $Z^s_{\tsf{f}}(H)$ the pullback of $Z_{_n}^s \subset Z_{_n}$ under $\tc{s}_{_*}$ \eqref{Hquottoartinquot}. This subset is analogous to the subset $M_{_{\tt f}}(H)$ of full holonomy bundles in \S\ref{towardsmain}. It is seen easily that 
$Z^s_{\tsf{f}}(H)$ is an open subset of $Z(H)^s$.

 For $g \geq 2$, observe that the subset $Z^s_{\tsf{f}}(H)$ is {\em non-empty}. This will use the same remark made in \eqref{towardsmain} for showing that the space of principal bundles with full holonomy is non-empty. 
 
 %Indeed, given a surjection $\Lambda_g\to \lie H$, we observe that the composite $\Lambda_g\to \lie H\to \mathfrak{gl}_n$ is an irreducible representation. This implies that the associated map $k\{X_1,\ldots, X_g\}\to \cm_{_n}(k)$ is an irreducible representation, i.e. a surjection, and this gives a point of $Z_{_n}^s$, proving the claim. 
\brem\label{basicoverfull} We remark that this principal $\overline{H}$-bundle \eqref{honestprbundle} is in fact the reduction of structure group of the principal $\PGL(n)$-bundle \eqref{basicpglbundle} on $Z_{_n}^{^s}$ pulled back to $Z^s_{\tsf{f}}(H)$.\erem
Let 
\beqa
\Gamma' := \overline{V}_{_n} \times_{_{Z_{_n}}} Z(H)
\eeqa
with the reduced induced structure. Observe that $\Gamma'$ is birational to $Z(H)$. Let $\cz_{_H}$ denote the  component of the normalization of $\Gamma'$ {\em which contains the locus $Z^s_{\tsf{f}}(H)$}. We therefore have a diagram
\beqa\label{befthm}
\begin{tikzcd}
	{\mathcal Z_{_H}} \\
	& {\Gamma'} & {\overline{V}_{_n}} \\
	& {Z(H)} & {Z_{_n}}
	\arrow[from=2-2, to=2-3]
	\arrow[from=2-2, to=3-2]
	\arrow[from=3-2, to=3-3]
	\arrow[from=2-3, to=3-3]
	\arrow[from=1-1, to=2-2]
	\arrow["{_\pi}", from=1-1, to=3-2]
	\arrow[from=1-1, to=2-3]
\end{tikzcd}
\eeqa
The projection $\pi$ is a proper, birational morphism, being obtained by a base change of one such morphism.
\bth\label{mainthappendix}
Let $H$ be as in Theorem \ref{mainintro}. The canonical projection \eqref{befthm}
\beqa\label{thedsingmaphere}
\pi:\cz_{_H} \to Z(H)
\eeqa
 gives a desingularisation of $Z(H)$. \eeth
The rest of the article is devoted to proving this theorem.

\subsubsection{On the structure of $\cz_{_H}$} 
Let $\xi\in \cz_{_H}$ be a closed point and  $\xi' \in \overline{V}_{_n}$ be its image. Let $R_{\xi}:=\widehat{\co_{\xi}}$, and $Z_{\xi}:=\spec R_{\xi}$. Let $U(\xi')_{_{\sf Q}}$ be as in \eqref{localversalhere}. Thus we have a diagram:
\beqa
\begin{tikzcd}
	{Z_{_\xi}} & {U(\xi')_{_{\sf Q}}} \\
	& {\sf Q}
	\arrow[from=1-1, to=1-2]
	\arrow["{^\vartheta}", from=1-2, to=2-2]
	\arrow["{^\vartheta{_H}}"', from=1-1, to=2-2]
\end{tikzcd}
\eeqa

Let $\vartheta_{_H}(\xi) = t$. Let $B := \hat{\co}_{_{\cq,t}}$ denote the completion of the local ring of $\tsf{Q}$ at $t$ as in \eqref{overq}. Whence \eqref{hbt1}, one has the group schemes $\bh_{_B}, \cu_{_B}$ and $\GL({n^2})_{_B}$ on $\spec B$ with the inclusions. These  group schemes pullback to $Z_{_\xi} $ by the morphism $\vartheta_{_H}$. 

Let $K$ denote the fraction field of $R_{\xi}$. By virtue of the map $\spec K\to V_{_n}$, we obtain an Azumaya algebra $ \ca_{_K}$ on $\spec K$, which is obtained as the associated fibration of the $\PGL(n)$-torsor   $\cm_{_n}^{g,s} \to Z_{_n}^{^s} \simeq V_{_n}$
pulled back from $V_{_n}$ (see \eqref{basicpglbundle}). The Azumaya algebra $\ca_{_K}$ extends as a \textit{degenerate} Azumaya algebra bundle on $Z_{_\xi}$ in the following manner. The map $Z_{_\xi} \to U(\xi')_{_{\sf Q}} \subset  \overline{V}_{_n}$ gives a $R_{\xi}$-valued point of $\Hilb_{n^2}$. By the versality of $\tsf{Q}$ this degenerate Azumaya bundle (with its induced Lie algebra structure) coincides with the Lie algebra bundle $\lie \cu_q$ pulled back by $\vartheta_{_H}:Z_{_\xi} \to \spec B$. We denote this algebra bundle on $Z_{_\xi}$ by $\ca_{_{R_{\xi}}}$. Since it arises from the Hilbert scheme, we in fact have that $\ca_{_{R_{\xi}}}$ arises as a quotient $R_{\xi}\{X_1,\ldots, X_g\}\twoheadrightarrow \ca_{_{R_{\xi}}}$. 

Our {\em claim is that this arrow factors} through the Lie algebra inclusion of $\lie~\bh_{_{Z_{_\xi}}}$ in $\lie \cu_{_{Z_{_\xi}}}$ (induced by \eqref{hbt2}). 

First we show that over the generic point the map $K\{X_1,\ldots, X_g\}\twoheadrightarrow \ca_{_K}$ factors through $\lie{H}_{_K}$. Observe that the $\PGL(n)$ torsor pulled back to $\spec K$ has a reduction of structure group to $\overline{H}$ in the following manner. The generic point $\spec K$ maps to $Z^s_{\tsf{f}}(H)$, and one pulls back the $\overline{H}$-torsor to $\spec K$. Let $E$ be this pulled back $\overline{H}$-torsor. Via the inclusion $\lie H\hookrightarrow \mathfrak{gl}_{n}$, we see that $E(\PGL(n))$ is the aforementioned $\PGL(n)$ torsor pulled back from $V_{_n}$ (see Remark \ref{basicoverfull}). 

Upon going to a finite extension of $K_{\xi}$ we can ensure a splitting of $E$. Let $K'$ denote this extension. On $\spec K'$ we thus have the data
% https:\parallelslantq.uiver.app/?q=WzAsMyxbMCwwLCJLJ1xce1hfMSxcXGxkb3RzLCBYX2dcXH0iXSxbMiwwLCJNX25cXG90aW1lcyBLJyJdLFsxLDEsIlxcbGllIEhfe0snfSJdLFswLDFdLFswLDIsIiIsMix7InN0eWxlIjp7ImJvZHkiOnsibmFtZSI6ImRvdHRlZCJ9fX1dLFsyLDFdXQ==
\[\begin{tikzcd}
	{K'\{X_1,\ldots, X_g\}} && {\cm_{_n}\otimes K'} \\
	& {\lie H_{_{K'}} \cong \lie H \otimes K' \cong \lie H_{_K} \otimes K'}
	\arrow[from=1-1, to=1-3]
	\arrow[dotted, from=1-1, to=2-2]
	\arrow[from=2-2, to=1-3]
\end{tikzcd}\]
Since the sub-bundle $\lie H \otimes K'$ is the pullback of $\lie \bh_{_K}$, it is easy to see that the sections of $\ca_{_K}$, namely $X_1,\ldots, X_g$, must have been in $\lie \bh_{_K}$ to begin with.% This follows from the fact that $\lie H_K = \lie H_{K}\otimes $

Next, recall that the sections $X_1,\ldots, X_g$ extend as those of the bundle $\ca_{_{R_{\xi}}}$. Over any height $1$ prime $\mathfrak p$ of $R_{\xi}$, with local ring $S$, we have that $X_i\in \ca_{_S}\cap \lie H_{_K} = \lie \bh_{_S}$. Now, by Hartogs' lemma, since the base is normal, we have that the sections $X_1,\ldots, X_g$ extend as sections of $\lie~ \bh_{_{R_{\xi}}}$ itself. In conclusion, we have the following data on $Z_{_\xi} = \spec R_{\xi}$. 
% https:\parallelslantq.uiver.app/?q=WzAsMyxbMCwwLCJSX3tcXHhpfVxce1hfMSxcXGxkb3RzLFhfZ1xcfSJdLFsyLDAsIlxcc2Nye0F9X3tSX3tcXHhpfX0iXSxbMSwxLCJcXGxpZSBcXGZyYWt7SH1fe1Jfe1xceGl9fSJdLFswLDFdLFswLDIsIiIsMix7InN0eWxlIjp7ImJvZHkiOnsibmFtZSI6ImRvdHRlZCJ9fX1dLFsyLDFdXQ==
\beqa\label{dataonrxi}
\begin{tikzcd}
	{R_{\xi}\{X_1,\ldots,X_g\}} && {\ca_{_{R_{\xi}}}} \\
	& {\lie~ \bh_{_{R_{\xi}}}}
	\arrow[from=1-1, to=1-3]
	\arrow[dotted, from=1-1, to=2-2]
	\arrow[from=2-2, to=1-3]
\end{tikzcd}
\eeqa

\subsubsection{The smoothness of $\cz_{_H}$}
Let $\spec D_0$ be an Artin local ring, $\spec D$ a complete regular local ring of which $D_0$ arises as a quotient. Consider a map $\phi_0:\spec D_0\to \cz_{_H}$ such that the closed point of $\spec D_0$ maps to $\xi$ as above. In order to prove smoothness, one wants to lift the map $\phi_0$ to $\phi:\spec D\to Z_{_\xi}$ such that the following diagram commutes.
% https:\parallelslantq.uiver.app/?q=WzAsMyxbMCwwLCJcXHNwZWMgRCJdLFswLDEsIlxcc3BlYyBEXzAiXSxbMSwwLCJcXEdhbW1hIl0sWzEsMCwiIiwwLHsic3R5bGUiOnsidGFpbCI6eyJuYW1lIjoiaG9vayIsInNpZGUiOiJ0b3AifX19XSxbMSwyLCJcXHBoaV8wIiwyXSxbMCwyLCJcXHBoaSIsMCx7InN0eWxlIjp7ImJvZHkiOnsibmFtZSI6ImRhc2hlZCJ9fX1dXQ==
\[\begin{tikzcd}
	{\spec D} & Z_{_\xi} \\
	{\spec D_0}
	\arrow[hook, from=2-1, to=1-1]
	\arrow["{\phi_0}"', from=2-1, to=1-2]
	\arrow["\phi", dashed, from=1-1, to=1-2]
\end{tikzcd}\]
Consider the further composite to $\spec B \hra \tsf{Q}$ via $\vartheta_{_H}:Z_{_\xi} \to \spec B$. By the smoothness of $\tsf{Q}$, we obtain a lift of $\vartheta_{_H} \circ \phi_0:\spec D_0\to \tsf{Q}$ to $\spec D$ as follows. 
% https:\parallelslantq.uiver.app/?q=WzAsNCxbMCwwLCJcXHNwZWMgRCJdLFswLDEsIlxcc3BlYyBEXzAiXSxbMSwwLCJcXEdhbW1hIl0sWzEsMSwiXFx0c2Z7UX0iXSxbMSwwLCIiLDAseyJzdHlsZSI6eyJ0YWlsIjp7Im5hbWUiOiJob29rIiwic2lkZSI6InRvcCJ9fX1dLFsxLDIsIlxccGhpXzAiLDFdLFswLDIsIlxccGhpIiwwLHsic3R5bGUiOnsiYm9keSI6eyJuYW1lIjoiZGFzaGVkIn19fV0sWzEsMywicFxcY2lyYyBcXHBoaV8wIiwyXSxbMiwzLCJwIl0sWzAsM11d
\[\begin{tikzcd}
	{\spec D} & Z_{_\xi} \\
	{\spec D_0} & {\tsf{Q}}
	\arrow[hook, from=2-1, to=1-1]
	\arrow["{\phi_0}"{description}, from=2-1, to=1-2]
	\arrow["\phi", dashed, from=1-1, to=1-2]
	\arrow["{\vartheta_{_H}\circ \phi_0}"', from=2-1, to=2-2]
	\arrow["\vartheta_{_H}", from=1-2, to=2-2]
	\arrow[from=1-1, to=2-2]
\end{tikzcd}\]
As in \S6, via this map, we can pullback the group scheme inclusions on $Z_{_\xi}$ \eqref{hbt2} to $\spec D$. Consequently, one also has the corresponding Lie algebras of these group schemes. Via $\phi_0$, one can pullback the data \eqref{dataonrxi} on $\spec R_{\xi}$ to $\spec D_0$, and obtain a diagram
% https:\parallelslantq.uiver.app/?q=WzAsMyxbMCwwLCJEXzBcXHtYXzEsXFxsZG90cyxYX2dcXH0iXSxbMiwwLCJcXHNjcntBfV97RF8wfSJdLFsxLDEsIlxcbGllIFxcZnJha3tIfV97RF8wfSJdLFswLDFdLFswLDIsIiIsMix7InN0eWxlIjp7ImJvZHkiOnsibmFtZSI6ImRvdHRlZCJ9fX1dLFsyLDFdXQ==
\[\begin{tikzcd}
	{D_0\{X_1,\ldots,X_g\}} && {\ca_{_{D_0}}} \\
	& {\lie~\bh_{_{D_0}}}
	\arrow[from=1-1, to=1-3]
	\arrow[dotted, from=1-1, to=2-2]
	\arrow[from=2-2, to=1-3]
\end{tikzcd}\]
Now, $\lie ~\bh_{_{D_0}}$ arises just from the free module $\lie ~\bh_{_D}$ over $D$. Thus, if we have any sections $X_1,\ldots, X_g$ of $\lie {\bh_{_{D_0}}}$, we can lift them to sections of $\lie {\bh_{_D}}$, and thus obtain a diagram, now on $\spec D$: 
% https:\parallelslantq.uiver.app/?q=WzAsMyxbMCwwLCJEXFx7WF8xLFxcbGRvdHMsWF9nXFx9Il0sWzIsMCwiXFxzY3J7QX1fe0R9Il0sWzEsMSwiXFxsaWUgXFxmcmFre0h9X3tEfSJdLFswLDFdLFswLDIsIiIsMix7InN0eWxlIjp7ImJvZHkiOnsibmFtZSI6ImRvdHRlZCJ9fX1dLFsyLDFdXQ==
\[\begin{tikzcd}
	{D\{X_1,\ldots,X_g\}} && {\ca_{_D}} \\
	& {\lie{\bh_{_D}}}
	\arrow[from=1-1, to=1-3]
	\arrow[dotted, from=1-1, to=2-2]
	\arrow[from=2-2, to=1-3]
\end{tikzcd}\]
Observe that the map $D\{X_1,\ldots, X_g\}\to \ca_{_D}$ is surjective, as it was at the closed point (this is the statement for $D_0$). This in fact arises from the lift of the map $\spec D_0 \to Z_{_\xi} \to \overline{V}_{_n}$, to the map $\spec D\to \overline{V}_{_n}$. It remains to find a \textit{compatible} map $\spec D\to Z(H)$. 

Let $K_{_D}$ denote the field of fractions of $D$.  Then by Lemma \ref{faprime1}, 
we have a map $\lie {\bh_{_{{K_{_D}}}}}\to \lie {H_{_{{K_{_D}}}}}$.  Since $H_{_{{K_{_D}}}}$ is a semisimple group, by going to a finite Galois extension $L_{_D}/K_{_D}$, we have an isomorphism $\lie {H_{_{{L_{_D}}}}}\cong  \lie {H} \otimes L_{_D}$. Again, we may also assume that by going to this extension the Azumaya algebra $\ca$ ``splits" - i.e. $\ca\cong \lie \cu_{_{K_{_D}}}\otimes L_{_D}$. 

%for which the following are satisfied.
%\begin{itemize}
%\item We have a map $\lie {\bh_{_{{K_{_D}}}}}\to \lie {H_{_{{K_{_D}}}}}$ induced by the morphism of group schemes, such that the map over the generic point is a conjugate of the canonical one pulled back from $\spec K_{_D}$. 
%\item We have an isomorphism $\lie {H_{_{{L_{_D}}}}}\cong  \lie {H_{_{K_{_D}}}} \otimes L_{_D}$. 
%\item $\ca$ ``splits" - i.e. $\ca\cong \lie \cu_{_{K_{_D}}}\otimes L_{_D}$. 
%\end{itemize}
Thus, on $\spec L_{_D}$ we have a map $L_{_D}\{X_1,\ldots, X_g\}\to \lie {\bh_{_{{L_{_D}}}}} \to\lie {H} \otimes L_{_D} \to \lie \cu_{_{K_{_D}}}\otimes L_{_D}$, where the last two arrows are really maps of Lie algebras, the first one being the one from the free Lie algebra on $g$ generators. Observe that the composite $L_{_D}\{X_1,\ldots, X_g\}\to \lie {\bh_{_{{L_{_D}}}}} \to\lie {H} \otimes L_{_D}$ gives a map $\spec L_{_D} \to (\lie {H})^{^g}$. Composition by the canonical projection from $(\lie {H})^{^g} \to Z(H)$, gives a map $\psi:\spec L_{_D} \to Z(H)$. Now, consider the action of $\Gal L_{_D}/K_{_D}$ on $\psi$. The action comprises two ``parts". It first acts on the arrow $L_{_D}\{X_1,\ldots, X_g\}\to \lie {\bh_{_{{L_{_D}}}}} $, which really comes as a pull-back and hence descends to $\spec K_{_D}$. For the latter part, we proceed as follows. 

Consider the diagram

% https://q.uiver.app/?q=WzAsNixbMywxLCJaX24oXFxQR0wpIl0sWzMsMCwiTV9uXmciXSxbMiwwLCIoXFxsaWUgSCleZyJdLFsyLDEsIlooSCkiXSxbMCwxLCJcXHNwZWMgUyciXSxbMCwyLCJcXHNwZWMgUyJdLFsyLDEsIiIsMCx7InN0eWxlIjp7InRhaWwiOnsibmFtZSI6Imhvb2siLCJzaWRlIjoidG9wIn19fV0sWzEsMCwiXFxQR0wiXSxbMiwzLCJcXG92ZXJsaW5le0h9IiwyXSxbMywwXSxbNCw1LCJcXEdhbCBMX0QvS19EIiwyXSxbNCwyXSxbNSwwXV0=

\[\begin{tikzcd}
	&& {(\lie H)^g} & {\cm_{_n}^g} \\
	{\spec L_{_D}} && {Z(H)} & {Z_{_n}(\PGL)} \\
	{\spec K_{_D}}
	\arrow[hook, from=1-3, to=1-4]
	\arrow["\PGL", from=1-4, to=2-4]
	\arrow["{\overline{H}}"', from=1-3, to=2-3]
	\arrow[from=2-3, to=2-4]
	\arrow["{\Gal L_{_D}/K_{_D}}"', from=2-1, to=3-1]
	\arrow[from=2-1, to=1-3]
	\arrow[from=3-1, to=2-4]
\end{tikzcd}\]
We replace $(\lie H)^g$ by the pullback of $\cm_{_n}^g$ to $Z(H)$. Let $Z'$ denote this pullback. Then $Z(H)$ is a quotient of $Z'$ by $\PGL$. Thus we have:
\beqa
\begin{tikzcd}
	{Z'} & {\cm_{_n}^g} \\
	{Z(H)} & {Z_{_n}(\PGL)}
	\arrow["\PGL", from=1-2, to=2-2]
	\arrow[from=1-1, to=1-2]
	\arrow["\PGL"', from=1-1, to=2-1]
	\arrow[ from=2-1, to=2-2]
\end{tikzcd}
\eeqa

 Since it is clear from the above diagram that the action of $\Gal L_{_D}/K_{_D}$ on $\spec L_{_D}$ is compatible with that of $\PGL$ on $\cm_{_n}^g$, we see that the  map $\spec L_{_D} \to Z'$ descends to a map $\spec K_{_D} \to Z(H)$.

Now let $\spec D_{_\mathfrak p}$ denote the local ring of a height $1$ prime $\mathfrak p$ in $\spec D$. By the discussion above we have the map $K_{_D}\{X_1,\ldots, X_g\}\to \lie \bh_{_{K_{_D}}}$. 

What proceeds may be interpreted in the background of the  ``semistable reduction theorem" for principal bundles (see \cite[Proposition 8]{bapa}).  As in Lemma \ref{faprime1},  one can extend the map over the height $1$ prime $\mathfrak p$ to take its values in $(\lie H)^g$, and again, as above, the map descends to a map $\spec D_{_\mathfrak p} \to Z(H)$. Consequently, we have a {\em rational map} $\spec D \rightarrow Z(H)$ which is defined over every height $1$ prime $\mathfrak p$. By normality and affineness of $\spec D$ and $Z(H)$, using Hartogs' lemma, we can extend this rational map to a morphism $\phi:\spec D\to Z(H)$ as required. 

To check that this gives us a lift $\phi:\spec D \to Z_{_\xi}$, we need to check that the composites, $\spec D\to Z(H)\to Z_{_n}$ and $\spec D\to \overline{V}{_{_n}}\to Z_{_n}$ agree. For this, it suffices to check that the further composites into $\Delta_{_n}:Z_{_n} \hookrightarrow Z_{_{n^2}}$ agree. This follows from the fact that the manner in which the map is constructed from $\spec K_{_D}$ to $Z(H)$ is compatible with the composite into $Z_{_{n^2}}$, and that moreover, any two maps which coincide over a dense open set (that being in our case $\spec K_{_D}$) from a reduced scheme to a separated scheme must in fact be the same.

What we have at this point is that the morphism $\pi$ \eqref{thedsingmaphere} is birational over  $Z^s_{\tsf{f}}(H)$. The last bit which remains to be checked is that this morphism is birational over the entire locus $Z(H)^s$. 

Let a point as in \eqref{pointsofzh} lie in $Z(H)^s$ in the complement of $Z^s_{\tsf{f}}(H)$.
The strategy is again as in the proof in \S\ref{secondpart}. We arrive at the situation in the last paragraph  after \eqref{locnontr},  and have the closed fibre $H'$ of a maximal parahoric group scheme mapping to $H$ and factoring via the Levi quotient $L_{_o}'$.
We {\em claim that} the element $\beta \in Z(L'_{_o})$ there now provides a non-trivial automorphism of the point in $Z(H)^s$. The data we have is as follows:
\[\begin{tikzcd}
	& {\lie{H'}} \\
	{k\{X_1, \ldots, X_g\}} & {\lie{L_{_o}}'} \\
	& {\lie{H}}
	\arrow[from=2-1, to=1-2]
	\arrow[from=1-2, to=2-2]
	\arrow[from=2-2, to=3-2]
	\arrow[from=2-1, to=2-2]
	\arrow[from=2-1, to=3-2]
\end{tikzcd}\]
Since the morphism $k\{X_1, \ldots, X_g\} \to {\lie{H}}$ comes via ${\lie{L_{_o}}'}$ and since $\beta \in Z(L'_{_o})$, it is easy to see that we get a non-trivial element in the isotropy subgroup at the point in ${\lie{H}}^g$ corresponding to the arrow $k\{X_1, \ldots, X_g\} \to {\lie{H}}$. This contradicts the assumption that the action of $\overline{H}$ on points of $({\lie{H}})^{g,s}$ is free.

%Now, let $S'$ be the integral closure of a local ring of a height $1$ prime in a finite extension as before. Then by the previous remarks, we see that the maps to $Z_{_{n^2}}$ induced by the inclusions $\lie \bh_{S'}\to \mathfrak{gl}_{n^2}$ and $\lie \bh_{S'}\to \mathfrak{l}\to \lie \bh_{S'}\to \mathfrak{gl}_{n^2}$ are in the same orbit closure. Combined with the fact that the images of $\lie \bh\to \mathfrak{l} \to \lie H \to \mathfrak{gl}_{n^2}$ and $\lie \bh_{S'}\to \mathfrak{l}\to \lie \bh_{S'}\to \mathfrak{gl}_{n^2}$ are $\GL_{n^2}$-conjugate, we see that the maps to the quotient $Z_n\parallelslant\PGL$ coincide, which finishes the proof. 

\end{document}